\documentclass[journal,onecolumn]{IEEEtran}
\usepackage{epsfig,color,amsmath,cite}

\IEEEoverridecommandlockouts
\usepackage{wrapfig}
\usepackage{setspace}  
\usepackage{xcolor}
\usepackage{epstopdf}
\usepackage{amssymb}
\usepackage{url}
\usepackage{mathtools}
\usepackage{enumitem}
\usepackage{multirow}
\usepackage{makecell}
\usepackage{amsmath}
\usepackage{stackrel}
\usepackage{booktabs}
\usepackage[flushleft]{threeparttable}
\usepackage{makecell}

\usepackage{subfig}
\usepackage{graphicx}

\usepackage[linesnumbered,lined,boxed,commentsnumbered,ruled,longend]{algorithm2e}


\setlist[itemize]{leftmargin=*}

\newcommand{\kron}{\otimes}
\DeclareMathOperator{\diag}{diag}

\DeclareMathOperator*{\minimize}{minimize}

\DeclareMathOperator{\subjectto}{subject\ to}

\makeatother

\newcommand{\bmat}[1]{\begin{bmatrix} #1 \end{bmatrix}}

\def\blue{\textcolor{black}}

\newcommand{\m}{\boldsymbol}
\allowdisplaybreaks[4]
\usepackage[colorlinks = true,
linkcolor = blue,
urlcolor  = blue,
citecolor = blue,
anchorcolor = blue]{hyperref}

\newcommand{\mbb}[1]{\mathbb{#1}}

\usepackage{amsthm}
\newtheorem{mypro}{Property}
\newtheorem{myrem}{Remark}
\newtheorem{asmp}{Assumption}

\newtheorem{exmpl}{Example}
\usepackage[usestackEOL]{stackengine}
\stackMath
\usepackage[framemethod=TikZ]{mdframed}
\usepackage{bm}
\mdfdefinestyle{MyFrame}{%
	linecolor=black,
	outerlinewidth=1.5pt,
	roundcorner=1.5pt,
	innerrightmargin=5pt,
	innerleftmargin=5pt,}

\usepackage{tikz}
\usetikzlibrary{matrix,positioning,decorations.pathreplacing}

\DeclareMathOperator{\Tr}{Tr}

\usepackage[export]{adjustbox}
\begin{document}
	%\title{\vspace{0.2cm}\LARGE  State Space Formulation and Model Predictive Control for Water Quality Control  Problem }
\title{How Effective is Model Predictive Control in Real-Time Water Quality Regulation? State-Space Modeling and Scalable Control}
	\author{Shen Wan$\text{g}^\dagger$, Ahmad F. Tah$\text{a}^{\dagger,*}$, and Ahmed A. Abokif$\text{a}^\ddagger$ \thanks{
		$^*$Corresponding author. 	$^\dagger$Department of Electrical and Computer Engineering, The University of Texas at San Antonio, TX 78249.   $^{\ddagger}$Department of Civil and Materials Engineering, The University of Illinois at Chicago. Emails: mvy292@my.utsa.edu, abokifa@uic.edu, ahmad.taha@utsa.edu. This material is based upon work supported by the National Science Foundation under Grants 1728629, 2015671, and 2015603. }} 
	\maketitle
	
	%\begin{abstract}
	%Real-time water quality control (WQC) in water distribution networks (WDN), the problem of regulating disinfectant levels,  is challenging due to lack of \textit{(i)} a proper control-oriented modeling considering complicated components (junctions, reservoirs, tanks, pipes, pumps, and valves) for water quality modeling in WDN and \textit{(ii)} a corresponding scalable control algorithm dealing with all kind of uncertainty (i.e., demand and parametric uncertainty, unknown disturbances, and contamination events).  In this paper, we solve WQC problem by \textit{(a)} proposing a novel state-space representation of the WQC problem that provides explicit relationship between inputs (chlorine dosage at booster stations) and outputs (chlorine concentrations in the entire network) and \textit{(b)} designing a highly scalable model predictive control (MPC) algorithm that showcases fast response time and resilience against various sources of uncertainty. 
	%%The broader impact includes providing the control theory community with thorough dynamic modeling 
	%\end{abstract}
	
	\begin{abstract}
		\large Real-time water quality control (WQC) in water distribution networks (WDN), the problem of regulating disinfectant levels,  is challenging due to lack of \textit{(i)} a proper control-oriented modeling considering complicated components (junctions, reservoirs, tanks, pipes, pumps, and valves) for water quality modeling in WDN and \textit{(ii)} a corresponding scalable control algorithm that performs real-time water quality regulation.
In this paper, we solve the WQC problem by \textit{(a)} proposing a novel state-space representation of the WQC problem that provides an explicit relationship between inputs (chlorine dosage at booster stations) and states/outputs (chlorine concentrations in the entire network) and \textit{(b)} designing a highly scalable model predictive control (MPC) algorithm that showcases fast response time and resilience against some sources of uncertainty.
	\end{abstract}


\markboth{WATER RESOURCES RESEARCH, IN PRESS, DECEMBER 2020}{}

\section{\large Introduction and Paper Contributions}~\label{sec:literature}

\IEEEPARstart{B}{efore} drinking water leaves water treatment plants, chemical disinfection is typically applied to ensure the microbiological safety of the treated water. Water utilities worldwide rely on chlorine-based disinfectants due to their strong antimicrobial activity and low cost. Excess chlorine is usually applied at the treatment plant to prevent microbial re-contamination of the treated drinking water as it moves through the pipes of water distribution networks (WDN).

Residual chlorine concentrations are routinely monitored to verify that a sufficient residual is maintained throughout WDN. Maintenance of a detectable residual is also typically mandated by state and federal regulations in many countries. For instance, water utilities in the US are required to preserve detectable chlorine residual throughout their WDNs under the Surface Water Treatment Rule (SWTR) \cite{haas1999benefits}, and many states have established even more stringent numerical thresholds on  the minimum residual concentration \cite{roth2018dbp}. 

Nevertheless, determining the appropriate chlorine dosage to ensure a sufficient residual, particularly at the far ends of WDNs where the water age is the highest, is rather challenging. Applying large doses of chlorine-based disinfectants at the treatment plant has been associated with multiple issues, including the excessive formation of disinfection byproducts as well as aesthetic issues with water taste and odor \cite{fisher2011evaluation,hua2015variable}. Alternatively, the disinfectant can be injected in smaller doses at multiple locations in the network, a practice commonly known as booster disinfection, to maintain a uniform disinfectant concentration throughout the WDN \cite{tryby1999booster}. Most recently, to solve the problem of low disinfectant concentrations at critical dead-end nodes with no need of increasing disinfectant dose at sources or installing additional booster stations, the modulation of nodal outflows in WDN is proposed~\cite{avvedimento2020modulating}. For more context of real-time control of water quality in WDN; see~\cite{creaco2019real}.

%The objective of this paper 

\subsection{Literature review}
Over the past two decades, many studies have investigated the water quality control problem (WQC) of optimizing the locations and/or dosing schedules of booster disinfection systems. 

A wide range of optimization-based methods was used to solve the WQC problem, including linear programming (LP), quadratic programming (QP), heuristic algorithms such as genetic algorithm (GA), and multi-objective optimization algorithms. Boccelli \textit{et al.}~\cite{boccelli1998optimal} solved the WQC problem to minimize the total disinfectant mass dose required to satisfy residual constraints. Applying the principle of linear superposition to disinfectant concentrations resulting from multiple injections, they formulated the WQC problem as an LP. Their work was later extended by~\cite{tryby2002facility} by including the locations of booster stations as decision variables resulting in a mixed-integer linear programming (MILP) problem. WQC problem was also solved by~\cite{propato2004linear} who proposed a linear least-squares method, which is a QP, to optimize disinfectant injection rates. Munavalli and Kumar~\cite{munavalli2003optimal} determined the injection rates of boosters whose locations are known in advance, and used a GA algorithm to solve the optimization problem. 

GA was also implemented by Ostfeld and Salomons \cite{ostfeld2006conjunctive} to simultaneously optimize pump scheduling, the layout, and operation of booster chlorination stations, and was combined with multi-species water quality simulations to incorporate {disinfection by-product (DBP)} levels in the constraints \cite{ohar2014optimal}. Furthermore, multi-objective optimization was also applied to solve the WQC problem, where Prasad \textit{et al.} \cite{prasad2004booster} used a multi-objective genetic algorithm (NSGA-II) to minimize the total disinfectant dose and simultaneously maximize the water demand within specified residual limits. A comprehensive review of the literature on the optimization of booster chlorination systems can be found in the recent works by Islam \textit{et al.} \cite{islam2017optimizing}, and Mala-Jetmarova \textit{et al.} \cite{mala2017lost}. 

From these studies, we summarize the  WQC problem as being comprised of two main components: \textit{water quality modeling} and \textit{water quality control}. Water quality modeling depicts  the decay and transport of chlorine-based disinfectants in the WDN. This can be expressed via the advection-reaction dynamics for which three different families of numerical schemes can be used to obtain the numerical solution:  Eulerian-based schemes~\cite{rossman1993discrete}, Lagrangian-based schemes~\cite{liou1987modeling,boulos1994event}, and hybrid Eulerian–Lagrangian schemes~\cite{basha2007eulerian}. {For example, EPANET~\cite{rossman1994modeling} is a widely used modeling software using a Lagrangian-based approach.}  Water quality control is the corresponding algorithm or mechanism applied control to water quality modeling to reach the control objectives. In theory, WQC can be performed through a plethora of control algorithms with different underlying concepts including feedback control, adaptive control, and model predictive control (MPC).

We observe that a common drawback exists in the majority of the aforementioned studies: the utilized approaches failed to write or formulate water quality modeling explicitly depicting the relationship between multiple network inputs (booster injection) and outputs (various critical junction concentrations). This implies that the majority of water quality modelings are not designed for control-theoretic algorithms to be faithfully applied---this hinders applying state-of-the-art control algorithms for water quality. {To overcome some aforementioned drawbacks, Zierolf \textit{et al.}~\cite{zierolf1998development} derived an input/output (I/O) water quality model giving  explicit relationship between inputs and outputs.}  However, this model does not consider storage components. This I/O model was then further extended by~\cite{shang2000input,shang2002particle} to allow storage tanks, multiple water sources, and quality inputs.  {The I/O model is friendly for control purposes since it can be applied in the WQC problem after combining with one of the aforementioned control algorithms including  adaptive control~\cite{polycarpou2002feedback,wang2005adaptive}, and MPC~\cite{chang2003robust,duzinkiewicz2005hierarchical}.} 

Besides the I/O model which creates relationships between specific inputs and outputs in a general dynamic system, another standard modeling paradigm is a state-space model, i.e., a mathematical model of a physical system as a set of input $\m u$, output $\m y$, and state variable $\m x$ related by first-order differential equations (in continuous time) or difference equations (in discrete time). In contrast with the rather simplistic I/O model, control-theoretic, state-space models also capture the evolution of all physical states in the system. This can either be too complex or cumbersome to incorporate within I/O models that only capture specific input/output relationships. {Additionally, numerous I/O models need to be developed based on the different combinations of locations of inputs (boosters) and outputs (sensors), but one space-state model can cover such a situation.} For linear dynamic systems, state-space representation can be written as
%\begin{linenomath*}
	\begin{equation}~\label{equ:general}
	\begin{aligned}
	{\m{x}}(t+1) &=\m A(t) \m{x}(t)+\m B(t) \m {u}(t), \\
	 \m{y}(t) &=\m C(t) \m{x}(t)+\m D(t) \m{u}(t)
	\end{aligned}
	\end{equation}
%\end{linenomath*}
where $\m x \in \mbb{R}^{n_x}$ is the vector of system states (i.e., concentrations in all network components and pipe segments), and we do not consider hydraulic states; $\m u \in \mbb{R}^{n_u}$ includes all control inputs (i.e., chlorine injection schedules); $\m y \in \mbb{R}^{n_y} $ models sensor measurements (i.e., concentrations at {selected network locations});  $\m A$, $\m B$, $\m C$, and $\m D$ are corresponding state-space matrices modeling the evolution of the states according to water chemistry and conservation of mass in space and time. 

The model~\eqref{equ:general} is absent from the WQC literature; all WQC studies do not report or propose thorough state-space dynamics akin to~\eqref{equ:general}, and all control studies~\cite{polycarpou2002feedback,wang2005adaptive,chang2003robust,duzinkiewicz2005hierarchical} built an I/O model that does not explicitly model the states, which results in a model that only captures the output performance rather than all state variables $\m x$ in the network.  In particular, Zierolf \textit{et al.}~\cite{zierolf1998development} claimed that 
\begin{quotation}
	\textit{Formulating the problem of chlorine concentration with a state-space model is rather intractable due to the dimension of the system \cite{zierolf1998development}.}
\end{quotation}
In this paper, we overcome this challenge by showcasing that the formulation of the state-space model is indeed possible and tractable, while demonstrating that optimal model predictive WQC can also be scalable.

\subsection{Paper objectives and contributions}

To that end, the objective of this study is to build a time-varying, control- and network-theoretic model that corresponds to the variation of chlorine concentrations in arbitrary water networks with all kinds of components. This model is then explicitly used to inform the near real-time, optimal injections of booster chlorination stations, while satisfying and optimizing water quality constraints and metrics. Rather than relying on simulation-based techniques, we propose a novel approach for solving the WQC problem by coupling an MPC algorithm with a scalable control-oriented water quality modeling approach. The designed model and optimal control formulations are informed by the network structure and the future evolution of chlorine concentrations.  The paper contributions are:
\begin{itemize}
	\item A novel state-space formulation of WQC problem is developed, akin to~\eqref{equ:general}. {In particular, this paper presents the very first thorough attempt to build a control-theoretic model for water quality simulation---without needing a water quality simulation toolbox.} This model is indeed large-scale with millions of state variables for mid-size networks, endowed with great detail and accuracy, its derivation is provided with thorough examples {in appendices}, and numerical case studies show that the derived state-space model produces nearly identical water quality simulations in comparison with simulation packages.  
	\item The derived state-space model is then used to formulate a real-time WQC problem as a constrained quadratic problem (QP), via a model predictive control (MPC) approach---a contemporary approach to solving control problems. Due to the large dimensionality of the formulated QP, and to allow for real-time implementation,  a simple transformation is utilized to show how analytical solutions of the QP can be obtained without needing to solve a large-scale optimization problem. In short, the proposed control algorithm is truly plug-and-play: the system operator measures the chlorine concentrations at certain junctions and the control algorithm \textit{immediately} provides optimal dosages of chlorine where booster stations are installed in a closed-loop, feedback control fashion.
	\item Thorough case studies show the potential of applying the MPC-based, WQC on water networks with some sources of uncertainty (i.e., demand uncertainty, reaction rate coefficients uncertainty, and unmodeled disturbances impacting chlorine concentrations). 
\end{itemize}
The rest of the paper is organized as follows. Section~\ref{sec:Ctrl-WQM} introduces control-oriented water quality modeling by presenting the models of each component in detail.  An abstract, nonlinear, state-space format for the water quality model is given first, then the corresponding linear form is presented after considering the first-order reaction model in Section~\ref{sec:abstract-WQM}.   The control algorithm, specifically, MPC, starts to be introduced in Section~\ref{sec:LP-WQC}, and a linear program is formed. {Section~\ref{sec:scalable-WQC} presents a scalable algorithm that transforms the linear program with millions of variables to a quadratic program with fewer variables.} Section~\ref{sec:test} presents case studies to corroborate the paper's findings, and the limitations and future research directions are given in Section~\ref{sec:limitations}. 
All mathematical proofs are given in the appendices, which also contain other important derivations.   The notation for this paper is introduced next.

%\vspace{1em}
\noindent \textit{\textbf{Paper's Notation}} --- Italicized, boldface upper and lower case characters represent matrices and column vectors: $a$ is a scalar, $\m a$ is a vector, and $\m A$ is a matrix. Matrix $\m I$ denotes a identity square matrix, whereas $\m O_{m \times n}$ denotes a zero matrix  with size $m$-by-$n$.
The notations $\mathbb{R}$ and $\mathbb{R}_{++}$ denote the set of  real and positive real numbers.  The notations $\mathbb{R}^n$ and $\mathbb{R}^{m\times n}$ denote a column vector with $n$ elements and an $m$-by-$n$ matrix in $\mathbb{R}$. For any vector $\m x \in \mathbb{R}^{n}$,  $\m x^{\top}$ is its  transpose. 
For $\m x, \m y \in\mbb{R}^n$, the element-wise product or Hadamard product is defined as $\m x \circ \m y := [x_1y_1,x_2y_2,\ldots,x_ny_n]^\top$, which can also be expressed by  matrix multiplication form, that is $\m x \circ \m y := \diag(\m x) \m y = \diag(\m y) \m x$;  similarly,  the element-wise division or Hadamard division  is defined as $\m x \oslash \m y := [\frac{x_1}{y_1}, \frac{x_2}{y_2},\ldots, \frac{x_n}{y_n}  ]^\top$ when $\m x\in\mbb{R}^n, \m y\in\mbb{R}_{++}^n$, and   $\m x \oslash \m y  = \diag(\m y)^{-1} \m x$. For any two matrices $\m A$ and $\m B$ with same number of columns, the notation $\{\m A, \m B\}$ denotes $[\m A^\top \  \m B^\top]^\top$ and $\m A \kron \m B$ stands for the Kronecker product of $\m A$ and $\m B$. The trace of square matrix $\m A$ is $\Tr(\m A)$.
\section{\large Control-Oriented Water Quality Modeling}~\label{sec:Ctrl-WQM}
We model WDNs by a directed graph $\mathcal{G} = (\mathcal{W},\mathcal{L})$.  Set $\mathcal{W}$ defines the nodes and is partitioned as $\mathcal{W} = \mathcal{J} \bigcup \mathcal{T} \bigcup \mathcal{R}$ where $\mathcal{J}$, $\mathcal{T}$, and $\mathcal{R}$ are collection of junctions, tanks, and reservoirs. Let $\mathcal{L} \subseteq \mathcal{W} \times \mathcal{W}$ be the set of links, and define the partition $\mathcal{L} = \mathcal{P} \bigcup \mathcal{M} \bigcup \mathcal{V}$, where $\mathcal{P}$, $\mathcal{M}$, and $\mathcal{V}$ represent the collection of pipes, pumps, and valves. Let downstream nodes of pumps and valves as a set $\mathcal{D}$. Hence, the nodes having no connection with pumps and valves and the upstream nodes of pumps and valves are  $\mathcal{W} \char`\\ \mathcal{D}$. For the $i$-th node, set $\mathcal{N}_i$ collects its neighboring nodes (any two nodes connected by a link) and is partitioned as $\mathcal{N}_i = \mathcal{N}_i^\mathrm{in} \bigcup \mathcal{N}_i^\mathrm{out}$, where $\mathcal{N}_i^\mathrm{in}$ and $\mathcal{N}_i^\mathrm{out}$ are collection of inflow and outflow nodes. The number of junctions, reservoirs, tanks, pipes, pumps, and valves is $n_{\mathrm{J}}$, $n_{\mathrm{R}}$, $n_{\mathrm{TK}}$, $n_{\mathrm{P}}$, $n_{\mathrm{M}}$, and $n_{\mathrm{V}}$. In this paper, we use the Lax Wendroff~\cite{lax1964difference} numerical scheme to solve the advection-reaction dynamics in pipes, and Pipe $ij$ with length $L_{ij}$ is split into $s_{L_{ij}}$ segments. {The total number of segment in all pipes are $n_{\mathrm{S}} = \sum\limits_{{ij} \in \mathcal{P}} s_{L_{ij}}$. Thus, the number of chlorine concentration variables at nodes and in links are $n_\mathrm{N} = n_{\mathrm{J}}+n_{\mathrm{R}}+n_{\mathrm{TK}}$ and $n_\mathrm{L} = n_{\mathrm{S}}  +n_{\mathrm{M}}+n_{\mathrm{V}}$. }

Before introducing water quality modeling for each WDN component, we introduce symbols and vector notations. We define $\m E^\mathcal{G}$ as the connectivity matrix of graph $\mathcal{G}$.  Notice that the assignment of direction to each link (and the resulting inflow/outflow node classification) is arbitrary.  Thus, $\m E^\mathcal{G}$ is comprised of $-1$, $0$, and $1$ representing negative connection, no connection, and positive connection. Besides, $\m E_\mathcal{G}$ has row vector form $[{\m E^\mathrm{N}_{\mathrm{P}}}^\top\  {\m E^\mathrm{N}_{\mathrm{M}}}^\top\  {\m E^\mathrm{N}_{\mathrm{V}}}^\top]$  and column vector form $\{\m E^\mathrm{L}_{\mathrm{J}},\; \m E^\mathrm{L}_{\mathrm{R}},\; \m E^\mathrm{L}_{\mathrm{TK}} \}$  from  different perspectives, and
\setlength\extrarowheight{3pt}
%\begin{linenomath*}
	\begin{align}~\label{equ:WFP-Incidence-matrix}
	{\large \m E^\mathcal{G}} = \begin{array}{*{4}{cccc}@{}c}
	&\textit{Pipe}&\textit{Pump}&\textit{Valve}\\
	\cline{2-4}
	\textit{Junction}&\multicolumn{1}{|c|}{\m E^{\mathrm{P}}_\mathrm{J}} &\multicolumn{1}{|c|}{\m E^{\mathrm{M}}_\mathrm{J}}&\multicolumn{1}{|c|}{{\m E^{\mathrm{V}}_\mathrm{J}}}& \multirow{1}{*}{$\left.\rule[1ex]{0pt}{1ex}\right\}\m E^\mathrm{L}_\mathrm{J}$}\\
	\cline{2-4}
	\textit{Reservoir}&\multicolumn{1}{|c|}{\m E^{\mathrm{P}}_\mathrm{R}}  &\multicolumn{1}{|c|}{{\m E^{\mathrm{M}}_\mathrm{R}}}&\multicolumn{1}{|c|}{{\m E^{\mathrm{L}}_\mathrm{R}}}& \multirow{1}{*}{$\left.\rule[1ex]{0pt}{1ex}\right\}\m E^\mathrm{L}_\mathrm{R}$}\\
	\cline{2-4}
	\textit{Tank}&\multicolumn{1}{|c|}{\m E^{\mathrm{P}}_\mathrm{TK}}  &\multicolumn{1}{|c|}{{\m E^{\mathrm{M}}_\mathrm{TK}}}&\multicolumn{1}{|c|}{{\m E^{\mathrm{L}}_\mathrm{TK}}}& \multirow{1}{*}{$\left.\rule[1ex]{0pt}{1ex}\right\}\m E^\mathrm{L}_\mathrm{TK}$}\\
	\cline{2-4}
	\noalign{\vspace{-11pt}}
	\multicolumn{1}{c}{}	& \multicolumn{1}{@{}c@{}}{\underbrace{\hspace*{2\tabcolsep}\hphantom{.......}}_{\textstyle {\m E^{\mathrm{N}}_\mathrm{P}}^\top}}  & \multicolumn{1}{@{}c@{}}{\underbrace{\hspace*{2\tabcolsep}\hphantom{.......}}_{\textstyle {\m E^{\mathrm{N}}_\mathrm{M}}^\top}} & \multicolumn{1}{@{}c@{}}{\underbrace{\hspace*{2\tabcolsep}\hphantom{.......}}_{\textstyle \ {\m E^{\mathrm{N}}_\mathrm{V}}^\top}}
	\end{array}
	\end{align}
%\end{linenomath*}
\noindent where $\m E^\mathrm{N}_{\mathrm{M}}$ represents a connectivity matrix from nodes ($\mathrm{N}$) to pumps ($\mathrm{M}$), and it lumps connectivity submatrices $\m E_\mathrm{J}^{\mathrm{M}}$, $\m E_\mathrm{R}^{\mathrm{M}}$, and $\m E_\mathrm{TK}^{\mathrm{M}}$ from junctions, reservoirs, and tanks to pumps. Moreover, ${\m E_\mathrm{J}^{\mathrm{M}}}^\top = \m E^\mathrm{J}_{\mathrm{M}}$, and the other symbols in~\eqref{equ:WFP-Incidence-matrix} have similar meaning. Similarly, we define the connectivity matrix from booster stations to nodes as
%\begin{linenomath*}
	{ 	\setlength\extrarowheight{3pt}
		\begin{align} ~\label{equ:WFP-Booster-matrix}
		{\large \m E^\mathrm{B}_\mathrm{N}} =\begin{array}{*{4}{cccc}@{}c}
		& \multicolumn{3}{c}{Booster}\\
		\cline{2-4}
		\textit{Junction}&\multicolumn{1}{|c|}{\m E^{\mathrm{b}}_\mathrm{J}} &\multicolumn{1}{|c|}{\m O}&\multicolumn{1}{|c|}{\m O}& \multirow{1}{*}{$\left.\rule[0ex]{-3pt}{0ex}\right\}\m E^\mathrm{B}_\mathrm{J}$}\\
		\cline{2-4}
		\textit{Reservoir}&\multicolumn{1}{|c|}{\m O}  &\multicolumn{1}{|c|}{\m E^{\mathrm{b}}_\mathrm{R}}&\multicolumn{1}{|c|}{\m O}& \multirow{1}{*}{$\left.\rule[0ex]{-3pt}{0ex}\right\}\m E^\mathrm{B}_\mathrm{R}$}\\
		\cline{2-4}
		\textit{Tank}&\multicolumn{1}{|c|}{\m O}  &\multicolumn{1}{|c|}{\m O}&\multicolumn{1}{|c|}{\m E^{\mathrm{b}}_\mathrm{TK}}& \multirow{1}{*}{$\left.\rule[0ex]{3pt}{0ex}\right\}\m E^\mathrm{B}_\mathrm{TK}$}\\
		\cline{2-4}
		\end{array}
		\end{align}
	}
%\end{linenomath*}

Note that \textit{(i)} boosters only installed at nodes (junctions, reservoirs, and tanks), \textit{(ii)} ${\m E^{\mathrm{b}}_\mathrm{J}}$, ${\m E^{\mathrm{b}}_\mathrm{R}}$, and  ${\m E^{\mathrm{b}}_\mathrm{TK}}$  are squared and diagonal submatrices composing of only $0$ and $1$, and it describes if a booster is installed or not at a specific location, \textit{(iii)} the summation of each row or column in each submatrix should be no larger than 1 to ensure each component only has at most one booster installed, \textit{(iv)} the  number of boosters at each component is $n_\mathrm{bJ} = \Tr{(\m E^\mathrm{b}_\mathrm{J})}$, $n_\mathrm{bR} = \Tr{(\m E^\mathrm{b}_\mathrm{R})}$, $n_\mathrm{bTK} = \Tr{(\m E^\mathrm{b}_\mathrm{TK})}$, and the total number of boosters is $n_\mathrm{B} = \Tr{(\m E^\mathrm{B}_\mathrm{N})} = n_\mathrm{bJ} + n_\mathrm{bR} + n_\mathrm{bTK}$.

All symbols and vector notations are defined in Tab.~\ref{table:vector}. The quantity $c$ represents disinfectant concentration, and $c^\mathrm{P}(x,t)$ is the solute concentration in a pipe at location $x$ and time $t$. The superscript represents the quantity for specific components, i.e., vector $\m c^\mathrm{J}$ collects concentrations at all junctions, and vector $\m q^\mathrm{P}$ stands for the flow rates in all pipes. 
%In Tab.~\ref{table:vector}, the flow rate is in $\mathrm{GPM}$, concentrations and reaction rates are in $\mathrm{mg/L}$, and volumes of tanks are in $\mathrm{ft}^3$.

Next, we present the high-level objective of Sections~\ref{sec:Ctrl-WQM} and~\ref{sec:abstract-WQM}, that is, obtaining a model: $\m x(t+\Delta t) = \m f(\m x(t), \m u(t), \m {p}(t))$ where $\m x(t)$ includes all concentrations for all points in the network, $\m u(t)$ is booster station control action (chlorine dosage), $\m {p}(t)$ are the time-dependent parameters (such as flow rates, velocities); vector valued function $\m f(\cdot)$ mapping the dynamics in time and space has linear and nonlinear components depending on the reaction model. We remind readers that  water quality modeling of one component might depend on water quality modeling of another, and we would like to  discuss this further at the end of Section~\ref{sec:Ctrl-WQM}. After presenting water quality modeling of  components as difference equations, in Section~\ref{sec:abstract-WQM} we lump  them into state-space formulation akin to~\eqref{equ:general}.

\begin{table}
	%    \scriptsize
	\fontsize{10}{10}\selectfont
	\centering
	\setlength\tabcolsep{1.2 pt}
	\renewcommand{\arraystretch}{1.8}
	\makegapedcells
	\setcellgapes{1.4pt}
	\caption{Vector Notation of Water Quality Modeling. }
	\vspace{-0.2cm}
	\begin{tabular}{ c|c|c} 
		\hline
		\hspace{-4pt}\textit{Symbol}\hspace{-4pt} & \textit{Description} & \textit{Dimension}  \\ \hline
		$\m q^\mathrm{P}$, $\m q^\mathrm{M}$, $\m q^\mathrm{V}$& Flow rates in pipes, pumps, valves & $\mbb{R}^{n_{\mathrm{P}}  }$, $\mbb{R}^{n_{\mathrm{M} }}$, $\mbb{R}^{n_{\mathrm{V} }}$\\ \hline
		$\m q^\mathrm{L}$&\makecell{$\m q^\mathrm{L} \triangleq \{ \m q^\mathrm{P}, \m q^\mathrm{M}, \m q^\mathrm{V} \}$ \\collects flows in all links} & $\mbb{R}^{n_{\mathrm{Q}}}$\\ \hline
		$\m q^\mathrm{D}$ & \makecell{Demand (leaving junctions)} &$\mbb{R}^{n_{\mathrm{J} }}$\\ \hline
		$\m q_\mathrm{J}^\mathrm{B}$, $\m q_\mathrm{R}^\mathrm{B}$, $\m q_\mathrm{TK}^\mathrm{B}$&  \makecell{Booster flow rates injected to\\ junctions, reservoirs, and tanks}& $\mbb{R}^{n_{\mathrm{J}}}$, $\mbb{R}^{n_{\mathrm{R} }}$, $\mbb{R}^{n_{\mathrm{TK} }}$\\ \hline
		$\m q^\mathrm{B}$&\makecell{$\m q^\mathrm{B} \triangleq \{ \m q_\mathrm{J}^\mathrm{B}, \m q_\mathrm{R}^\mathrm{B},\m q_\mathrm{TK}^\mathrm{B} \}$ collects \\ booster flows rates injected to all nodes} & $\mbb{R}^{n_{\mathrm{N}}}$\\ \hline
		%$\m q^\mathrm{O}$, $\m m^\mathrm{O}$& \makecell{Net outflow rates (from connected links) leaving  \\ junctions or tanks and the corresonding mass rates.} \\ \hline
		$\m c^\mathrm{J}$, $\m c^\mathrm{TK}$, $\m c^\mathrm{R}$&  \makecell{Concentrations at junctions,\\ reservoirs, and tanks }& $\mbb{R}^{n_{\mathrm{J}}}$, $\mbb{R}^{n_{\mathrm{R} }}$, $\mbb{R}^{n_{\mathrm{TK} }}$\\ \hline
		$\m c^\mathrm{N}$&\makecell{$\m c^\mathrm{N} \triangleq \{\m c^\mathrm{J}, \m c^\mathrm{R}, \m c^\mathrm{TK} \}$ \\ collects concentrations at all nodes} & $\mbb{R}^{n_{\mathrm{N} }}$\\ \hline
		$\m c^\mathrm{P}$, $\m c^\mathrm{M}$, $\m c^\mathrm{V}$&  Concentrations in pipes, pumps, valves & {$\mbb{R}^{n_{\mathrm{S} }}$}, $\mbb{R}^{n_{\mathrm{M} }}$, $\mbb{R}^{n_{\mathrm{V} }}$\\ \hline  
		$\m c^\mathrm{L}$&\makecell{$\m c^\mathrm{L} \triangleq \{ \m c^\mathrm{P}, \m c^\mathrm{M}, \m c^\mathrm{V} \}$ \\ collects concentrations in  all links } & $\mbb{R}^{n_{\mathrm{L} }}$\\ \hline
		$\m c_\mathrm{J}^\mathrm{B}$, $\m c_\mathrm{R}^\mathrm{B}$, $\m c_\mathrm{TK}^\mathrm{B}$&  \makecell{Booster concentrations injected to\\ junctions, reservoirs, and tanks}& $\mbb{R}^{n_{\mathrm{J}}}$, $\mbb{R}^{n_{\mathrm{R} }}$, $\mbb{R}^{n_{\mathrm{TK} }}$\\ \hline
		$\m c^\mathrm{B}$&\makecell{$\m c^\mathrm{B} \triangleq \{ \m c_\mathrm{J}^\mathrm{B}, \m c_\mathrm{R}^\mathrm{B}, \m c_\mathrm{TK}^\mathrm{B} \}$ collects \\  booster concentrations at  all nodes } & $\mbb{R}^{n_{\mathrm{N} }}$\\ \hline
		$\m V^\mathrm{B}$&\makecell{$\m V^\mathrm{B} \triangleq \{ \m V_\mathrm{J}^\mathrm{B}, \m V_\mathrm{R}^\mathrm{B},\m V_\mathrm{TK}^\mathrm{B} \}$ collects \\ booster volumes injected to all nodes} & $\mbb{R}^{n_{\mathrm{N}}}$\\ \hline
		$\m V^\mathrm{TK}$&\makecell{ Volumes at all tanks }& $\mbb{R}^{n_{\mathrm{TK}}}$ \\ \hline
		$\m  r^\mathrm{P} $, $\m  r^\mathrm{TK} $ &  Reaction rate  for pipes, tanks &$\mbb{R}^{n_{\mathrm{P}}}$, $\mbb{R}^{n_{\mathrm{TK}}}$   \\ \hline 
		$ k^{b}$, $ k^{w}$ & Bulk and wall reaction rate constant & ---    \\ \hline \hline
		\multicolumn{3}{l}{\footnotesize{ \makecell{$^*$ Flow rate is in $\mathrm{GPM}$; concentration is in $\mathrm{mg/L}$; reaction rates are  \\in $\mathrm{mg\cdot L^{-1} \cdot h^{-1}}$; volumes of tanks are in $\mathrm{ft}^3$.}}}
	\end{tabular} 
	\label{table:vector}
	\vspace{-0.1cm}
\end{table}

\subsection{Conservation of mass}\label{sec:conservation}

The water quality modeling represents the movement of all chemical and/or microbial species (contaminant, disinfectants, DBPs, metals, etc.) within a WDN as they traverse various components of the network. This movement or time-evolution is based on three principles: \textit{(i)} mass balance in pipes, which is represented by solute transport in differential pipe lengths by advection in addition to its decay/growth due to reactions; \textit{(ii)} mass balance at junctions, which is represented by complete and instantaneous mixing of all in-flows; and \textit{(iii)} mass balance in tanks, which is represented by a continuously stirred tank reactors (CSTRs) model with complete and instantaneous mixing and  growth/decay reactions. {Note that simulation software package EPANET is also built on these three basic principles, which allows our method to be directly compared with it.}  The modeling of each component is introduced next.     
\begin{figure}[t]
	\centering
	\includegraphics[width=0.65\linewidth]{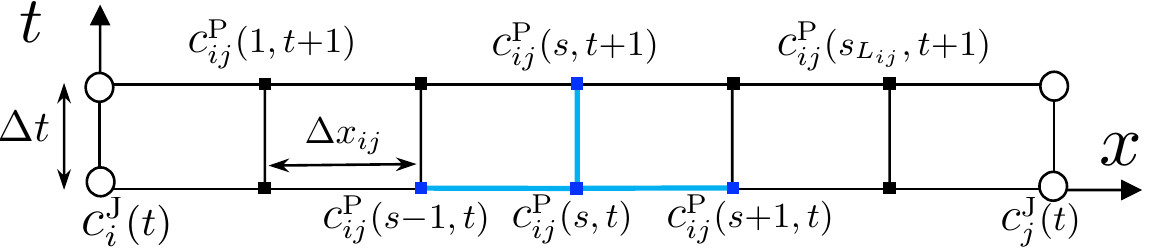}
	\caption{Time-space discretization of Pipe $ij$ based on the L-W scheme. }
	\label{fig:lax}
	%\vspace{-1.5em}
\end{figure}

\subsubsection{Solute transport and reaction in pipes}  The water quality modeling for pipes involves modeling the solute transport and reaction by  the 1-D advection-reaction (A-R) equation. For any Pipe $ij$, where $i$ and $j$ are indices of its upstream and downstream nodes, and the 1-D A-R equation is given by a PDE as
%\begin{linenomath*}
	\begin{equation} ~\label{equ:adv-reac}
	\textit{{Solute-Pipes}:} \; \; \; {\partial_t c^\mathrm{P}} =  -v_{ij}(t) {\partial_x c^\mathrm{P}} + r^\mathrm{P}_{ij} (c^\mathrm{P}) ,
	\end{equation}
%\end{linenomath*}
\noindent where $v_{ij}(t)$ is  flow velocity, that equals  flow rate ${q_{ij}(t)}$ divided by its cross-sectional area; and $r^\mathrm{P}_{ij}(c^\mathrm{P})=k_{ij}^\mathrm{P}(c^\mathrm{P})^n$ is the $n$-th order solute reaction rate, and $k_{ij}^\mathrm{P}$ is  the rate constant.

Here, the Lax-Wendroff (L-W) scheme~\cite{lax1964difference} shown in Fig.~\ref{fig:lax} is used to solve \eqref{equ:adv-reac}--this
model has been used and accepted in~\cite{rossman1996numerical,morais2012fast,fabrie2010quality}. Pipe $ij$ with length $L_{ij}$ is split into $s_{L_{ij}}$ segments and  the discretized form {of} any segment except the first and last one is given by~\eqref{equ:adv-reac-lax}  in Tab.~\ref{table:wqp}, and the coefficients for previous, current, and next segment  are
%\begin{linenomath*}
	\begin{subequations}~\label{equ:alphas}
		\begin{align}
		\underline{\alpha}_{ij}(t) &= 0.5\alpha_{ij}(t)\left(1+\tilde{\alpha}_{ij}(t)\right), \\
		\alpha_{ij}(t) &=  1-\tilde{\alpha}_{ij}^2(t),  \\
		\overline{\alpha}_{ij}(t) &= -0.5\alpha_{ij}(t)\left(1-\tilde{\alpha}_{ij}(t)\right),
		\end{align}
	\end{subequations}
%\end{linenomath*}
\noindent where $\tilde{\alpha}_{ij}(t)= {v_{ij}(t) \Delta t}(\Delta x_{ij})^{-1}$,  and $\Delta t$ and $\Delta x_{ij}$ are the time-step and the space-discretization step in Fig.~\ref{fig:lax}. The stability condition of  L-W scheme is $\tilde{\alpha}_{ij}(t) \in \left(0,1\right]$. Hence, the water quality time-step is constrained by $\Delta t \leq {\Delta x_{ij}}/|v_{ij}(t)|$. 
\begin{myrem}~\label{rem:Deltat}
	{To make the space discretization in L-W scheme  stable and accurate, we calculate $\Delta t_{ij} = {\Delta x_{ij}}/|{v_{ij}(t)|}$, $\forall ij \in \mathcal{P}$, and find the minimum value as the final $\Delta t$, that is, $\Delta t = \min(\Delta t_{ij})$.}
\end{myrem}

%%\begin{asmp}\label{asmp:leaving}
%	The concentration leaving a node, either by user demands or in  downstream links (a pump, a pipe, or a valve), is the same. 
%\end{asmp}
%Assumption~\ref{asmp:leaving} is reasonable and intuitive. That is, for Pipe $ij$ with a upstream Junction $i$, $c^\mathrm{J}_{i}(t)=c^\mathrm{P}_{ij}(1,t)$.  Hence, 

The first segment $c^{\mathrm{P}}_{ij}(1,t+\Delta t)$ and last segment $c^{\mathrm{P}}_{ij}(s_{L_{ij}},t+\Delta t)$, which are the special cases of~\eqref{equ:adv-reac-lax},  can be written as 
%\begin{linenomath*}
	\begin{subequations}~\label{equ:node2pipe}
		\begin{align}
		\hspace{-1.1em}	\begin{split} 	c^{\mathrm{P}}_{ij}(1,t+\Delta t) &= \underline{\alpha}_{ij}(t) c^{\mathrm{J}}_{i}(t)     +{\alpha}_{ij}(t) 	c^{\mathrm{P}}_{ij}(1,t) + \overline{\alpha}_{ij}(t) c^{\mathrm{P}}_{ij}(2,t) + r^\mathrm{P}_{ij}(c^{\mathrm{P}}_{ij}(1,t)), \label{equ:node2pipeA} 
		\end{split} \\
		\hspace{-1.1em}	\begin{split} 
		c^{\mathrm{P}}_{ij}(s_{L_{ij}},t+\Delta t)  &= \underline{\alpha}_{ij}(t) c^{\mathrm{P}}_{ij}(s_{L_{ij}}\hspace{-3pt}-\hspace{-3pt}1,t)   +{\alpha}_{ij}(t) 	c^{\mathrm{P}}_{ij}(s_{L_{ij}},t) + \overline{\alpha}_{ij}(t) c^{\mathrm{J}}_{j}(t) + r^\mathrm{P}_{ij}(c^{\mathrm{P}}_{ij}(s_{L_{ij}},t)). \label{equ:node2pipeB}
		\end{split}
		\end{align}
	\end{subequations}
%\end{linenomath*}
\normalcolor

\begin{figure}[t]
	\centering
	\includegraphics[width=0.5\linewidth]{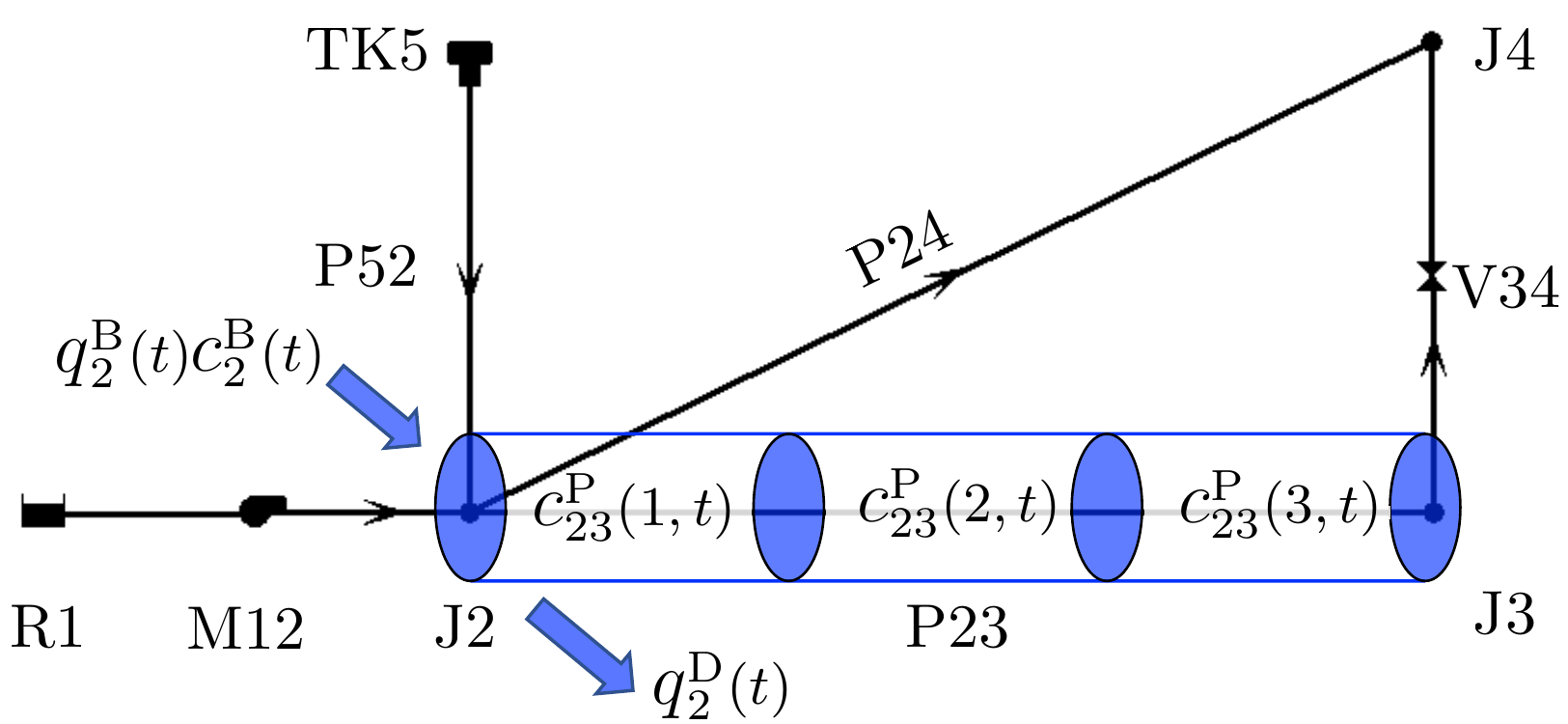}
	\caption{An illustrative example with pipes are divided into three segments (only $\mathrm{P}23$ is split and displayed).}
	\label{fig:example3}
\end{figure}   

{Note that Equation~\eqref{equ:node2pipe} is different from~\eqref{equ:adv-reac-lax} because we use the connected junction as the previous (next) segment of the first (last) segment  of a pipe.} After concentration of all segments are derived, we can lump them into $ \m c^{\mathrm{P}}_{ij}(t +\Delta t)$. {To help readers understand the details of~\eqref{equ:adv-reac-lax} and~\eqref{equ:node2pipe}, a pipe split into  three segments in Fig.~\ref{fig:example3} is shown as Example~\ref{exmp:pipe} in Appendix~\ref{app:example}.}   

%\textit{(i)} the first segment $c^{\mathrm{P}}_{ij}(1,t+\Delta t)$ depends on the previous time-step $c^{\mathrm{J}}_{i}(t)$, $c^{\mathrm{J}}_{i}(1,t)$ and $c^{\mathrm{J}}_{i}(2,t)$; \textit{(ii)} the last segment $c^{\mathrm{P}}_{ij}(s_{L_{ij}},t+\Delta t)$ relies on the previous time-step $c^{\mathrm{J}}_{j}(t)$, $c^{\mathrm{J}}_{j}(s_{L_{ij}},t)$, and $c^{\mathrm{J}}_{j}(s_{L_{ij}}-1,t)$ in Pipe $ij$. 

\begin{table*}[t]
	%\resizebox{\textwidth}{!}{%
%	\vspace{-0.3cm}
	%    \scriptsize
	\fontsize{8.3}{8.3}\selectfont
	\centering
	\setlength\tabcolsep{1 pt}
	\renewcommand{\arraystretch}{1.3}
	\caption{Water quality modeling and corresponding matrix forms. }
	\begin{tabular}{c|c|c}
		\hline
		\textit{\makecell{Compo-\\nent}} & {\textit{Solute balance equation}} & {\textit{Matrix form}}\\ \hline
		\textit{Pipe} & \parbox{9.2cm}{
			\vspace{-0.7em}
			\begin{align}~\label{equ:adv-reac-lax}
			\hspace{-0.2em}    c^\mathrm{P}_{ij}(s,t+\Delta t) &= \underline{\alpha}_{ij}(t) c_{ij}^\mathrm{P}(s- 1,t) 
			+ {\alpha}_{ij}(t)c^\mathrm{P}_{ij}(s,t)\\
			&+\overline{\alpha}_{ij}(t) c^\mathrm{P}_{ij}(s+ 1,t)
			+ r^\mathrm{P}_{ij}(c^\mathrm{P}_{ij}(s,t)) \notag
			\end{align}
			\vspace{-1.1em}
		}  
		&\parbox{6.8cm}{
			\vspace{-0.7em}
			\begin{align} ~\label{equ:abs-pipe-mass}
			\m c^\mathrm{P}(t+ \Delta t) = &\m A^\mathrm{N}_\mathrm{P}(t) \m c^\mathrm{N}(t)+ \m A^\mathrm{P}_\mathrm{P}(t) \m c^\mathrm{P}(t)  \\ 
			&+ \m r^\mathrm{P}(\m c^\mathrm{P}(t)) \notag
			\end{align}
			\vspace{-1.1em}
		}\\ \hline
		\textit{\makecell{Junc-\\tion}} & {    \parbox{9.2cm}{
				\vspace{-0.7em}
				\begin{align} ~\label{equ:mass-balance-junc}
				\hspace{-1.5em}
				q^\mathrm{B}_i(t) c^\mathrm{B}_i(t)     
				&+ \textstyle \sum_{k = 1}^{|\mathcal{N}_i^\mathrm{in}|} q_{ki}(t)c_{ki}(t)  \\
				&=
				q^\mathrm{D}_i(t) c^\mathrm{J}_{i}(t) + \textstyle \sum_{j = 1}^{|\mathcal{N}_i^\mathrm{out}|} q_{ij}(t)  c_{ij}(t) \notag
				\end{align}
				\vspace{-1.1em}
		}} 
		& \parbox{6.8cm}{
			\vspace{-0.7em}
			\begin{align}~\label{equ:abs-nodes-mass}
			\hspace{-2em}  &\m c^\mathrm{J}(t\hspace{-2pt}+\hspace{-2pt}\Delta t)  =  \m A_\mathrm{J}^\mathrm{J}(t \hspace{-2pt}+\hspace{-2pt} \Delta t) \m c^\mathrm{J}(t)  + \m A_\mathrm{J}^\mathrm{L}(t\hspace{-2pt}+\hspace{-2pt}\Delta t) \m c^\mathrm{L}(t)  \notag \\
			&+  \m B^\mathrm{J}(t\hspace{-2pt}+\hspace{-2pt}\Delta t)   \m c^\mathrm{B}(t\hspace{-2pt}+\hspace{-2pt}\Delta t) + \m R^\mathrm{J}(t \hspace{-2pt}+\hspace{-2pt} \Delta t) \m r(\m x(t)) 
			\end{align}
			\vspace{-1.1em}
		}
		\\ \hline
		\textit{Tank} &\parbox{9.2cm}{ 
			\vspace{-0.7em}
			\begin{align} ~\label{equ:mass-balance-tank}
			\hspace{-0.5em} &V^{\mathrm{TK}}_{i}(t + \Delta t) c_{i}^{\mathrm{TK}}(t +  \Delta t) \hspace{-2pt}=\hspace{-2pt}( V^{\mathrm{TK}}_i(t)- \Delta t  \textstyle \sum_{j = 1}^{|\mathcal{N}_i^\mathrm{out}|} q_{ij}(t)) c_i^{\mathrm{TK}}(t)  \\
			&+ V^\mathrm{B}_i(t+  \Delta t) c_{i}^{\mathrm{B}}(t+  \Delta t) +  \Delta t  (\textstyle \sum_{k = 1}^{|\mathcal{N}_i^\mathrm{in}|} q_{ki}(t)c_{ki}(s_{L_{ki}},t) +  r^{\mathrm{TK}}_{i} (c^{\mathrm{TK}}(t))  )  \notag 
			\end{align}
			\vspace{-0.6em}
		} 
		& \parbox{6.8cm}{
			\vspace{-0.7em}
			\begin{align} ~\label{equ:abs-tank-mass}
			\hspace{-0.8em}   &\m c^{\mathrm{TK}}(t\hspace{-2pt}+\hspace{-2pt}\Delta t) =  \m A_{\mathrm{TK}}^\mathrm{TK}(t) \m c^{\mathrm{TK}}(t) + \m A_{\mathrm{TK}}^\mathrm{P}(t) \m c^{\mathrm{P}}(t) \notag \\ 
			&+ \m B^\mathrm{TK}(t\hspace{-2pt}+\hspace{-2pt}\Delta t) \m c^{\mathrm{B}}(t\hspace{-2pt}+\hspace{-2pt}\Delta t) \hspace{-2pt}+\hspace{-2pt}\m R^\mathrm{TK}(t) \m r(\m x(t))
			\end{align}
			%\vspace{-0.6em}
		}
		\\ \hline \hline
	\end{tabular}
	\label{table:wqp}
\end{table*}
Equation~\eqref{equ:adv-reac-lax} can be lumped in a matrix-vector form for all the pipes as given by~\eqref{equ:abs-pipe-mass} in Tab.~\ref{table:vector}, where vector $\m c^\mathrm{P}(t) \triangleq \Bigl \lbrace \m c_1^\mathrm{P}(t),\ldots, \m c_{n_\mathrm{P}}^\mathrm{P}(t) \Bigr \rbrace$ collects the concentrations  at all grid points in all pipes with $n_\mathrm{P}$ being the number of pipes in the water network; $\m c^\mathrm{N}(t)$ defined in Tab.~\ref{table:vector} collects all concentration from nodes; matrix $\m A_\mathrm{P}^\mathrm{N}(t)$ derived from~\eqref{equ:node2pipe} stands for the concentration contribution from nodes to the segments in pipes; matrix $\m A_\mathrm{P}^\mathrm{P}(t)$ lumps the L-W coefficients defined in~\eqref{equ:alphas} representing the contribution from pipes at previous time step to current one. 

%\begin{myrem}~\label{rem:Pipe}
%	The concentration difference equations from the non-first segment of pipes ($\mathcal{P} \char`\\ \mathcal{F}$) are obtained so far; concentration in the first segment of pipes ($\mathcal{F}$) is available after concentrations in their upstream nodes are obtained in next subsection. In the end, the matrix form is formed as~\eqref{equ:abs-pipe-mass}.
%\end{myrem}

\begin{myrem}~\label{rem:Pipe}
	The concentrations difference equations  for pipes $\mathcal{P}$ are obtained.
\end{myrem}

\subsubsection{Solute mass balance at junctions}~\label{sec:mbjunction} Mass conservation of the disinfectant for Junction $i$ at time $t$ can be described by~\eqref{equ:mass-balance-junc} in Tab.~\ref{table:wqp} when Assumption~\ref{asmp:mixing} is made, where $ \{ ki : k \in \mathcal{N}_i^\mathrm{in} \} $ and  $ \{ ij : j \in \mathcal{N}_i^\mathrm{out} \} $ represent the sets of links with inflows and outflows of Junction $i$; $c_{ki}(t)$ and $c_{ij}(t)$ are the concentrations in Links $ki$ and $ij$.
\begin{asmp}\label{asmp:mixing}
	The mixing of the solute is complete and instantaneous  at junctions and  in tanks with a continuously stirred tank reactors (CSTRs) model.
\end{asmp}
{Assumption~\ref{asmp:mixing} is widely used in~\cite{helbling2009modeling,shang2008epanet,boulos2006comprehensive,rossman1996numerical}. It implies \textit{(i)} that the solute injected from boosters takes effect to junctions and tanks immediately instead of delays, and \textit{(ii)} the concentration leaving a node (a junction or a tank), either by user demands or in  downstream links (a pump, a pipe, or a valve), is the same. }

Note that we only can obtain the concentration difference equation at junctions in $\mathcal{W} \char`\\ \mathcal{D}$ so far (such as $\mathrm{J}3 \in \mathcal{W} \char`\\ \mathcal{D}$ which is the upstream node of $\mathrm{V}34$ in Fig.~\ref{fig:example3}); for the junctions in $\mathcal{D}$ (such as $\mathrm{J}2, \mathrm{J}4 \in \mathcal{D}$ which are the downstream nodes  $\mathrm{M}12$ and $\mathrm{V}34$ in Fig.~\ref{fig:example3}),  we can not obtain them until their dependence links are available in Section~\ref{sec:PumpandValve}; {see Example~\ref{exmp:junction} in Appendix~\ref{app:example} for the details of difference equation at $\mathrm{J}3$ in Fig.~\ref{fig:example3}}. 

%We next show how  to derive $c_{3}^\mathrm{J}(t + \Delta t)$ via the example shown in Fig.~\ref{fig:example3}. 

Now, we consider for all {junctions}, that is, Equation~\eqref{equ:mass-balance-junc} can be lumped in a matrix form as
	\begin{align} ~\label{equ:tempJunc}
& \m q^\mathrm{B}_\mathrm{J}(t \hspace{-1.5pt}+\hspace{-1.5pt}  \Delta t) \circ \m c^\mathrm{B}_\mathrm{J}(t \hspace{-1.5pt}+\hspace{-1.5pt}  \Delta t)  +\m q^\mathrm{in}_\mathrm{J}(t \hspace{-1.5pt}+\hspace{-1.5pt}  \Delta t)\circ (\diag{(\m S^\mathrm{in}_\mathrm{J})}\m c^\mathrm{L}(\m s_{L},t \hspace{-1.5pt}+\hspace{-1.5pt}  \Delta t))\notag
\\
&= \left(\m q^{\mathrm{out}}_\mathrm{J}(t \hspace{-1.5pt}+\hspace{-1.5pt}  \Delta t)+\m q^\mathrm{D}(t \hspace{-1.5pt}+\hspace{-1.5pt}  \Delta t)\right)  \circ \m c^\mathrm{J}(t \hspace{-1.5pt}+\hspace{-1.5pt}  \Delta t) 
\end{align}
where $\m q^\mathrm{D}$ is demand ; $\m q^{\mathrm{in}}_\mathrm{J}(t+\Delta t)$ and $\m q^{\mathrm{out}}_\mathrm{J}(t+\Delta t)$ are inflows and outflows defined by ${\diag{(\m S^\mathrm{in}_\mathrm{J})}} \m q^\mathrm{L}(\m s_{L},t+\Delta t)$ and $\m S^\mathrm{out}_\mathrm{J} \m q^\mathrm{L}(1,t+\Delta t)$, and {note that each element in vector $\m q^{\mathrm{out}}_\mathrm{J}(t+\Delta t)$ is the summation of all outflows from a corresponding junction}; selection matrix $\m S^\mathrm{in}_\mathrm{J}$ and $\m S^\mathrm{out}_\mathrm{J}$ is from  $\m E^\mathrm{L}_\mathrm{J}$ in~\eqref{equ:WFP-Incidence-matrix}, {we give Example~\ref{exmp:selectmatrix4J} in Appendix~\ref{app:example} to explain the detail}.  Note that $\m c^\mathrm{L}(\m s_{L},t+\Delta t) \in \mathbb{R}^{n_\mathrm{Q}}$ only includes concentration of last segments in pipes, and $\m q^\mathrm{L}(1,t+\Delta t)$ includes the first segments; $\m q^\mathrm{B}_\mathrm{J} = \m E^\mathrm{B}_\mathrm{J} \m q^\mathrm{B}$ and $ \m c^\mathrm{B}_\mathrm{J} =\m E^\mathrm{B}_\mathrm{J} \m c^\mathrm{B}$ are the flow rate vectors and concentration vectors injected by boosters at junctions, where $\m E^\mathrm{B}_\mathrm{J}$ shown in~\eqref{equ:WFP-Booster-matrix} describes booster locations at junctions. Here, $\m c^\mathrm{L}(\m s_{L},t+\Delta t)$, that is a subvector of  $\m  c^\mathrm{L}(t+\Delta t)$,  includes concentration of the last segment of each pipe $\m c^\mathrm{P}(\m s_{L},t+\Delta t)$, and concentrations of pumps $\m c^\mathrm{M}(t+\Delta t)$ and valves $\m c^\mathrm{V}(t+\Delta t)$.

After expressing the result of Hadamard product and division as matrix multiplication form, we obtain matrix form~\eqref{equ:abs-nodes-mass} in Tab.~\ref{table:wqp}, where  matrices $\m A_\mathrm{J}^\mathrm{J}$, $\m A_\mathrm{J}^\mathrm{L}$, and  $\m B^\mathrm{J}$ represent the contribution to $\m c^\mathrm{J}(t+\Delta t)$ from corresponding junctions, links, and boosters. Moreover, $\m R^\mathrm{J}$ stands for the impact from  nonlinear  reaction item $\m r(\m x(t))$ defined in later Section~\ref{sec:nonlinearSpace}. The  detail of matrix derivation for junctions are in Appendix~\ref{app:derivation_Junction}.

\subsubsection{Solute mass balance at tanks}   Tanks are typically filled during low demand periods, and are drained when going through high demand periods. Mass conservation of  disinfectants in Tank $i$ can be expressed by~\eqref{equ:mass-balance-tank}, where variables are defined in Tab.~\ref{table:vector}. The physical meaning of~\eqref{equ:mass-balance-tank} is that the solute mass at $t+\Delta t$ equals the current mass plus the change of mass due to boosters, inflows, outflows, and reactions in tanks. Note that  the mixing in tanks is instantaneous, see Assumption~\ref{asmp:mixing}. Hence, the $c_{i}^{\mathrm{B}}(t +  \Delta t)$ instead of $c_{i}^{\mathrm{B}}(t)$  has effects to  $c_{i}^{\mathrm{TK}}(t +  \Delta t)$. %We give a specific tank in Fig.~\ref{fig:example3} to illustrate~\eqref{equ:mass-balance-tank} further. 

After listing equations for all tanks, matrix form of~\eqref{equ:mass-balance-tank} can be expressed as
%\begin{linenomath*}
\begin{align} 
 &\m V^{\mathrm{TK}}(t+\Delta t) \circ \m c^{\mathrm{TK}}(t+\Delta t)~\label{equ:abs-tank-temp} \\
 &= \hspace{-2pt} (\m V^{\mathrm{TK}}(t) \hspace{-2pt}-\hspace{-2pt} \Delta t \,  \m q^{\mathrm{out}}_\mathrm{TK}(t) ) \hspace{-2pt}\circ \hspace{-2pt} \m c^{\mathrm{TK}}(t) \hspace{-2pt}+\hspace{-2pt}  \m V^\mathrm{B}_\mathrm{TK}(t\hspace{-2pt}+\hspace{-2pt}\Delta t)   \hspace{-2pt}\circ \hspace{-2pt} \m c^{\mathrm{B}}_\mathrm{TK}(t+\Delta t)  +\Delta t \left(\m q^\mathrm{in}_\mathrm{TK}(t)  \circ \m c^{\mathrm{P}}(t) +  \m r^\mathrm{TK}(\m c^\mathrm{TK}(t))\right), \notag
\end{align}
%\end{linenomath*}
where all variables have similar meaning defined for junctions, except that $\m q^{\mathrm{in}}_\mathrm{TK}$ and $\m q^{\mathrm{out}}_\mathrm{TK}$ are inflow and outflow subvectors, and can be obtained from multiplication of $\m q^\mathrm{L}$ with the selection matrices from $\m E^\mathrm{L}_\mathrm{TK}$; and  $\m c^\mathrm{B}_\mathrm{TK}  = \m E^\mathrm{B}_\mathrm{TK} \m c^\mathrm{B}$, and $\m E^\mathrm{B}_\mathrm{TK}$ in~\eqref{equ:WFP-Booster-matrix} defines the booster location at tanks; {see Example~\ref{exmp:tank} in Appendix~\ref{app:example} for the details of difference equation at $\mathrm{TK}5$ in Fig.~\ref{fig:example3}}.

Akin to dealing with~\eqref{equ:tempJunc} to simplify \eqref{equ:abs-tank-temp}, we obtain~\eqref{equ:abs-tank-mass} in Tab.~\ref{table:wqp} for tanks; see Appendix~\ref{app:derivation_Tank} for the detail of derivation. 
\subsubsection{Solute mass balance at reservoirs} 
We  assume that the concentration at a reservoir is constant, that is 
%\begin{linenomath*}
	\begin{align}\label{equ:abs-reservoir-i}
	c_i^\mathrm{R}(t + \Delta t) &= c_i^\mathrm{R}(t) ,
	\end{align}
%\end{linenomath*}
and the corresponding matrix form for all reservoirs is
%\begin{linenomath*}
	\begin{align}\label{equ:abs-reservoir}
	\m c^\mathrm{R}(t + \Delta t) &= \m c^\mathrm{R}(t) .
	\end{align}
%\end{linenomath*}
Example~\ref{exmp:reservoir} in Appendix~\ref{app:example} shows the difference equation at $\mathrm{R}1$ in Fig.~\ref{fig:example3}.

\begin{myrem}~\label{rem:Junctions}  
	The concentrations difference equations  from \textit{(i)} nodes having no connection with pumps and valves and \textit{(ii)}  the upstream nodes of pumps and valves, that are junctions, tanks, and reservoirs in $\mathcal{W} \char`\\ \mathcal{D}$, are obtained so far; for the downstream nodes in ${\mathcal{D}}$, we can obtain them after  concentrations in pumps and valves are derived in next section. In the end, the matrix forms of junctions, tanks, and reservoirs are formed as~\eqref{equ:abs-nodes-mass},~\eqref{equ:abs-tank-mass}, and~\eqref{equ:abs-reservoir}.
\end{myrem}

\subsubsection{Solute transport in pumps and valves}~\label{sec:PumpandValve}
The lengths of pumps and valves are assumed to be zeros, and they do not store any water. Therefore, we assume that the concentration at pumps or valves equal concentration of the upstream nodes they connect. The corresponding matrices form for pumps and valves can be written as 
%\begin{linenomath*}
	\begin{align}
	\m c^\mathrm{M}(t+\Delta t) &= \m S^\mathrm{N}_\mathrm{M} \m c^\mathrm{N}(t+\Delta t) \label{equ:pump}\\
	\m c^\mathrm{V}(t +\Delta t) &=  \m S^\mathrm{N}_\mathrm{V} \m c^\mathrm{N}(t+\Delta t), \label{equ:valve}
	\end{align}
%\end{linenomath*}
where all vectors defined in Tab.~\ref{table:vector}, and $\m S^\mathrm{N}_\mathrm{M}$ and $\m S^\mathrm{N}_\mathrm{V}$ are selection matrices and  can be obtained by  changing $-1$ to $0$ in the node-pump and node-valve connectivity matrices, that are $\m E^\mathrm{N}_\mathrm{M}$ and $\m E^\mathrm{N}_\mathrm{V}$ in~\eqref{equ:WFP-Incidence-matrix}; {see Example~\ref{exmp:selectmatrix} in Appendix~\ref{app:example}}. 

When the upstream node of Pump $ij$ is a reservoir, and it means the concentration in Pump $ij$ equals the one in Reservoir $i$. In this case, Equation~\eqref{equ:pump} is rewritten as
	\begin{subequations}~\label{equ:abs-pump}
		\begin{align}~\label{equ:pumpi}
		c_{ij}^\mathrm{M}(t+\Delta t) =  c_i^\mathrm{R}(t + \Delta t) =  c_i^\mathrm{R}(t) =  c_{ij}^\mathrm{M}(t).
		\end{align}
		When the upstream node of Pump $ij$ is a junction, and the concentration in Pump $ij$ equals the one at Junction $i$. Hence, Equation~\eqref{equ:pump} is rewritten as
			\begin{align}~\label{equ:pumpjunction}
			c_{ij}^\mathrm{M}(t+\Delta t) =  c_i^\mathrm{J}(t + \Delta t),
			\end{align}
	\end{subequations}
where each $c_i^\mathrm{J}(t + \Delta t)$ can be obtained from~\eqref{equ:abs-nodes-mass} directly in Section~\ref{sec:mbjunction}. There are no differences between pumps and valves from the point of view of solution transport except that valves connect to junctions and tanks. {Therefore, we present  Example~\ref{exmp:pump} to demonstrate it further in Appendix~\ref{app:example}.} %of the upstream node of a pump (a valve) is a reservoir (a junction) .

Considering all pumps for all cases,  $\m c^\mathrm{M}(t+\Delta t)$ can always be expressed by the concentrations at nodes and in links from previous time step $t$. In practice,  pumps  connect to either junctions or reservoirs, and the corresponding connectivity matrices are $\m E^\mathrm{M}_\mathrm{J}$ and $\m E^\mathrm{M}_\mathrm{R}$ in~\eqref{equ:WFP-Incidence-matrix}. Besides that, the matrix forms of junction and reservoir solute equations are known as~\eqref{equ:abs-nodes-mass} and~\eqref{equ:abs-reservoir}.  Hence, the matrix form of pump solute equation can be obtained directly by multiplying selection matrices ($\m S^\mathrm{M}_\mathrm{J}$ and $\m S^\mathrm{M}_\mathrm{R}$) with corresponding matrix forms given in~\eqref{equ:abs-nodes-mass} and~\eqref{equ:abs-reservoir}.

Similarly, the matrix form of valve solute equations can be expressed by multiplying selection matrices $\m S^\mathrm{V}_\mathrm{J}$ and $\m S^\mathrm{V}_\mathrm{TK}$ with their corresponding  matrix forms~\eqref{equ:abs-nodes-mass}  and~\eqref{equ:abs-tank-mass}. The notations helping to present them as matrix form would be introduced in~Section~\ref{sec:nonlinearSpace}. Hence, we prefer to show the matrix forms for pumps and valves in that section.

\begin{myrem}~\label{rem:PumpValve}
	The concentrations difference equations  for pumps and valves $\mathcal{M} \cup \mathcal{V}$ are obtained.
\end{myrem}

We clearly see the concentration difference equation of a component  might depend on concentration equation of another from Remarks~\ref{rem:Pipe},~\ref{rem:Junctions}, and~\ref{rem:PumpValve}, and we refer this property in water quality modeling as \textit{concentration dependence}. For example, concentration in links, such as valves and pumps, depend on their upstream nodes while their downstream nodes rely on the links in turn. This implies that there exists a specific order to list concentration equations for all components, and the order can be defined by a \textit{concentration dependence tree or forest}. 
Details are discussed next.  

\subsection{Steps to find difference equations for components based on  dependence forest }~\label{sec:dependence}
In this section, we organize the dependence relationship presented in Section~\ref{sec:conservation} and summarize steps to find difference equation of each component as Algorithm~\ref{alg:matrix}, and we use the illustrative example in Fig.~\ref{fig:example3} to explain the detail, {see Example~\ref{exmp:dependence} in Appendix~\ref{app:example}}.
%\begin{algorithm}[t]
%	\begin{algorithmic}[1]
%		\STATE \textbf{Input:} WDN typologies and parameters
%		\STATE \textbf{Output:} Water quality modeling
%		\STATE Run hydraulic simulation at $t$ for flow direction 
%		\STATE Find  concentrations  in set $ \mathcal{P} \char`\\ \mathcal{F}$ via~\eqref{equ:adv-reac-lax}, that are in non-first segments for  $ij \in \mathcal{P}$
%		\STATE Find $c_i(t+\Delta t)$ via~\eqref{equ:mass-balance-junc},~\eqref{equ:mass-balance-tank}, or~\eqref{equ:abs-reservoir-i} for $i \in \mathcal{W} \char`\\ \mathcal{D}$  
%		\STATE Find $c_{ij}(t+\Delta t)$ via~\eqref{equ:abs-pump} for Pump or Valve $ij \in \mathcal{M} \cup \mathcal{V}$  
%		\STATE Find $c_i(t+\Delta t)$ via ~\eqref{equ:mass-balance-junc} or~\eqref{equ:mass-balance-tank} for $i \in \mathcal{D}$
%		\STATE Find  concentrations  in set $\mathcal{F}$ via~\eqref{equ:adv-reac-lax}, that are in the first segments for  $ij \in \mathcal{P}$
%	\end{algorithmic}
%	\caption{Obtain water quality modeling for components}
%	\label{alg:matrix}
%\end{algorithm}

\begin{algorithm}[t]
	\small	\DontPrintSemicolon
	\KwIn{WDN typologies and parameters}
	\KwOut{Difference equations of all components}
	\While {  $t \leq T_d$ }{
		Run hydraulic simulation at $t$ for flow direction \;
		Update $\m E_\mathcal{G}$ according to the true flow direction in links\;
		{Obtain   $c_{ij}(s,t+\Delta t)$  via~\eqref{equ:adv-reac-lax} in set $\mathcal{P}$}\;
		Obtain $c_i(t+\Delta t)$ via~\eqref{equ:mass-balance-junc},~\eqref{equ:mass-balance-tank}, or~\eqref{equ:abs-reservoir-i} for $i \in \mathcal{W} \char`\\ \mathcal{D}$ \;
		Obtain $c_{ij}(t+\Delta t)$ via~\eqref{equ:abs-pump} for Pump or Valve $ij \in \mathcal{M} \cup \mathcal{V}$, and   construct matrix form~\eqref{equ:pump} and~\eqref{equ:valve} \; 
		Obtain $c_i(t+\Delta t)$ via ~\eqref{equ:mass-balance-junc} or~\eqref{equ:mass-balance-tank} for $i \in \mathcal{D}$, and construct matrix form \eqref{equ:abs-nodes-mass}, \eqref{equ:abs-tank-mass}, \eqref{equ:abs-reservoir}\;
		%Obtain  concentrations  $c_{ij}(1,t+\Delta t)$  in set $\mathcal{F}$ via~\eqref{equ:adv-reac-lax}, that are in the first segments for  $ij \in \mathcal{P}$, and construct matrix form~\eqref{equ:abs-pipe-mass} \;
		$ t = t + \Delta t$
	}
	\caption{Offline, time-varying difference equation formulation of all WDN components}
	\label{alg:matrix}
	With difference equations of all component, matrices \blue{$\m A(t)$ and $\m B(t)$}  for state-space form are derived in Section~\ref{sec:abstract-WQM}.
\end{algorithm}

From this dependence forest in Example~\ref{exmp:dependence}, we see that the order to find difference equations is unique and fixed for components, that is, {$\mathcal{P} \rightarrow \mathcal{\mathcal{W} \char`\\ \mathcal{D}} \rightarrow \mathcal{M \cup V}  \rightarrow \mathcal{D}$, see Fig.~\ref{fig:Example2findMatrix}.} This order can be mapped to  steps in Algorithm~\ref{alg:matrix}. In this algorithm, the total duration time is $T_d$, and for each time-step $t$, we need to update $E_\mathcal{G}$ and replace the assigned direction with true flow direction, see Steps 2--3. When the concentration equations of all components are available, they can be lumped together into matrix form (just as the one in Tab.~\ref{table:wqp}) via Steps 4--7.  Next, we discuss how to derive the state-space form for overall network based on  matrix form of these equations.
\section{\large Water Quality Modeling in State-Space Form}~\label{sec:abstract-WQM}
This section presents an abstract, state-space model that can be useful for a wide range of control theoretic studies in water quality control studies.  In particular, we showcase a mathematical model comprising a set of inputs, outputs, and state variables related by first-order difference equations. First,  connectivity matrix $ \m E^\mathrm{B}_\mathrm{N} $ is assumed to be pre-determined and time-independent. Second, for a water quality control time-horizon, it is customary to assume that all hydraulic variables are assumed to be known~\cite{wang2005adaptive}, including the flow rate $\m q^{\mathrm{L}}$, demand $\m q^\mathrm{D}$, tank volumes $\m V^{\mathrm{TK}}$.  This is due to hydraulic simulation has a much slower time-scale  compared with the one in a water quality simulation. {Note that this assumption in practice is a big issue, since the hydraulic variables are impossible to be measured directly. Fortunately, hydraulic state estimation can be used to obtain the $\m q^{\mathrm{L}}$, demand $\m q^\mathrm{D}$, tank volumes $\m V^{\mathrm{TK}}$, and even pipe roughness coefficients; see~\cite{tshehla2017state,wang2019state}.}  {Besides that the flow rate and corresponding volume injected from  boosters are known, and assumed as $\m V^\mathrm{B}= \m q^\mathrm{B}$, that is, it takes boosters unit time to inject solute into networks}.

We define the state and input vectors first in ensuing sections, and then show how matrix form~\eqref{equ:abs-pipe-mass}, \eqref{equ:abs-nodes-mass}, \eqref{equ:abs-tank-mass}, \eqref{equ:abs-reservoir},~\eqref{equ:pump}, and~\eqref{equ:valve} in previous section can form the state-space form.

\subsection{Nonlinear state-space form}~\label{sec:nonlinearSpace}
System state $\m x(t)$ is a vector collecting concentrations at all nodes and in all links defined as
%\begin{linenomath*}
	\begin{eqnarray}
	\m x(t) \triangleq \{\m c^\mathrm{N}(t),\m c^\mathrm{L}(t) \}  =  \{\m c^\mathrm{J}(t), \m c^\mathrm{R}(t), \m c^\mathrm{TK}(t),  \m c^\mathrm{P}(t), \m c^\mathrm{M}(t), \m c^\mathrm{V}(t) \}, ~\label{equ:xVector}
	\end{eqnarray}
%\end{linenomath*}
{where the  dimension of $\m x$ is the summation of number of concentration variables in nodes and at links, that is, $n_x = n_\mathrm{N}  + n_\mathrm{L}$.}

{We define $\m u(t)$ standing for the input of the system at time $t$ as $\m c^\mathrm{B}(t+\Delta t)$ instead of $\m c^\mathrm{B}(t)$ in our paper due to Assumption~\ref{asmp:mixing}. That is, if there is no delay when mixing of the solute, the injections or inputs happening at $t$ and $t + \Delta t$ at junctions or tanks can be considered as the same ``moment". Hence, the control variable is still defined as $\m u(t) \triangleq  \m c^\mathrm{B}(t+\Delta t)$, and the dimension of $\m u$ is the same as the dimension of $\m c^\mathrm{B}$, that is, $n_u = n_{\mathrm{N}}$.}

%The input $\m u(t)$ affects  on a small proportion of states in $\m x(t)$, that is $\m c^\mathrm{J}$ and $\m c^\mathrm{TK}$. However, states in links, especially $\m c^\mathrm{P}$, occupying the majority of states in $\m x(t)$ would not be changed or affected by inputs instantaneously since it takes  time for disinfectants to transport in pipes.

% Normally, control variable $\m u(t)$ stands for the input of the system at time $t$, that is, its impact only  starts to show up at the next time-step $t+\Delta t$. But in our case, the situation is different due to Assumption~\ref{asmp:mixing},  

The reaction rate vector collecting concentrations in pipes and tanks is defined as
%\begin{linenomath*}
	\begin{align}~\label{equ:reactionVector}
	\m r(\m x(t)) \triangleq\{\m  r^\mathrm{TK}(\m c^\mathrm{TK}(t)),\m r^\mathrm{P}(\m c^\mathrm{P}(t))\}.
	\end{align}
%\end{linenomath*}
We also define $\m R_r$ as an  indicator matrix encodes how the nonlinearity $\m r(\m x(t))$ is impacting the dynamics of water quality, i.e.,  tanks and pipes having reaction rate could impact $\m x$, {each component has its own indicator submatrix, and the detail of elements in $\m R_r$ would be discussed in ensuing paragraphs.}

With the above definitions, we next show how~\eqref{equ:abs-pipe-mass}, \eqref{equ:abs-nodes-mass}, \eqref{equ:abs-tank-mass}, \eqref{equ:abs-reservoir},~\eqref{equ:pump}, and~\eqref{equ:valve} can yield the compact, state-space, and control-oriented formulation~\eqref{equ:fullMatrix} which is a Nonlinear Difference Equation (NDE).

% \begin{figure*}
	\begin{equation}~\label{equ:fullMatrix}
	\hspace{-1.8em}\includegraphics[width=0.85\linewidth,valign=c]{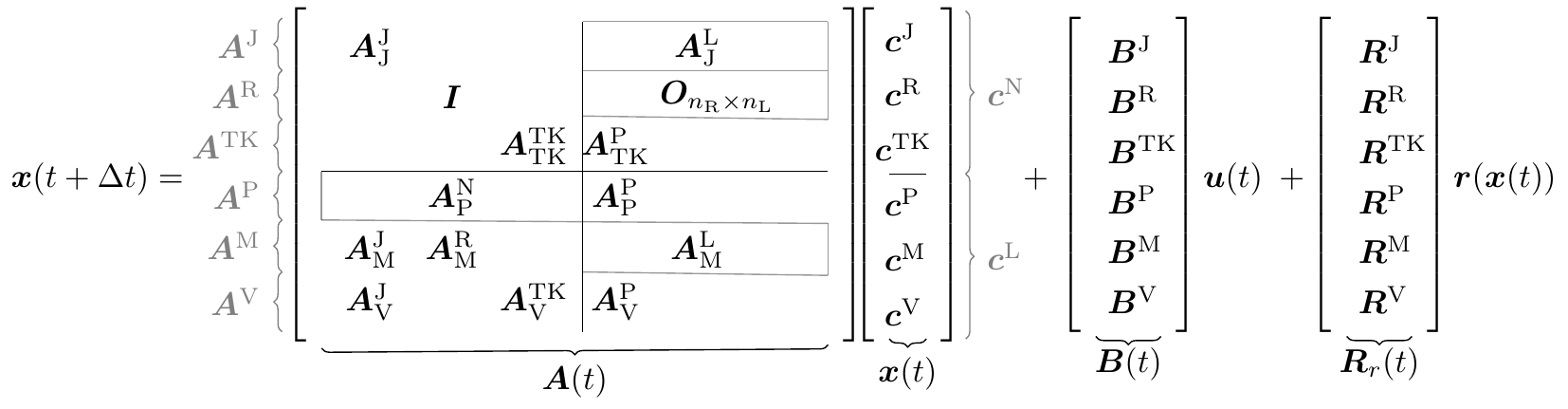}
	\end{equation}
%	\hrulefill
%%	\vspace{-1.5em}
%\end{figure*}

We observe that all matrix forms except~\eqref{equ:pump} and~\eqref{equ:valve} are in difference equation form from Section~\ref{sec:Ctrl-WQM}, so we just need to organize and put  their contribution matrices into the right block in~\eqref{equ:fullMatrix}.  Taking concentration difference equations at junctions as an example, Equation~\eqref{equ:abs-nodes-mass} in Tab.~\ref{table:wqp} can be rewritten as
%\begin{linenomath*}
	\begin{align*}
	\m c^\mathrm{J}(t +\Delta t) = \underbrace{\begin{bmatrix}
		\m A^\mathrm{J}_\mathrm{J} & \m O& \m O & \m A^\mathrm{L}_\mathrm{J} 
		\end{bmatrix}}_{\large \m A^\mathrm{J}}\m c^\mathrm{J}(t) +
\m B^\mathrm{J} \m u(t) + \underbrace{
		\begin{bmatrix} 
		\m O &  \m R^\mathrm{P}_\mathrm{J} 
		\end{bmatrix}
	}_{\large \m R^\mathrm{J}} \m r(\m x(t)),
	\end{align*}
%\end{linenomath*}
and we put it as  the first  block-row in~\eqref{equ:fullMatrix}. {Note that the elements in matrix $\m R^\mathrm{P}_\mathrm{J}  \in \mbb{R}^{n_\mathrm{J} \times n_\mathrm{S}}$ are the reaction rates passing from segments in the pipes to junctions; see $R_3^\mathrm{J}(t + \Delta t)$ which is a typical element of $\m R^\mathrm{P}_\mathrm{J}$ in Example~\ref{exmp:junction} for the detail. Similarly, reservoirs~\eqref{equ:abs-reservoir}, tanks~\eqref{equ:abs-tank-mass}, and pipes concentration equation~\eqref{equ:abs-pipe-mass} are put in  the second, third, and fourth block-rows in~\eqref{equ:fullMatrix}. The corresponding indicator matrix for reservoirs $\m R^\mathrm{R}$ are zero matrix since no decay is assumed in reservoirs, and no components can impact on it either;  indicator matrix for tanks and pipes are $\m R^\mathrm{TK} =[\m I_{n_{TK}} \   \m O]$ and $\m R^\mathrm{P} = [\m O \ \m I_{n_{S}} )]$.} 

As for pumps~\eqref{equ:pump}, we know that it depends on concentration of the upstream node it connects from~\eqref{equ:abs-pump}. Hence, $\m A^\mathrm{M}$ is formed by  selecting  specific block-rows from $\m A^\mathrm{J}$ and $\m A^\mathrm{R}$, and this action can be done by  selection matrix $\m S^\mathrm{N}_\mathrm{M}$, that is
%\begin{linenomath*}
	\begin{align*}
	\m A^\mathrm{M} 
	&\triangleq \begin{bmatrix}
	\m A^\mathrm{J}_\mathrm{M} & \m A^\mathrm{R}_\mathrm{M} & \m O & \m A^\mathrm{L}_\mathrm{M} \end{bmatrix} = \begin{bmatrix}
	\m S^\mathrm{J}_\mathrm{M}  & \m S^\mathrm{R}_\mathrm{M}
	\end{bmatrix}
	\begin{bmatrix}
	\m A^\mathrm{J} \\
	\m A^\mathrm{R}
	\end{bmatrix}  = \begin{bmatrix}
	\m S^\mathrm{J}_\mathrm{M} \m A^\mathrm{J}_\mathrm{J} & \m S^\mathrm{R}_\mathrm{M} & \m O & \m S^\mathrm{J}_\mathrm{M} \m A^\mathrm{L}_\mathrm{J} 
	\end{bmatrix},\\
	\m B^\mathrm{M} &= \begin{bmatrix}
	\m S^\mathrm{J}_\mathrm{M}  & \m S^\mathrm{R}_\mathrm{M}
	\end{bmatrix} \begin{bmatrix}
	\m B^\mathrm{J} \\
	\m B^\mathrm{R}
	\end{bmatrix},
	\m R^\mathrm{M} = \begin{bmatrix}
	\m S^\mathrm{J}_\mathrm{M}  & \m S^\mathrm{R}_\mathrm{M}
	\end{bmatrix} \begin{bmatrix}
	\m R^\mathrm{J} \\
	\m R^\mathrm{R}
	\end{bmatrix}.
	\end{align*}
%\end{linenomath*}
%and $$. 

Akin to pumps,  submatrices $\m A^\mathrm{V}$, $\m B^\mathrm{V}$, and $\m R^\mathrm{V}$ for valves~\eqref{equ:valve} can be expressed as
%\begin{linenomath*}
	\begin{align*}
	&\m A^\mathrm{J}_\mathrm{V} = \m S^\mathrm{J}_\mathrm{V} \m A^\mathrm{J}_\mathrm{J}, \hspace{5pt} \m A^\mathrm{TK}_\mathrm{V} = \m S^\mathrm{TK}_\mathrm{V} \m A^\mathrm{TK}_\mathrm{TK}, \hspace{5pt} \m A^\mathrm{P}_\mathrm{V} = \m S^\mathrm{P}_\mathrm{V} \m A^\mathrm{P}_\mathrm{TK},\\
	& \m B^\mathrm{V} = \begin{bmatrix}
	\m S^\mathrm{J}_\mathrm{V}  & \m S^\mathrm{TK}_\mathrm{M}
	\end{bmatrix} \begin{bmatrix}
	\m B^\mathrm{J} \notag\\
	\m B^\mathrm{TK}
	\end{bmatrix}, \m R^\mathrm{V} = \begin{bmatrix}
	\m S^\mathrm{J}_\mathrm{V}  & \m S^\mathrm{TK}_\mathrm{M}
	\end{bmatrix} \begin{bmatrix}
	\m R^\mathrm{J} \notag\\
	\m R^\mathrm{TK}
	\end{bmatrix}.
	\end{align*}
%\end{linenomath*}

Furthermore, $\m A(t)$ and $\m B(t)$ are time-dependent matrices---due to their slow change with hydraulics---that depend on the WDN topology and the aforementioned mass balance equations.
Finally, the {{NDE}} model~\eqref{equ:fullMatrix} is possibly nonlinear due to a potential nonlinear reaction rate model $\m r(\cdot)$.

\subsection{Single species reaction model }~\label{sec:reaction}
In this section,  the first-order reaction for single species that describes disinfectant decay both in the bulk flow and at the pipe wall~\cite{rossman1996numerical,basha2007eulerian,shang2008epanet} are used, and 
%\begin{linenomath*}
	\begin{align}~\label{equ:reaction}
	k_{ij}^\mathrm{P} = k^{b}_{ij}+\frac{k^{w}_{ij}k^{f}_{ij}}{D_{ij}(k^{w}_{ij}+k^{f}_{ij})},\,\,\, k_{ij}^\mathrm{TK} = k^{b}_{ij},
	\end{align}
%\end{linenomath*}
where $ k^{b}_{ij}$ and $k^{w}_{ij}$ are defined in Tab.~\ref{table:vector}; $k^{f}_{ij}$ is the the mass transfer coefficient between the bulk flow and the pipe wall; $D_{ij}$ is the hydraulic diameter. For tanks,  only bulk reaction constant exists. Hence, the reaction rate  can be expressed as
%\begin{linenomath*}
	\begin{align}~\label{equ:reactionPTK}
	\hspace{-1.05em}\m r^\mathrm{P}(\m c^\mathrm{P}(t)) = \m k^\mathrm{P} \circ  \m c^\mathrm{P}(t),\,\m r^\mathrm{TK}(\m c^\mathrm{TK}(t)) = \m k^\mathrm{TK}  \circ   \m c^\mathrm{TK}(t).
	\end{align}
%\end{linenomath*}

{After using first-order reaction model, the elements in $\m R_r(t) \m r(\m x(t))$ in the NDE~\eqref{equ:fullMatrix} can be combined into $\m A$ matrix. For example,}
%\begin{linenomath*}
	$$\m A^\mathrm{P}_\mathrm{P}(t) =\m A^\mathrm{P}_\mathrm{P}(t)  + \diag(\m k^\mathrm{P} ), \;\; \m A^\mathrm{TK}_\mathrm{TK}(t) =\m A^\mathrm{TK}_\mathrm{TK}(t)  + \diag(\m k^\mathrm{TK} ).$$
%\end{linenomath*}
in~\eqref{equ:abs-pipe-mass} and~\eqref{equ:abs-tank-mass}, and we obtain a Linear DE (LDE) model as
%\begin{linenomath*}
	\begin{align}~\label{equ:de-abstract1}
	\mathrm{{LDE}:}\;\;\;    \m x(t + \Delta t) &= \m A(t) \m x(t) + \m B(t) \m u(t).
	\end{align}
%\end{linenomath*}
%where  the solution of LDE in state-space Form~\cite{friedland2012control} is
%$\m x(t)=\m \phi(t, 0) \m x(0)+\sum_{j=0}^{t-1} \m \phi(t, j+1) \m B(j) u(j)$ and $\m \phi(n, t)=\prod_{j=t}^{n-1} \m A(j)$. 
To validate the effectiveness of the proposed LDE~\eqref{equ:de-abstract1}, we compare the results of LDE~\eqref{equ:de-abstract1}  model with results of the water quality model embedded in  EPANET \cite{rossman2000epanet}. The results are presented in Section \ref{sec:test}. Given this linear state-space  form of the water quality model, the next section presents a  model control algorithm to control water quality. 
\begin{myrem}
	\textit{Matrices $\m A(t)$ written in~\eqref{equ:fullMatrix} is different from the one in~\eqref{equ:de-abstract1}, seeing that the former does not incorporate the reaction rate, whereas the latter embeds the linear, first-order reaction rate model. Examples are given in the case studies section. }
\end{myrem}

\section{\large Plug-and-Play MPC for WQC }~\label{sec:LP-WQC}
The objective of a water quality control (WQC) problem is to ensure safe drinking water at all nodes, while minimizing the amount of injected chlorine, {seeing that it reflects} higher operational costs. Model predictive control (MPC) technology is used here, due to its robustness to parametric uncertainty and forecast errors.

\begin{figure}[t]
	\centering
	\includegraphics[width=0.69\linewidth]{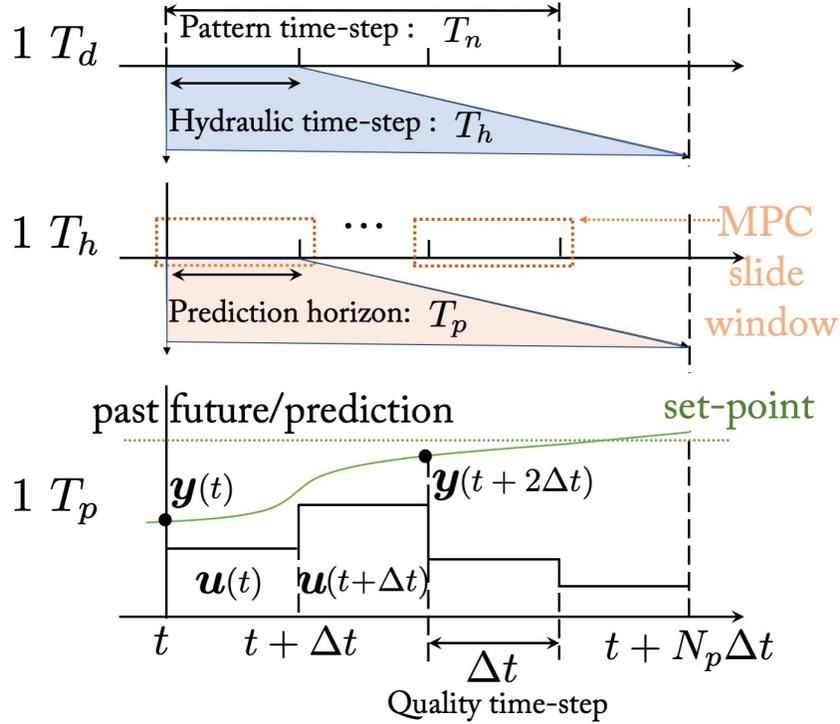}
	\caption{{Relationship among different time-scales for WQC problem and a discrete MPC scheme ($T_d = N_1 T_h = N_2 T_p = N_p \Delta t$, where $N_1$, $N_2$, and $N_p$ are integers).}}
	\label{fig:Timescale}
\end{figure}

\begin{table}[t]
	%    \scriptsize
	\fontsize{10}{10}\selectfont
	
	\centering
	\renewcommand{\arraystretch}{2}
	\makegapedcells
	\setcellgapes{1.4pt}
	\caption{{Time-scales in Water Quality Control problem.} }  
	\vspace{-0.2cm}
	\begin{tabular}{ c|c} 
		\hline
		\hspace{-4pt}\textit{Symbol}\hspace{-4pt} & \textit{Description}    \\ \hline
		$T_d$, $T_n$, $T_h$& \makecell{Duration of simulation (WQC problem time period); \\ demand pattern time-step; hydraulic time-step  } \\ \hline
		$T_p$, $T_c$, $\Delta t$& \makecell{MPC prediction horizon; MPC control horizon; \\ water quality time-step (in L-W scheme)}  \\ \hline
	\end{tabular} 
	\label{table:timescale}
	\vspace{-0.1cm}
\end{table}
\normalcolor

{Before discussing the objectives and constraints of  WQC problem, we would like to clarify different time-scales in WDN and MPC; see Tab.~\ref{table:timescale} for descriptions for different time-scales and Fig.~\ref{fig:Timescale} for the relationship among them. We note that the definitions of duration $T_d$, pattern time-step $T_n$, and hydraulic times step $T_h$ are exactly the same as the ones in EPANET. For example, the settings $1 \ T_d = 12\  T_n = 24\  T_h$ and $T_h = 1$ hour indicate that we are solving a hydraulic problem for 24 hours, demand changes every 2 hours, and the hydraulic time-step is 1 hour.}

We consider a WQC problem from $t = 0$ to $t=T_d$ which is constrained by lower and upper bounds constraints on chlorine concentrations in all links and nodes $c^{\min} = 0.2\ \mathrm{mg/L}$ and $ c^{\max} = 4\ \mathrm{mg/L}$, that is
%\begin{linenomath*}
	\begin{align} ~\label{equ:constr-abcstracted} 
	%    \hspace{-1.1em}    \textit{{Linear-Constraints}:}& \notag\\
	\m x(t) &\in [\m c^{\mathrm{min}}, \m c^{\mathrm{max}}].
	\end{align}
%\end{linenomath*}
Constraint~\eqref{equ:constr-abcstracted} which can succinctly be formulated as
%\begin{linenomath*}
	\begin{equation}
	\mathrm{{Constraints:}}\;\;\; \m x^{\min} \leq \m x(t) \leq \m x^{\max},~\label{equ:constr-abcstract-new} 
	\end{equation}
%\end{linenomath*}
where $\m x^{\min}$ and $\m x^{\max}$ model upper and lower bounds.

As for the objectives, the control problem to be addressed  is to minimize the total injected chlorine mass or the booster chlorination operational injection cost (BCI), while maintaining the chlorine concentration in all links and nodes $\m x(t)$ within pre-specified bounds discussed above. Given these considerations, the WQC problem in a prediction horizon $T_p$ can be written as
%\begin{linenomath*}
	\begin{align}\label{equ:MPCLP} 
	\hspace{-1em}\mathrm{{WQC-MPC:}} \;\;\;\; \minimize_{\m x(t), \m u(t)} \;\;\;&   J(\m u(t))=\lambda \sum_{t=0}^{N_p-1}  {\m q^\mathrm{B}(t)}^{\top} \m u(t) \\
	\subjectto\;\;\; \;& \mathrm{{LDE}}~\eqref{equ:de-abstract1},\;\;\; \mathrm{{Constraints}}~\eqref{equ:constr-abcstract-new},  \notag
	\end{align}
%\end{linenomath*}
where $\lambda$ the unit chlorine injection cost (in $\$/\mathrm{mg}$); ${\m q^\mathrm{B}(t)}$ is a vector of the flow rates matching the dimension of $\m u(t)$. 

We note the following. First, the water quality time-step $\Delta t$, decided by Remark~\ref{rem:Deltat}, varies in different hydraulic time-steps, but during the same hydraulic time-step \blue{$T_h$}, $\Delta t$ remains the same, and \blue{$T_h =  \frac{N_p}{N_1} \Delta t$} (see Fig.~\ref{fig:Timescale}). Similarly, matrices $\m A(t)$ and $\m B(t)$ in LDE varies along with hydraulic time-scale.   Second, the initial chlorine concentration $\m x(t=0)$ is considered to be known, and as a result the optimization variables are 
%\begin{linenomath*}
	$$ \{\m x(t)\}_{t=1}^{N_p}, \;\;  \{\m u(t)\}_{t=0}^{N_p-1},$$
%\end{linenomath*}
that is, $\m u(N_p-1)$ enables the evolution of the water quality model to the constrained $\m x(N_p)$. 
{Third, there is no measurement model involved in \eqref{equ:MPCLP} which is a default assumption in the MPC scheme. That is, all states are measured using sensors. Fourth, {Constraint}~\eqref{equ:constr-abcstract-new} is a hard constraint, and potentially makes $\mathrm{WQC-MPC}$ infeasible. This issue is solved in Section~\ref{sec:scalable-WQC} via softening hard constraints. Fifth, after $\mathrm{WQC-MPC}$ is solved in the current prediction horizon $T_p$, the MPC window slides to the next one and resolve this problem (see Fig.~\ref{fig:Timescale}). The obtained control sequence $\m u(t)$ is a $n_u \times N_p$  matrix containing $N_p$ time-steps needs to be condensed and lumped into a larger time-step one to match the capacity of booster stations such as injection frequency. Additionally,} $\mathrm{WQC-MPC}$~\eqref{equ:MPCLP} is a linear program (LP), and the number of optimization variables is decided by $N_p$,  $n_x$, and $n_u$ when each node is installed with a booster station (the worst case). That is, 
%\begin{linenomath*}
	\begin{align}\label{equ:numberofvar}
	N_p (n_x + n_u) = N_p  (2n_{\mathrm{N}} + n_\mathrm{L}), 
	\end{align}
%\end{linenomath*}
where $ n_\mathrm{N} = n_{\mathrm{J}}+n_{\mathrm{R}}+n_{\mathrm{TK}}$ and $n_\mathrm{L} = n_{\mathrm{S}}  +n_{\mathrm{M}}+n_{\mathrm{V}}$.

{From~\eqref{equ:numberofvar}, we can see for any given network, the number of variables in~\eqref{equ:MPCLP} is mainly decided by  prediction horizon $T_p$ and the number of pipes $n_{\mathrm{P}} $, and the number of segments of a pipe $s_{L}$.   For example, if we have a network  with $n_{\mathrm{P}} = 15$ pipes, $ s_{L} = 500$ segments for each pipe,  and prediction horizon $T_p = 5$ minutes when $\Delta t = 1$ second or equally $N_p = 300$ time-steps, then the number of optimization variables of $\mathrm{WQC-MPC}$ for this small network would be more than $2,250,000$ which is hard to solve for most solvers even if it is an LP. Next, we introduce a plug-and-play MPC which reduces the computational burden through an analytical solution and remove the default assumption built-in MPC which is all states (including each segment in pipes) are measured via sensors.}

 \begin{figure*}
	\begin{equation}~\label{equ:uMPCdetailed}
	\resizebox{0.9\hsize}{!}{
		$\underbrace{\bmat{\m C_a \m x_a(t+1 )\\
				\m C_a \m x_a(t+2 )\\
				\vdots \\
				\m C_a \m x_a(t+N_p ) }}_{\large \m y_{p}} = \underbrace{
			\bmat{
				\m C_a  \m     \Phi_a \\
				\m C_a  \m     \Phi_a^2 \\
				\vdots \\
				\m C_a  \m     \Phi_a^{N_p}}}_{\large \m W}\m x_a(t) +
		\underbrace{\bmat{\m C_a \m \Gamma_a \\
				\m C_a \m \Phi_a\m \Gamma_a &\m C_a \m  \Gamma_a \\
				\vdots & & \ddots \\
				\m C_a \m \Phi_a^{N_p-1}\m \Gamma_a & \ldots &\m C_a \m \Phi_a\m \Gamma_a & \m C_a \m \Gamma_a}}_{\large \m Z}     \underbrace{\bmat{\Delta \m u(k) \\ \Delta \m u(t+1) \\ \vdots \\ \Delta \m u(t+N_p-1)}}_{\large \Delta \m u_{{p}}} $
	}
	\end{equation}  
	\hrulefill
\end{figure*}

\section{\large Scaling $\mathrm{WQC-MPC}$ for Real-Time Implementation}~\label{sec:scalable-WQC}
Due to the space (and time) discretization of the advection-reaction dynamics, the dimension of $\mathrm{WQC-MPC}$ explodes with the number of pipes and segments. Furthermore, the time-step for water quality simulation is in seconds (decided by L-W scheme), thereby necessitating near real-time control actions that react according to the water quality status of the water network. 
Even for small networks with tens of pipes and junctions, the LP~\eqref{equ:MPCLP} has millions of variables. Most of these variables are state variables, $\m x(t)$. In short, since the time-step is in seconds, then an instant of~\eqref{equ:MPCLP} should be solved in a fraction of a second---to make MPC implementation realizable in near real-time. 

To that end, the objective of this section is to show how to transform $\mathrm{WQC-MPC}$~\eqref{equ:MPCLP} from a constrained LP with millions of variables to either a analytical solution or a quadratic program (QP) with orders of magnitude fewer variables. The main trick here is to eliminate the dependence on the state variables  $\m x(t)$ while modeling state-constraints~\eqref{equ:constr-abcstract-new} as soft constraints in a quadratic objective function. This approach is a classical one in control-theoretic textbooks~\cite{wang2009model}; however, its adaptation to water quality control problem is the main novelty here.

\subsection{Ultra-fast, analytical solution to WQC-MPC}
First, we rewrite the difference equation as
%\begin{linenomath*}
	\begin{align}~\label{equ:de-abstract3}
	\m x(t + 1) = \m A\m x(t) + \m B \m u(t), \;\; \m y(t) = \m C \m x(t).
	\end{align}
%\end{linenomath*}

In~\eqref{equ:de-abstract3}, {we made two simplifications for ease of brevity and presenting the formation of MPC. These simplifications are $\Delta t = 1$ and $\m A(t):=\m A$ and $\m B(t):=\m B$. Note that this does not mean $\Delta t$ is 1 second, instead, it simply means a unit step.} Furthermore, we include a measurement model where $\m C$ is a binary matrix that encodes the location of water quality sensors measuring chlorine concentrations $\m y(t)$ in real-time. However, we cannot install sensors for each segment of pipes. Hence, the  concentrations in pipes can be estimated by LDE~\eqref{equ:de-abstract1}  model.
Next, we define $$\Delta \m x(t+1) = \m x(t+1) - \m x(t), \;\; \Delta \m u(t+1) = \m u(t+1) - \m u(t)$$ to be new auxiliary states and input variables that quantify corresponding  change rate. Hence, we can rewrite~\eqref{equ:de-abstract3}
as
%\begin{linenomath*}
		\begin{align*}
		{\bmat{\Delta \m x(t+1)\\ \m y(t+1)}} & =  {\bmat{\m A & \m 0 \\ \m C\m A & \m I_{n_y}}}    {\bmat{\Delta \m x(t) \\ \m y(t)}} + {\bmat{\m B\\ \m C\m B}}  \Delta \m u(t)\\
		\m y(t) & =  {\bmat{\m 0 & \m I_{n_y}}}\bmat{\m \Delta x(t) \\ \m y(t)},
		\end{align*}
%\end{linenomath*}
or more abstractly as
%\begin{linenomath*}
	\begin{subequations}~\label{equ:de-abstract5}
		\begin{align}
		\m x_a(t+1)&=\m \Phi_a\m x_a(t)+\m \Gamma_a\Delta \m u(t)\\
		\m     y(t)&=\m C_a\m x_a(t),
		\end{align}
	\end{subequations}
%\end{linenomath*}
where $  \m x_a(t) =  \{\Delta \m x(t), \m y(t)\} \in\mathbb{R}^{n_x+n_y}, \; \m \Gamma_a \in\mathbb{R}^{n_x+n_x \times n_u},\m C_a\in\mathbb{R}^{n_y \times n_x+n_y}.$
In~\eqref{equ:de-abstract5}, $\m x_a(t)$ represents the \textit{augmented state}, which will allow us in this sequel to eliminate the variables $\m x(t)$ from $\mathrm{WQC-MPC}$. As in traditional MPC schemes, we consider here that $\m x_a(t)$ is given, and as a result we compose the following equality for a prediction/MPC horizon of $N_p$ time-steps as~\eqref{equ:uMPCdetailed} 
or 
%\begin{linenomath*}
	\begin{equation}~\label{equ:uMPC}
	\m y_p = \m W \m x_a+\m Z \Delta \m u_p,
	\end{equation}
%\end{linenomath*}
where $\m y_p \in \mbb{R}^{n_y N_p}$ and $\Delta \m u_p \in \mbb{R}^{n_uN_p }$.
Equation~\eqref{equ:uMPC} essentially lumps the difference equation dynamics into a single equality constraint with a known $\m x_a$ and predetermined matrices $\m Y,\m W,$ and $\m Z$; vector $\m u_p$ is the variable to be computed/predicted over $N_p$ time-steps. 

\blue{Note that $ \m x_a(t) = {\bmat{\Delta \m x(t) \\ \m y(t)}} $ comprised of two parts can be easily estimated with the assistance of developed LDE model and the sensors installed. The first part $\Delta \m x(t)$ or equivalently $\m x(t)$ can be obtained from $\m x(t + 1) = \m A\m x(t) + \m B \m u(t)$ in LDE model~\eqref{equ:de-abstract1} as long as our LDE model is accurate and the initial value $\m x(0)$ is given; the second part $ \m y(t)$ is the sensor output. Hence, the $\m x_a$ is viewed a known in our proposed method.}

Another needed tool here is to append the objective function given in~\eqref{equ:MPCLP} with soft constraints mimicking the upper and lower bound constraints on $\m x(t)$ given in~\eqref{equ:constr-abcstract-new}. To do so, we define  $\m y^{\mathrm{ref}}$ to be a vector of a constant, predetermined reference that sets the desired chlorine concentrations to be tracked at the sensor locations. For example, $\m y^{\mathrm{ref}}$ can  be average of minimum and maximum concentrations. 

Given these developments, we can append objective functions to be written as an unconstrained quadratic optimization:
%\begin{linenomath*}
	\begin{align}~\label{equ:minquad}
	\min_{\Delta \m u_p} \; \tilde{J}(\Delta \m u_p) &= \frac{1}{2}(\m y^{\mathrm{ref}}-\m y_p)^{\top}\m Q(\m y^{\mathrm{ref}}-\m y_p) + \frac{1}{2} \Delta \m u_p^\top \m R \Delta \m u_p +\m b^{\top} \Delta \m u_p, 
	\end{align}
%\end{linenomath*}
where the appended cost function minimizes deviations from reference chlorine concentrations, smoothness of control input, and the booster chlorination operation injection cost; $\m b =   \lambda  {\m q^\mathrm{B}(t)}  \kron \m 1_{N_p} $ is akin to the cost function in $\mathrm{WQC-MPC}$, $\kron$ denotes the Kronecker tensor product and $\m 1_{N_p}$ is an $N_p$ dimensional vector of ones; matrices $\m Q=\m Q^{\top}$ and $\m R=\m R^{\top}$ are  weight matrices that dictates the preference of measurements and the smoothness of control input. {We note that matrices $\m Q$, $\m R$, and $\m b$ have impact on the final solutions, and the value of these three matrices needs to be adjusted according to the importance of the three objectives. For example, setting a relatively large $\m Q$ helps to reach the reference quickly and accurately, but the control smoothness and the cost of injected chlorine can not be guaranteed as expected.}

The only optimization variable in~\eqref{equ:minquad} is $\Delta \m u_p$ since $\m y_p$ can be substituted in~\eqref{equ:minquad} with $\m W\m x_a + \m Z \Delta \m u_p$ from~\eqref{equ:uMPC}. Upon this substitution, we obtain 
% \begin{subequations}
%	\begin{align}
%	\tilde{J}(\Delta \m u_p) &= \frac{1}{2}(\m y^{\mathrm{ref}}-\m y_p)^{\top}\m Q(\m y^{\mathrm{ref}}-\m y_p)+ \m b^{\top} \Delta \m u_p\\
%	&=\frac{1}{2}(\m y^{\mathrm{ref}}-\m W\m x_a - \m Z \Delta \m u_p)^{\top}\m Q \;\; \cdot \\
%	& \;\;\;\;\;\;\;\;(\m y^{\mathrm{ref}}-\m W\m x_a - \m Z \Delta \m u_p)+ \m b^{\top} \Delta \m u_p\notag\\
%	&    = \dfrac{1}{2} (-\m y^{\mathrm{ref}\top}\m Q\m Z \Delta \m u_p+ \\
%	& \notag  \;\;\;\;\; \;\;\;\;\; \m x_a^{\top}\m W^{\top}\m Q \m Z \Delta\m  u_p -\Delta \m u_p^{\top} \m Z^{\top}\m Q\m y^{\mathrm{ref}} +  \\
%	&\notag \;\;\;\;\; \;\;\;\;\; \Delta \m u_p^{\top} \m Z^{\top}\m Q \m W \m x_a +  \Delta \m u_p^{\top} \m Z^{\top}\m Q\m Z \Delta \m u_p)\\
%	& \notag\;\;\;\;\; + \m b^{\top}  \Delta \m u_p + \mathrm{ConstantTerms}\\
%	\tilde{J}(\Delta \m u_p)&= \dfrac{1}{2} \left( \Delta \m u_p^{\top} \m Z^{\top}\m Q\m Z \Delta \m u_p\right)+ \\
%	& \;\; \left(\m b^{\top}-\m y^{\mathrm{ref}\top}\m Q\m Z+\m x_a^{\top}\m W^{\top}\m Q \m Z  \right)\Delta \m u_p + \notag\\
%	&\;\;\;\;\;\;\; \mathrm{ConstantTerms}.\notag
%	\end{align}
%\end{subequations}
%\begin{linenomath*}
	\begin{align*}
	\tilde{J}(\Delta \m u_p) = \dfrac{1}{2}   \Delta \m u_p^{\top} (\m Z^{\top} \m Q\m Z + \m R) \Delta \m u_p   + \left(\m b^{\top}-\m y^{\mathrm{ref}\top}   \m Q\m Z+\m x_a^{\top}\m W^{\top} \m Q \m Z  \right) \Delta \m u_p  + \mathrm{constants}.
	\end{align*}
%\end{linenomath*}
Optimizing for $\Delta \m u_p$ via deriving the first order necessary conditions for unconstrained minimization through setting $\frac{\partial \tilde{J}(\Delta \m u_p)}{\partial \Delta \m u_p} = \m 0$, we obtain
%\begin{linenomath*}
	\begin{equation}~\label{equ:controllaw}
	\Delta \m u_p^* = (\m Z^{\top}\m Q\m Z + \m R)^{-1} \left( \m b- \m Z^{\top}\m Q (\m y^{\mathrm{ref}}-\m W\m x_a) \right).
	\end{equation}
%\end{linenomath*}
The optimal control law given in~\eqref{equ:controllaw} is indeed plug-and-play: it can be \textit{immediately implemented} without optimization. In particular, and after measuring $\m x_a$ at time-step $t$, the operator can compute $\Delta \m u_p^*$ and extract $\m u^*(t)$ from it, without needing to solve a large-scale optimization. The only computation needed is the large-scale matrix-vector product in~\eqref{equ:controllaw}  which can be computed efficiently using sparse matrix, rather than using iterative, interior point methods to solve the LP in $\mathrm{WQC-MPC}$. Furthermore, a large portion of the RHS of control law~\eqref{equ:controllaw}  can be computed and stored, seeing the only changing quantity---considering fixed hydraulic quantities---is $\m x_a$.   
\begin{myrem}
	\textit{Since $\m Q=\m Q^{\top}$ and $\m R=\m R^{\top}$ are positive definite weight matrices, the second-order necessary conditions for optimality of unconstrained minimization in~\eqref{equ:minquad} are also satisfied since $\m Z^{\top} \m Q \m Z + \m R$ is an invertible, positive definite matrix for all values of $\m Z$.  }
\end{myrem}

\subsection{Guaranteeing satisfaction of state {and input} bounds}
A question that begs itself here is whether the operator can guarantee that upper and lower bound constraints for chlorine concentrations {in the state vector $\m x$ and control input vector $\m u$}, given that  optimal control law~\eqref{equ:controllaw} only minimizes the deviations from the reference concentrations. Unfortunately, and although the highly-scalable control law~\eqref{equ:controllaw} has its advantages, one cannot guarantee that the constraints~\eqref{equ:constr-abcstract-new} will not be violated. Case studies investigate this point, {but we still can ease these issues by converting the constraints on states and inputs into the constraints on optimization variable $\Delta \m u_p$.}

{First, the issue of bounding chlorine concentrations in system state $\m x$} can be addressed by imposing upper and lower bound constraints on $\m y_p$, which is again reflective of the only optimization variable $\Delta \m u_p$. To examine how this works in practice, we remind readers that $\m y_p$ collects all $\m x_a(t)$ for the entire control horizon, where $\m x_a(t)$ collects $\Delta \m x(t)$ and measurements $\m y(t)$.  To that end, it is conceivable to impose constraint
$\m y_p^{\min} \leq \m y_p \leq \m y_p^{\max}$
which mimics~\eqref{equ:constr-abcstract-new}. The constraint on $\m y_p$ can be written as
$$ \bmat{-\m y_p\\ \m y_p} \leq \bmat{
	-\m y_p^{\min} \\ \m y_p^{\max}
}.$$
Recalling that $\m y_p = \m W\m x_a + \m Z\Delta \m u_p$, we obtain
%\begin{linenomath*}
	\begin{align}~\label{equ:qpconstraints}
	\bmat{-\m Z\\ \m Z}\Delta \m u_p \leq \bmat{
		-\m y_p^{\min} +\m W\m x_a\\ \m y_p^{\max}-\m W\m x_a}.
	\end{align}
%\end{linenomath*}
 
Second, the bounds on control input $\m u$ can be also solved similarly as the way of bounding system state $\m x$. Recalling that $u(t) = u(t-1) + \Delta u(t) = u(t-1) + \bmat{\m I_{n_u} & \m O & \ldots \m O} \Delta \m u_p$ and $u(t+1) = u(t) + \Delta u(t+1) = u(t) + \bmat{\m I_{n_u} & \m I_{n_u} & \ldots \m O} \Delta \m u_p$, we obtain
%\begin{linenomath*}
	\begin{equation*}
	\resizebox{0.7\hsize}{!}{
		$\underbrace{\bmat{\m u(t)\\\m u(t+1)\\ \vdots \\ \m u(t+N_p -1 ) }}_{\large \m u_p} = \underbrace{
			\bmat{
				\m I_{n_u} \\
				\m I_{n_u} \\
				\vdots \\
				\m  I_{n_u}
			}
		}_{\large \m H_1}  u(t-1)+
		\underbrace{\bmat{\m I_{n_u} \\
				\m I_{n_u} &\m I_{n_u} \\
				\vdots & & \ddots \\
				\m I_{n_u} & \ldots &\m I_{n_u} &\m I_{n_u} }}_{\large \m H_2}     \underbrace{\bmat{\Delta \m u(k) \\ \Delta \m u(t+1) \\ \vdots \\ \Delta \m u(t+N_p-1)}}_{\large \Delta \m u_{{p}}}.$
	}
	\end{equation*}  
%\end{linenomath*}

Suppose that we have control input constraints $\m u_p^\mathrm{min} \leq \m u_p \leq \m u_p^\mathrm{max}$, and it can be rewritten  as
%\begin{linenomath*}
	\begin{align}~\label{equ:inputconstraints}
	\bmat{-\m H_2\\ \m H_2}\Delta \m u_p \leq \bmat{
		-\m u_p^{\min} +\m H_1 u(t-1)\\ \m u_p^{\max}-\m H_1 u(t-1)}.
	\end{align}
%\end{linenomath*}

We observe that all bounds on system state $\m x$ and control input $\m u$ are expressed as the constraints on optimization variable $\Delta \m u_p$ [see~\eqref{equ:qpconstraints} and~\eqref{equ:inputconstraints}], and these inequality constraints \normalcolor can then be incorporated into a QP formulation for WQC, That is,
%\begin{linenomath*}
	\begin{eqnarray}
	\minimize_{\Delta \m u_p} & & J(\Delta \m u_p) = \frac{1}{2} \Delta \m u_p^\top \m R \Delta \m u_p 
	+ \m b^{\top} \Delta \m u_p \label{equ:mpcscalable-2} \\ % \frac{1}{2} \Delta \m u_p^\top \m R \Delta \m u_p +
	\subjectto & &  ~\eqref{equ:qpconstraints},{~\eqref{equ:inputconstraints} }\notag.
	\end{eqnarray}     
%\end{linenomath*}   

Note that unlike $\mathrm{WQC-MPC}$~\eqref{equ:constr-abcstract-new}, optimization problem~\eqref{equ:mpcscalable-2} has far less optimization variables and constraints, seeing that $\m x(t)$, the vector with the majority of optimization variables in~\eqref{equ:constr-abcstract-new}, is not present in~\eqref{equ:mpcscalable-2}. The number of optimization variables of~\eqref{equ:minquad} and~\eqref{equ:mpcscalable-2} are 
%\begin{linenomath*}
	\begin{align}\label{equ:numberofvar2}
	N_p n_u = N_p n_{\mathrm{N}},
	\end{align}
%\end{linenomath*}
and the percentage reduction is $1- \frac{N_p n_{\mathrm{N}}}{N_p  (2n_{\mathrm{N}} + n_\mathrm{L})} = \frac{n_{\mathrm{N}} + n_\mathrm{L}}{2n_{\mathrm{N}} + n_\mathrm{L}}$; see specific examples in Tab.~\ref{tab:scale} in Section~\ref{sec:test} for details. 
\normalcolor

\subsection{Real-time implementation}
\begin{figure*}[t]
	\centering
	\includegraphics[width=0.7\linewidth]{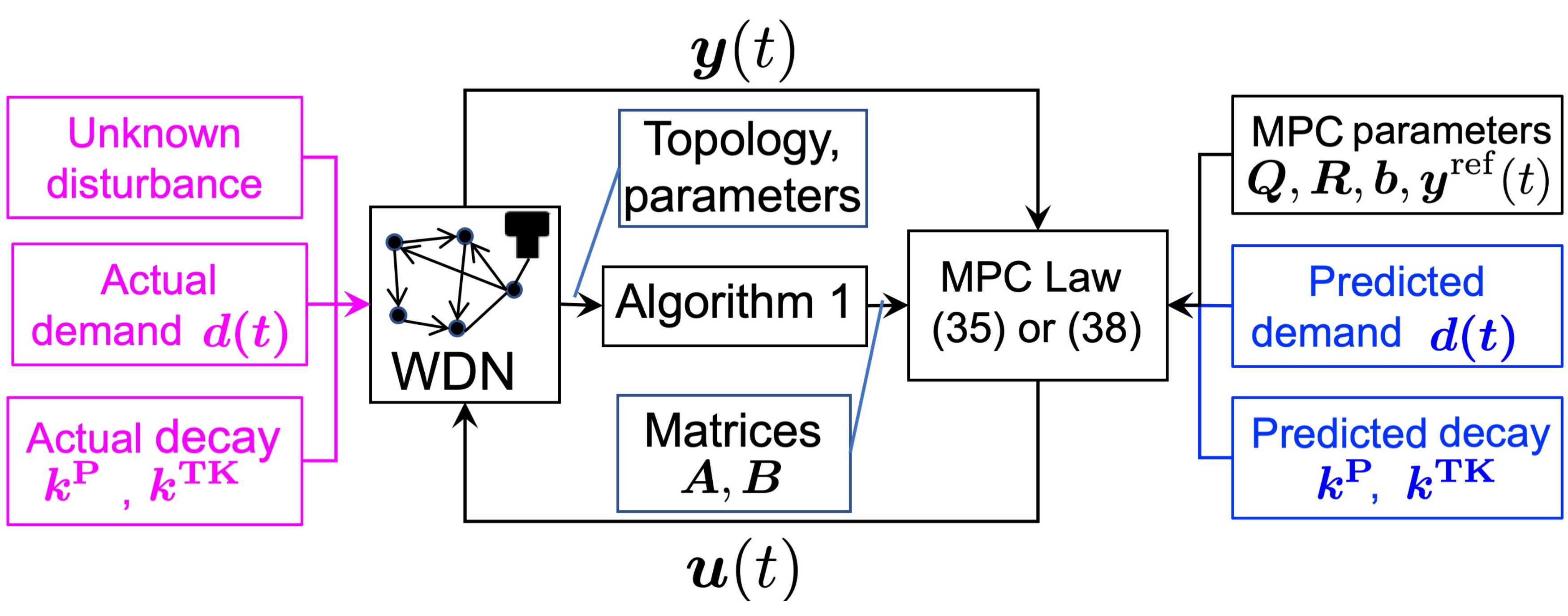}
	\caption{Control diagram for WQC problem. }
	\label{fig:control}
\end{figure*}
This section summarizes how the proposed computational methods can be implemented in real-world WQC applications under uncertainty. The approach is depicted in Fig.~\ref{fig:control}. The blue boxes showcase quantities that are given ahead of time while the red boxes are real-time, actual quantities that are \textit{not} furnished for the control algorithm, considering that perfect parameters for the MPC controller is unrealistic. Our approach considers three sources of uncertainty~\cite{basupi2019flexible,boosteruncertainty} for extended period simulation: \textit{(i)} unknown disturbance (i.e., the chlorine concentration drops suddenly and unexpectedly due to the reaction with unknown intrusive chemicals--sodium arsenite or contamination events happen), \textit{(ii)} demand uncertainty, and \textit{(iii)} parameter uncertainty (i.e., estimated decay coefficients of chlorine in pipes and tanks $\m k^{\mathrm{P}}$ and $\m k^{\mathrm{TK}}$ is not accurate). The details and test results are in Section~\ref{sec:test}.

First, the WDN topology and parameters are provided to Algorithm~\ref{alg:matrix} which is the first basic step in obtaining the state-space difference equation model. The WDN parameters include predictions of demand, and hence predictions of heads and flows for the entire network for an extended period of simulation (e.g., for $T_d = 24$ hours). {This indicates our proposed method is based on the solution of steady-state hydraulics}. Second, the state-space matrices $\m A(t)$ and $\m B(t)$ are computed for the entire simulation time via the derivations in~\eqref{equ:xVector}--\eqref{equ:de-abstract1} that result in the state-space model~\eqref{equ:fullMatrix}. These computations, albeit require little computational overhead, are performed offline and hence can be stored to be used sequentially for the control law.  

Third, and depending on the operator preference of either applying an analytical (more computationally tractable) control law as in~\eqref{equ:controllaw} or the control law via the constrained quadratic program~\eqref{equ:mpcscalable-2}, the control law is generated by extracting the first $m$ elements of vector $\Delta \m u_p(t)$ resulting in $\Delta \m u(t)$. {This is then followed by computing $\m u(t) = \Delta \m u(t) - \m u(t-1)$, where $\m u(t-1)$ can be initialized as zero vector.} We note that this control law, $\m u(t)$, is indeed a function of the current augmented state of the network, namely $\m x_a(t)$ which includes the sensor measurements $\m y(t)$.  Fourth, this control law is then used to perform water quality control via applying it to the WDN simulation---this simulation considers real-time quantities and parameters that are not provided for the MPC law as described above. Finally, this procedure is repeated after the MPC law obtains new measurements $\m y(t)$, thereby resulting in a closed-loop, feedback-driven real-time implementation. Case studies are given in the next section.

\section{\large Case Studies}~\label{sec:test}

We present three simulation examples (illustrative three-node, Net1, and Net3 networks~\cite{rossman2000epanet,shang2002particle}) to illustrate the applicability of our approach. The three-node network is designed to illustrate the details of the proposed method,  then we test the Net1 network with looped topology, and the final Net3 is for scalability test. All test cases are simulated via EPANET Matlab Toolkit~\cite{Eliades2016} on Windows 10 Enterprise with an Intel(R) Xeon(R) CPU E5-1620 v3 @3.50 GHz.  
All codes, parameters and tested networks are available in open data repository~\cite{wang_2020}. 

\subsection{Three-node network}\label{sec:3-node}
\begin{figure}[t]
	\centering
	\includegraphics[width=0.851\linewidth]{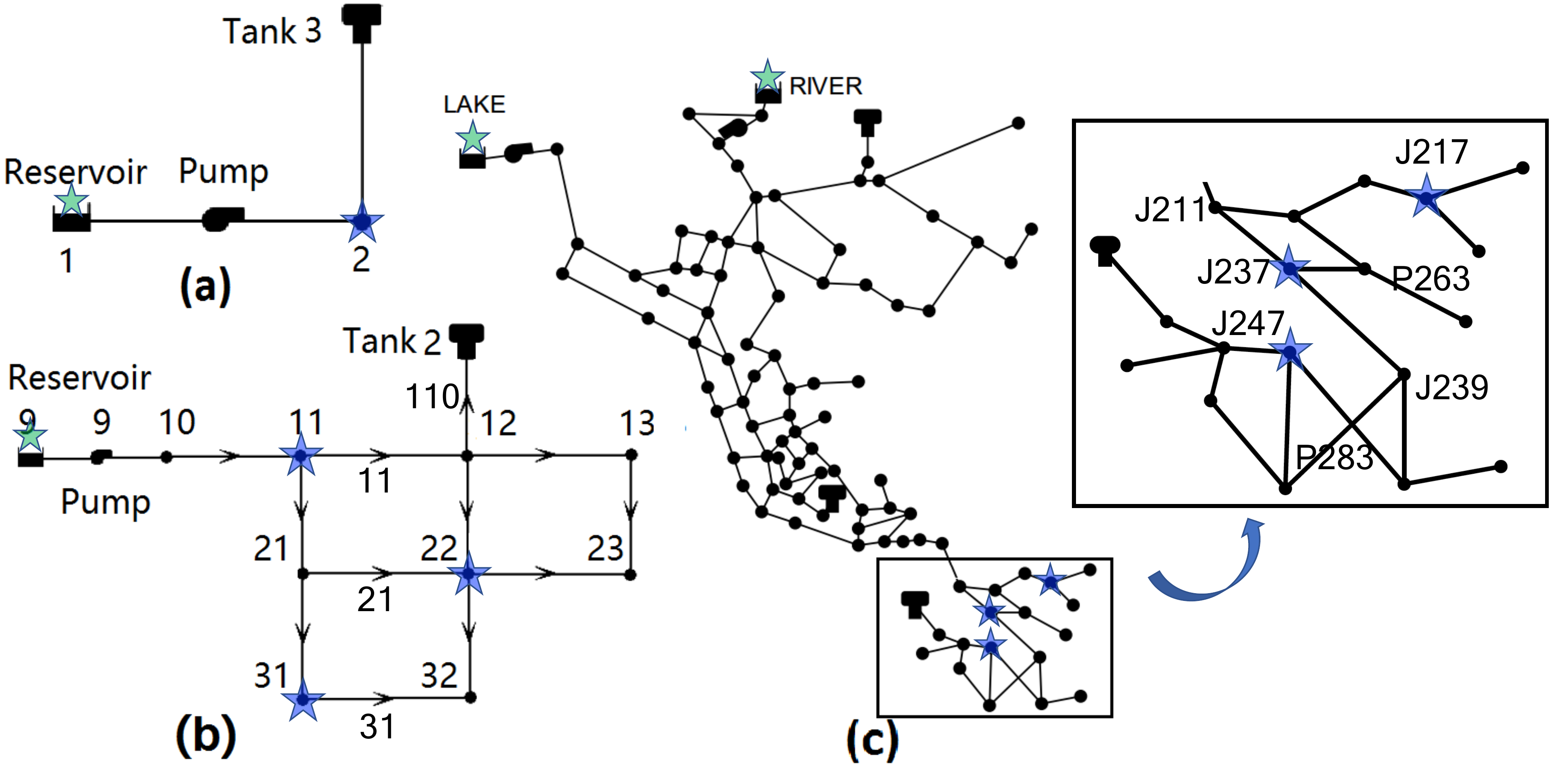}
	\caption{{Tested networks installed with mass boosters  (marked as blue stars) and concentration sources (marked as green stars): (a) Three-node network, (b) Net1, and (c) Net3 and its zoomed-in control area.}}
	\label{fig:setup}
\end{figure}

Three-node network shown in Fig.~\ref{fig:setup}a includes one junction, one pipe, one pump, one tank, and one reservoir. Only Junction 2 consumes water, and the difference between actual and predicted demand is a random number in $[-10\%, 10\%]$, which is viewed as demand uncertainty. A chlorine concentration source ($ c^\mathrm{R}_1 = 0.8$ mg/L) is installed at Reservoir 1 (marked as a green star in Fig.~\ref{fig:setup}a) and provides fixed $0.8$ milligram per liter chlorine concentration, and the initial chlorine concentrations at or in the other components are $0$ mg/L. The pipe is split into $s_{L_{23}} = 100$ segments. At first, we need to verify the effectiveness of proposed LDE model~\eqref{equ:de-abstract1}.
%\begin{figure}[t]
%	\centering
%	\subfloat[\label{fig:setup_a}]{\includegraphics[keepaspectratio=true,scale=0.4]{3node.pdf}}{} \hspace{1em}
%	\subfloat[\label{fig:setup_b}]{\includegraphics[keepaspectratio=true,scale=0.3]{Demand1day_3node.pdf}}{}\vspace{-0.1cm}\hspace{-0.0cm}\vspace{-0.02cm}
%	\caption{(a) Three-node network with a booster station is installed at Junction 2 (marked as a star), (b)  its predicted and actual demand pattern.}
%	\label{fig:setup}
%\end{figure}

\begin{figure}[t]
	\centering
	\subfloat[\label{fig:Error_EPANET_LDE}]{\includegraphics[keepaspectratio=true,scale=0.3]{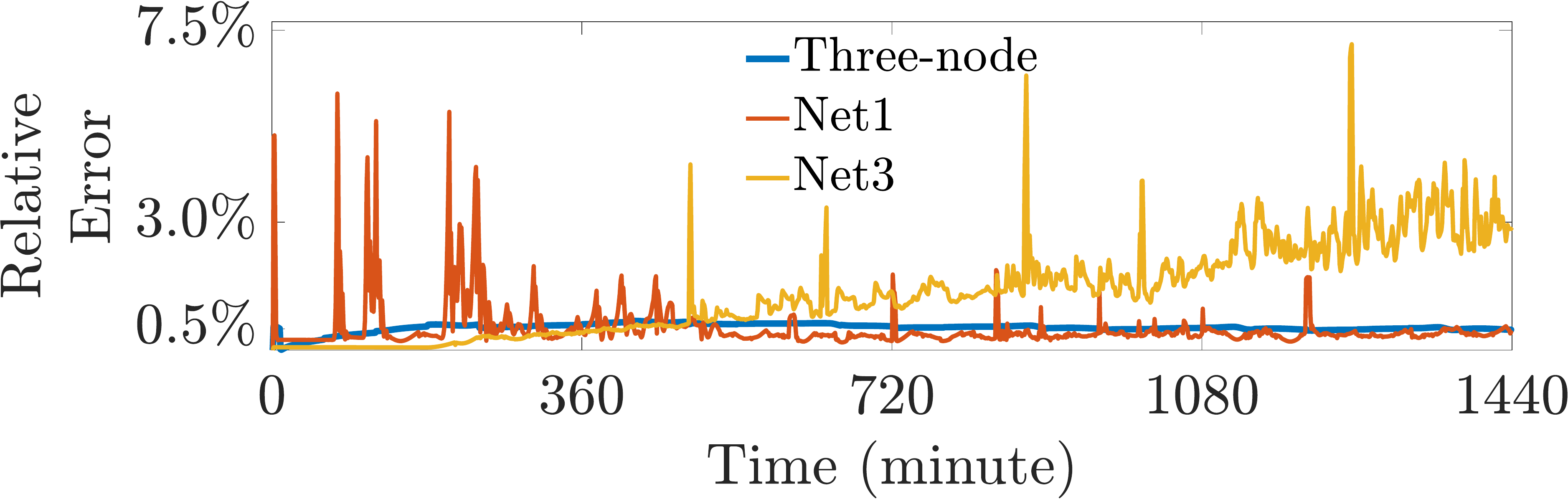}}{}\hspace{0.62cm}
	\subfloat[\label{fig:P21}]{\includegraphics[keepaspectratio=true,scale=0.35]{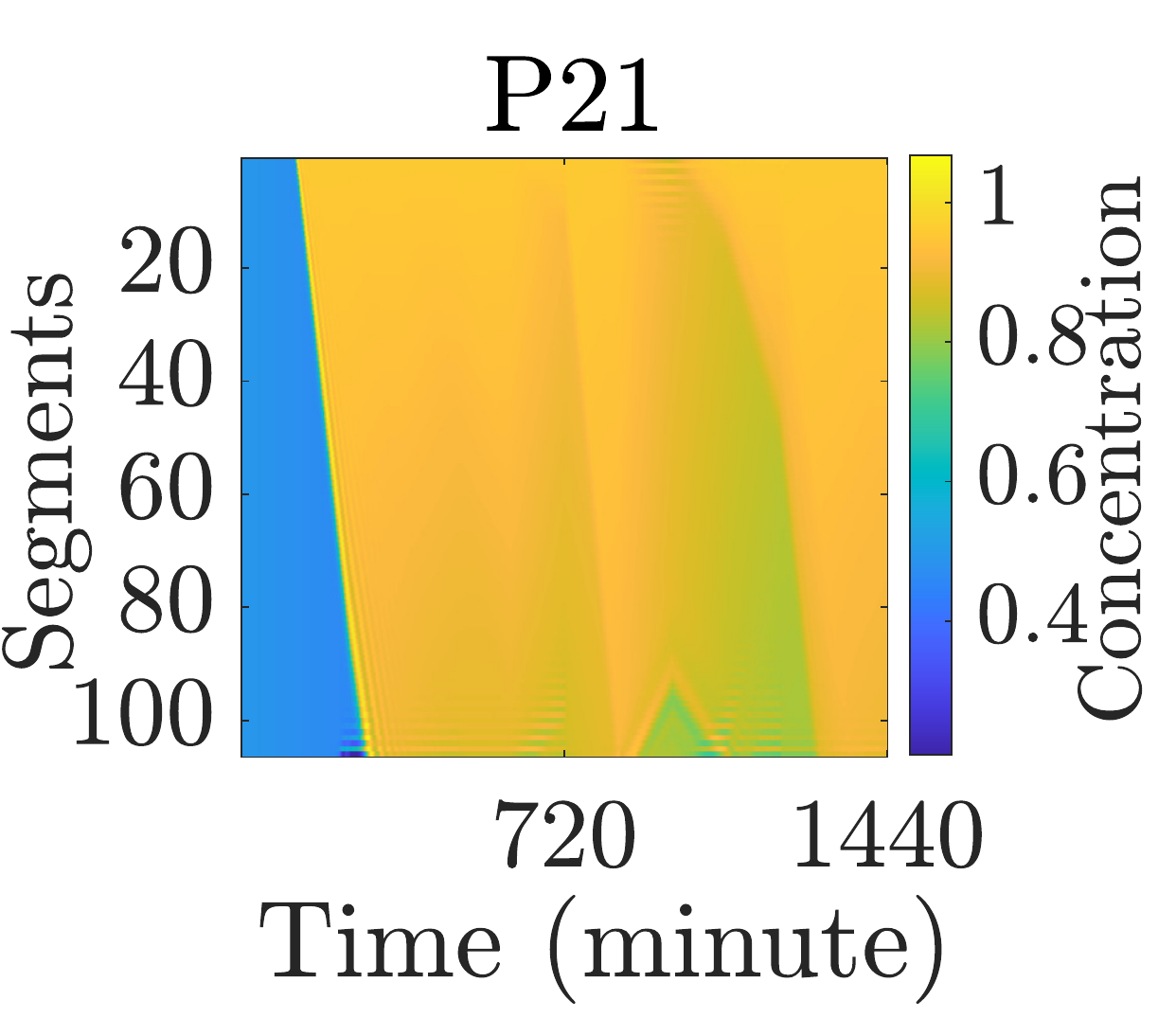}}{}  \vspace{-0.1cm}
	\caption{{(a) relative error between EPANET and LDE, (b) time-space discretization results of Pipe $21$ in Net1 based on the L-W scheme.}}
	\label{fig:timespaceResult}
\end{figure}

%\begin{figure}[t]
%	\centering
%	\includegraphics[width=0.6\linewidth]{Error_EPANET_LDE.eps}
%	\caption{{Relative error between EPANET and LDE.}}
%	\label{fig:Error_EPANET_LDE}% THIS figure is in the LDE branch on github.
%\end{figure} 
\subsubsection{Accuracy of the proposed LDE model}\label{sec:eff-3-node}
We need to verify the correctness and usability of the LDE model via simulation software package EPANET. To this end, we keep settings and parameters exactly the same in LDE and EPANET, such as keeping {a concentration source} $ c^\mathrm{R}_1 = 0.8$ mg/L, the same demands. The results are shown in Fig.~\ref{fig:Error_EPANET_LDE}, and we can see the relative error between LDE model and EPANET hovers around 0.5\% at steady-state, but exhibits larger oscillations for the bigger network. Regardless, this model is only used to inform the control law---the control law is then applied to the simulation toolbox.
% even if there are no measurements from EPANET to LDE \comment{this "no measurement" is absolutely confusing. I know what you mean, but the reader won't. }. But for a practical system, according to Assumption~\ref{asmp:meauremnet}, measurements would be sent from actual WDN (using EPANET here) to LDE to adjust the model, which means the relative error can be reset after a period time. \comment{Again, this is utterly confusing unfortunately.}

We only show a specific example of LDE model~\eqref{equ:de-abstract1}  with the number of pipe segments is $s_{L_{23}} = 3$. When $\Delta t = 1$ second, $t=0:300$ seconds {which is the first prediction horizon}, %see Remark~\ref{rem:AB}, 
and $\m x(t+\Delta  t)$ can be expressed by
%\begin{linenomath*}
	\begin{equation} 
	\hspace{-0.7em}\includegraphics[width=0.65\linewidth,valign=c]{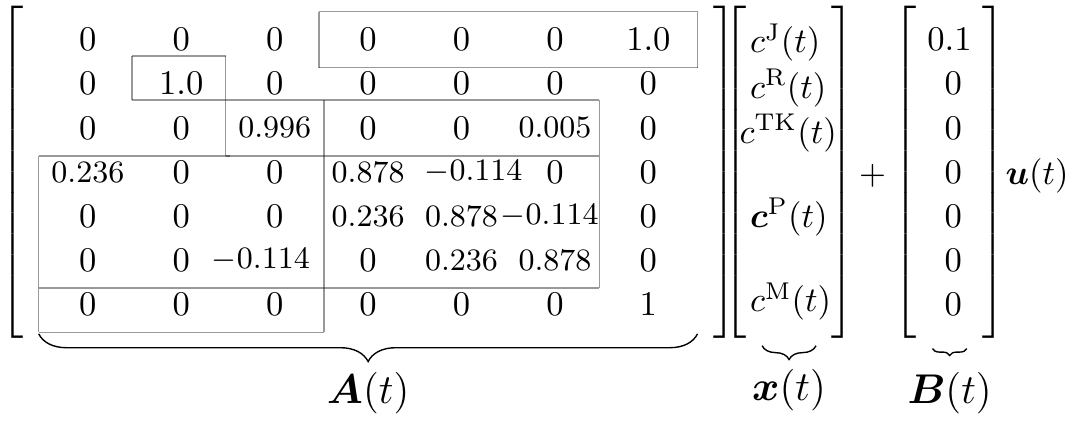}\notag
	\end{equation}
%\end{linenomath*}
where $\m c^\mathrm{P}(t)	\triangleq \{ c^\mathrm{P}(1,t), c^\mathrm{P}(2,t), c^\mathrm{P}(3,t) \}$,
and $\m u(t) \triangleq \{c^\mathrm{B}_\mathrm{J}(t+\Delta  t) \}$ are known. Next, we apply the MPC algorithm to implement water quality control considering uncertainty. 
%$\left(\begin{array}{ccccccc} 0 & 0 & 0 & 0 & 0 & 0 & 1.0\\ 0 & 1.0 & 0 & 0 & 0 & 0 & 0\\ 0 & 0 & 0.9955 & 0 & 0 & 0.004526 & 0\\ 0.2356 & 0 & 0 & 0.878 & -0.1136 & 0 & 0\\ 0 & 0 & 0 & 0.2356 & 0.878 & -0.1136 & 0\\ 0 & 0 & -0.1136 & 0 & 0.2356 & 0.878 & 0\\ 0.5 & 0.5 & 0 & 0 & 0 & 0 & 0 \end{array}\right)$
%
%$\left(\begin{array}{ccc} 0.1082 & 0 & 0\\ 0 & 0 & 0\\ 0 & 0 & 0\\ 0 & 0 & 0\\ 0 & 0 & 0\\ 0 & 0 & 0\\ 0 & 0 & 0 \end{array}\right)$
%
%$\begin{bmatrix}
%0 & 0 & 0 & 0 & 0 & 0 & 1.0\\ 0 & 1.0 & 0 & 0 & 0 & 0 & 0\\ 0 & 0 & 0.9955 & 0 & 0 & 0.004526 & 0\\ 0.2356 & 0 & 0 & 0.878 & -0.1136 & 0 & 0\\ 0 & 0 & 0 & 0.2356 & 0.878 & -0.1136 & 0\\ 0 & 0 & -0.1136 & 0 & 0.2356 & 0.878 & 0\\ 0.5 & 0.5 & 0 & 0 & 0 & 0 & 0
%\end{bmatrix}$
\subsubsection{WQC-MPC for three-node network}\label{sec:mpc-3-node} The MPC controller or \eqref{equ:controllaw} is used to decide the minimum dosage of chlorine injected by booster stations to maintain the proper concentration $\m y^{\mathrm{ref}}$ in the network. The reference or desired value vector in~\eqref{equ:de-abstract3} is  $\m y^{\mathrm{ref}} = 2$ mg/L, and the price $\lambda = 0.001 \$/\mathrm{mg}$.  The  $\m Q$ and $\m R$ in MPC~\eqref{equ:controllaw} are identity matrices.
{The parameters for this simulation are $T_d = 24$ hours ($1440$ minutes), $T_n = T_h = 1$ hour, and $T_p = 5$ minutes. Note that $s_{L_{23}} = 100$, and the water quality time-step $\Delta t$ is dynamically decided by Remark~\ref{rem:Deltat}.}
%It implies after solving hydraulic problem considering demand pattern Fig.~\ref{fig:setup_b} each hour, the~\eqref{equ:controllaw} are solved to obtain the control actions $\m u(t)$ and . 

\begin{figure}[t]
	\centering
	\subfloat[\label{fig:concentration_a}]{\includegraphics[keepaspectratio=true,scale=0.21]{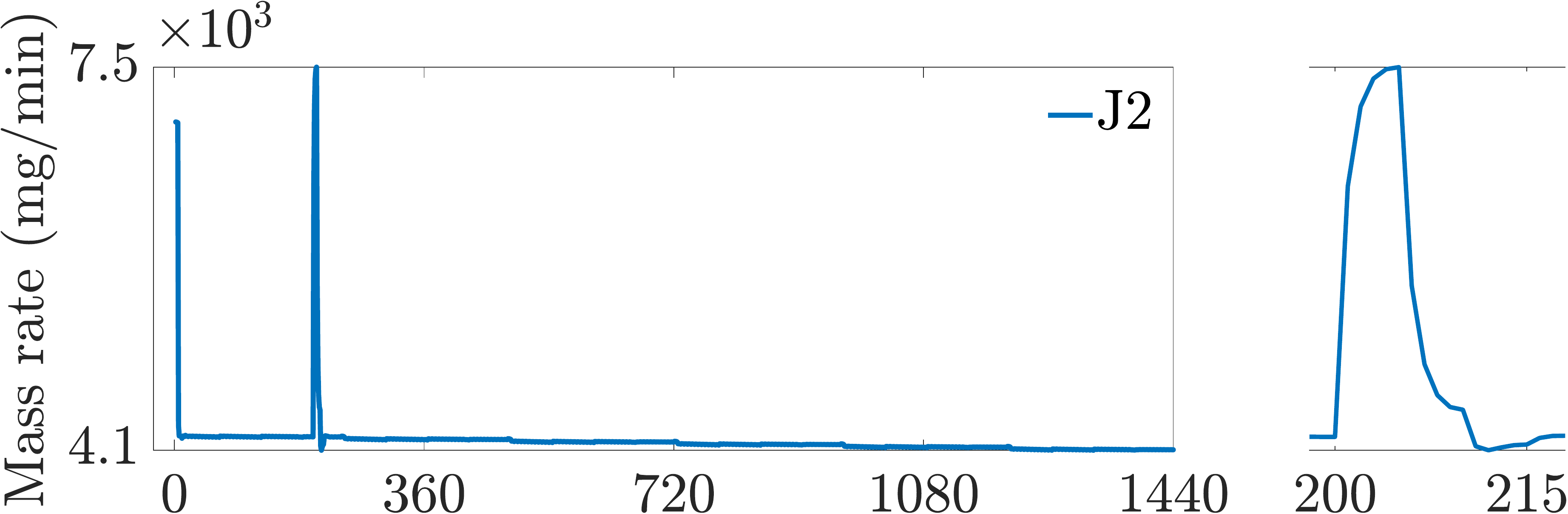}}{}  \hspace{1cm}
	\subfloat[\label{fig:concentration_b}]{\includegraphics[keepaspectratio=true,scale=0.21]{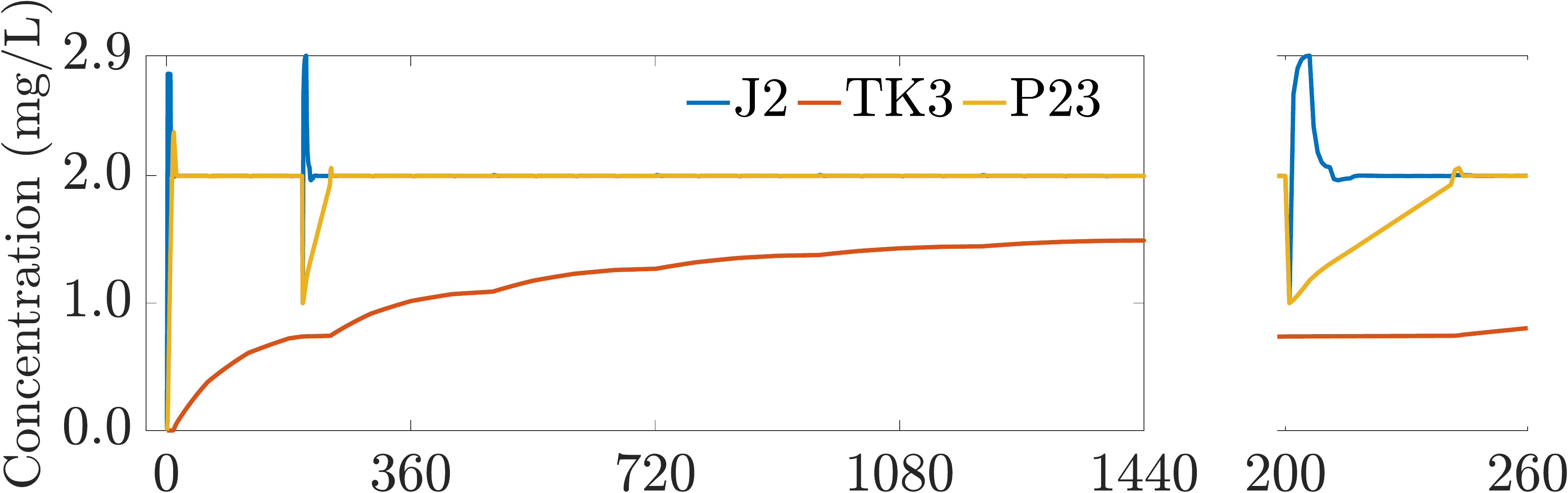}}{}\vspace{-0.1cm}
	\caption{(a) MPC control action $\m u$ during 1440 minutes for three-node network, (b) chlorine concentration at Junction 2, Tank 3, and Pipe 23 under three uncertainty sources, and unknown disturbance happens at 200-$\mathrm{th}$ minute.}
	\label{fig:concentration}
\end{figure}

Our simulation considers three sources of uncertainty~\cite{basupi2019flexible,boosteruncertainty}  for multi-period simulation. In this case, unknown disturbance can be simulated by hijacking or changing the concentrations at specific nodes or links suddenly. That is, concentration at Junction 2 and Pipe 23 at 200-$\mathrm{th}$ minute is forced to $1$ mg/L suddenly shown in Fig.~\ref{fig:concentration_b}.   The difference between actual and predicted demand is a random number in $[-10\%, 10\%]$,  which is the demand uncertainty. Similarly,  we add 10\% uncertainty to $k_{23}^b$ and $k^{w}_{23}$ in~\eqref{equ:reaction} which is reaction rate uncertainty.

The final control action is shown in Fig.~\ref{fig:concentration_a}, and the corresponding control effects under uncertainty are presented in Fig.~\ref{fig:concentration_b}. The initial chlorine concentration at Junction 2 and Pipe 23 are zeros, and in order to reach the reference value, which is $2.0$ mg/L as quickly as possible, MPC law injects {$6994.9$} mg chlorine at the very first minute. This control action results in overshoot at Junction 2, and the overshoot value is {$2.8$} mg/L as shown in Fig.~\ref{fig:concentration_b}. We can increase $\m R$ in~\eqref{equ:controllaw} to avoid overshoots, but it is not necessary here since {$2.8$} mg/L is still in $[0.2, 4]$ mg/L and acceptable. 

The mass rate injected by booster station drops to  $4288.3$ mg/minute five minutes later in Fig.~\ref{fig:concentration_a}, and the corresponding concentrations at Junction 2 and in Pipe 23 reach desired  $2.0$ mg/L in Fig.~\ref{fig:concentration_b}. Seeing that this result for 1440 minutes is obtained three sources of uncertainty, and they are properly handled by MPC controller since only tiny oscillation around $2.0$ mg/L is observed in Fig.~\ref{fig:concentration_b}.  

As for unknown disturbance uncertainty, we make $c^\mathrm{J}_{2}$ and $\m c^\mathrm{P}_{23}$ drop to $1.0$ mg/L at at the 200-$\mathrm{th}$ minute, and see the magnified plot in right part of Fig.~\ref{fig:concentration_b}.  MPC controller responds immediately and injects 7500 mg to improve the concentration level. It takes 15 minutes for  $c^\mathrm{J}_{2}$ to recover back to the desired value, but for the $\m c^\mathrm{P}_{23}$, this process takes 50 minutes.

%\begin{figure}[t]
%	\centering
%	\includegraphics[width=0.4\linewidth]{net1.pdf}
%	\caption{Net1 network scheme and its booster station locations (marked as stars).}
%	\label{fig:net1}
%\end{figure}

%\begin{figure}[t]
%	\centering
%	\subfloat[\label{fig:net1setup_a}]{\includegraphics[keepaspectratio=true,scale=0.2]{Demand4days}}{}
%	\subfloat[\label{fig:net1setup_b}]{\includegraphics[keepaspectratio=true,scale=0.2]{InterestedPipeFlowRate}}{}\vspace{-0.1cm}\hspace{-0.0cm}\vspace{-0.02cm}
%	\caption{(a) Predicted and actual demand, and (b) pipe flow rates.}
%	\label{fig:net1setup}
%\end{figure}

\begin{figure}[t]
	\centering
	\includegraphics[width=0.5\linewidth]{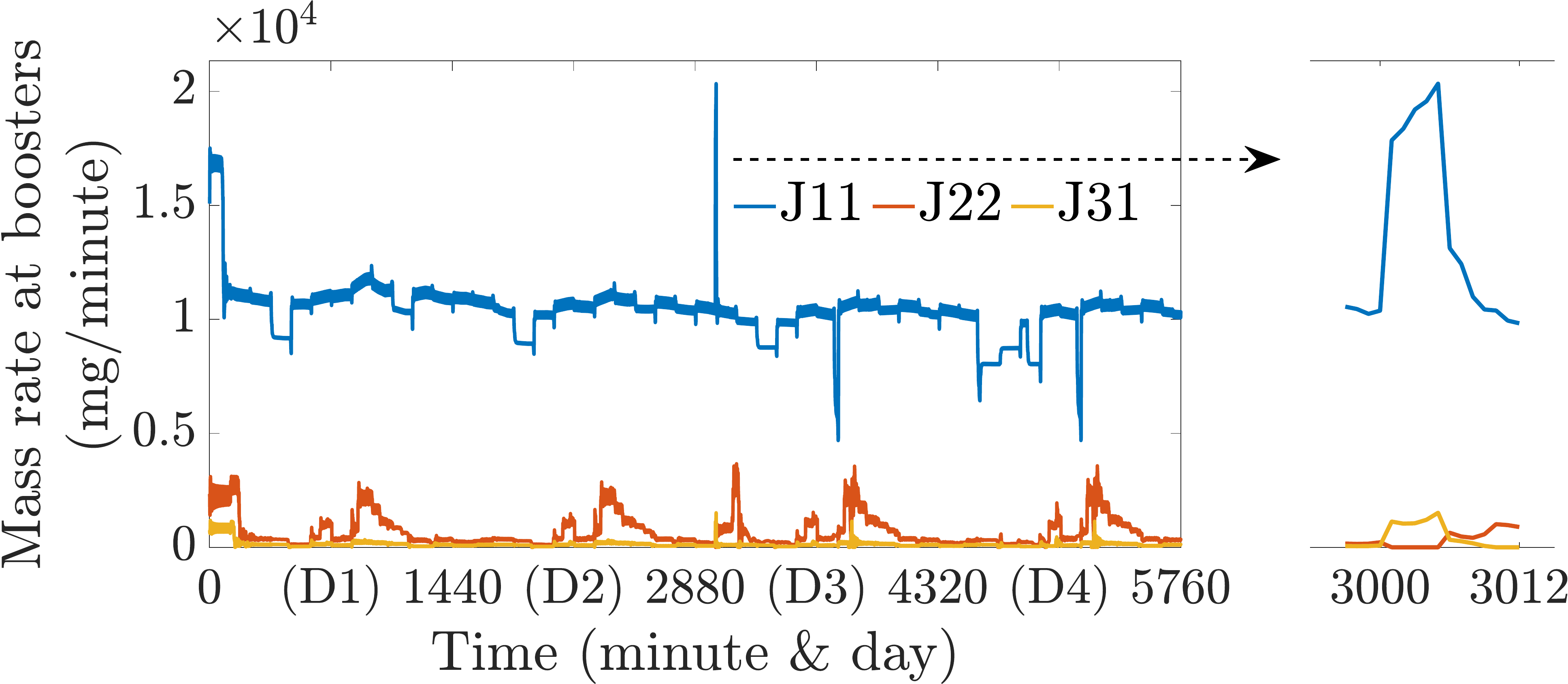}
	\caption{{Control action $\m u$ at three mass booster locations in Net1 (4 days).}}
	\label{fig:ControlActionU_net1}
\end{figure}

\begin{figure}[t]
	\centering
	\subfloat[\label{fig:concentration2_a}]{\includegraphics[keepaspectratio=true,scale=0.21]{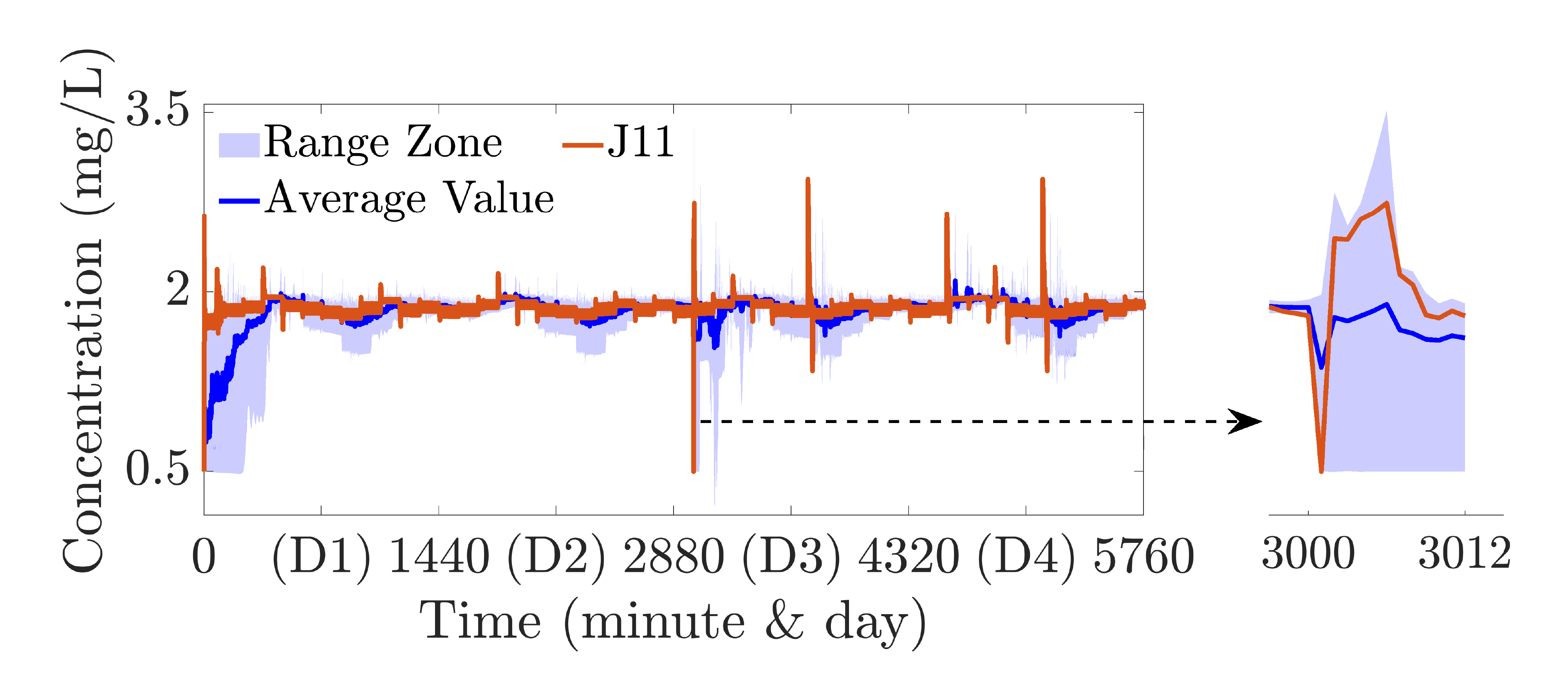}}{}\vspace{-0.1cm} \hspace{0.5cm}
	\subfloat[\label{fig:concentration2_b}]{\includegraphics[keepaspectratio=true,scale=0.21]{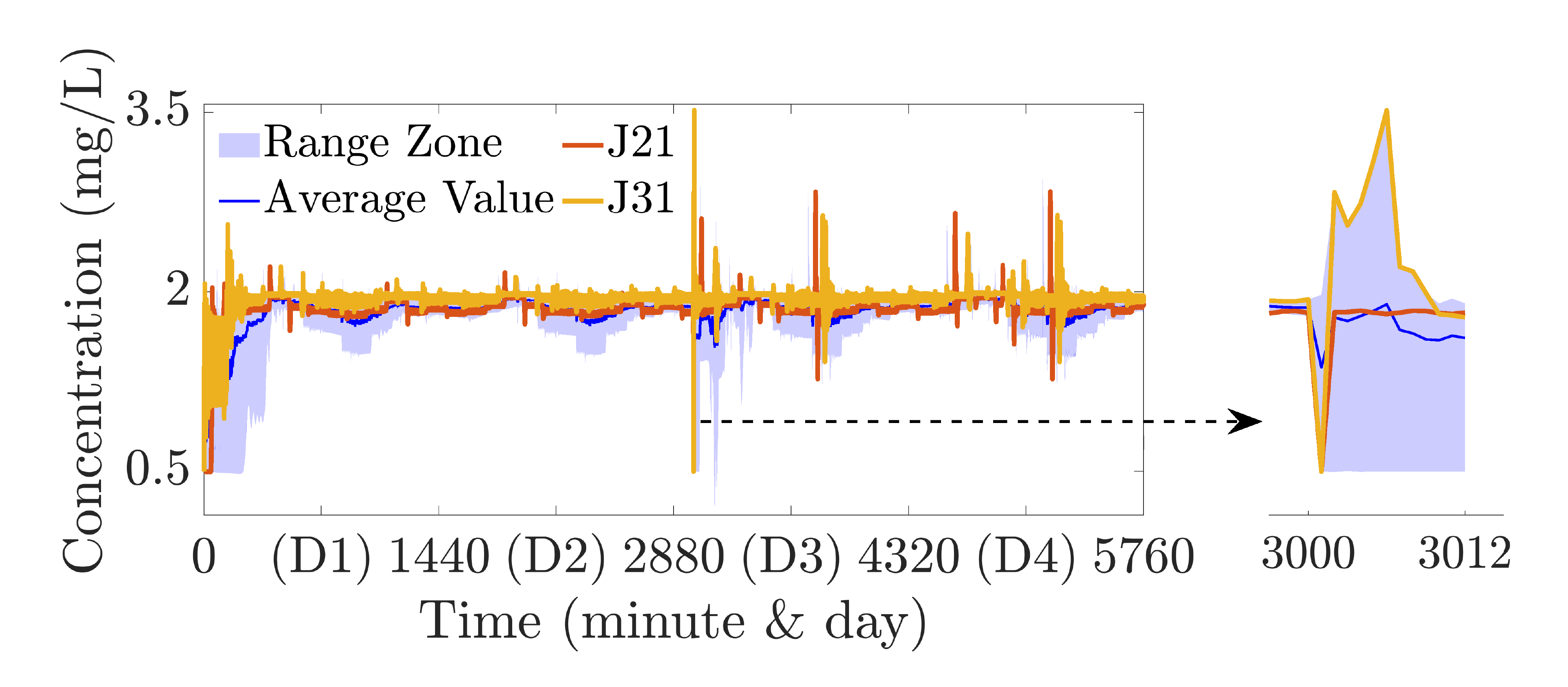}}{}\vspace{-0.1cm}\hspace{-0.2cm}
	\caption{Chlorine concentration at Junction 11 (a), Junctions 21 and 31 (b) for 4 days under three uncertainty sources (unknown disturbance happens at 3000-$\mathrm{th}$ minute and it takes 12 minutes to recover).}
	\label{fig:concentration2}
\end{figure} 

\subsection{Net1 network}\label{sec:net1} 
Net1 network~\cite{rossman2000epanet,shang2002particle}  shown in Fig.~\ref{fig:setup}b is composed of 9 junctions, 1 reservoir, 1 tank, 12 pipes, and 1 pump. 
The effectiveness of proposed LDE model is verified and shown in Fig.~\ref{fig:Error_EPANET_LDE}. {The maximum error is 7\%, that is due to the ``overshooting" (drawback) of L-W scheme. But when the L-W scheme is stable, the relative error is in [0.5\%, 1\%]. Additionally,  the time-space discretization  results (with 100 segments and 24 hours) of Pipe $21$ in Net1 based on the L-W scheme is presented as Fig.~\ref{fig:P21}.}

\begin{figure}[t]
	\centering
	\subfloat[\label{fig:concentration3_a}]{\includegraphics[keepaspectratio=true,scale=0.22]{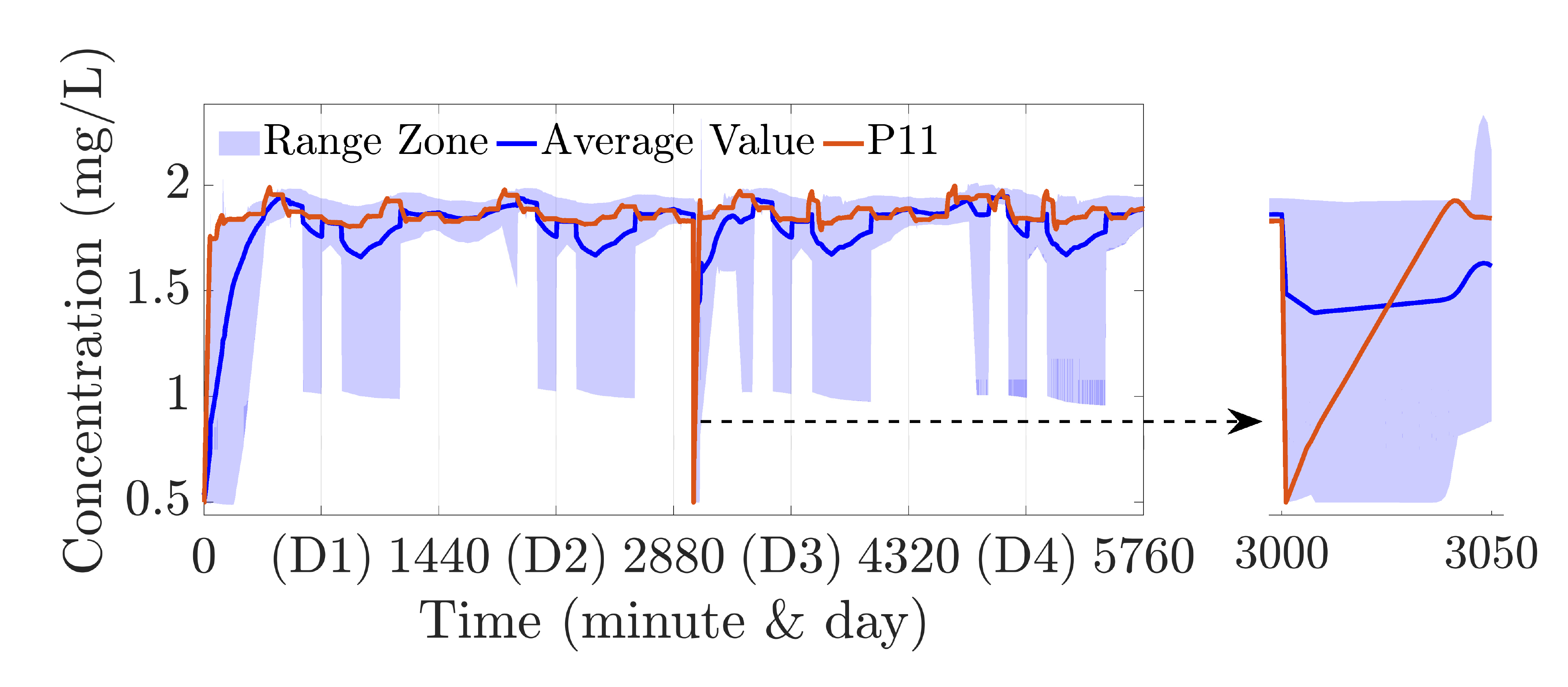}}{} 
%	\vspace{-0.1cm} 
	\subfloat[\label{fig:concentration3_b}]{\includegraphics[keepaspectratio=true,scale=0.22]{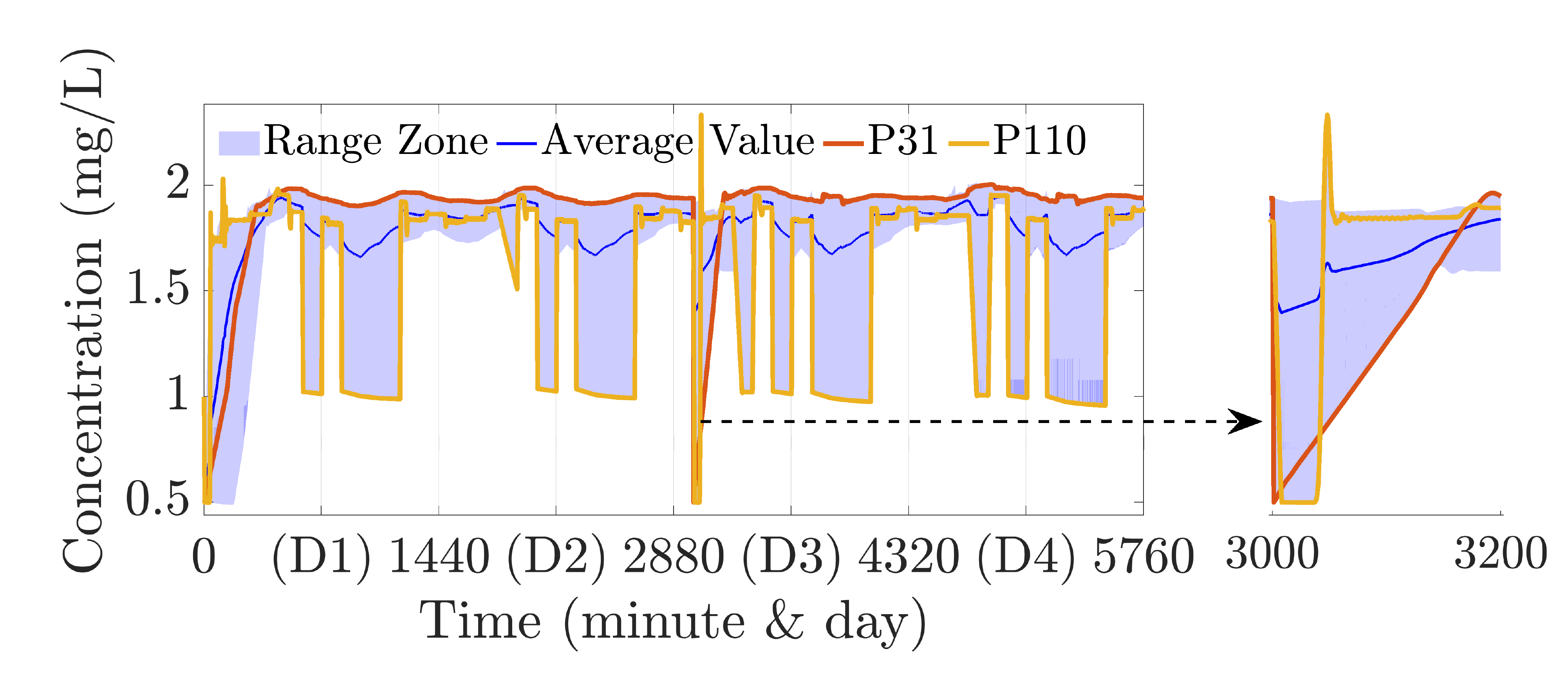}}{} \hspace{-0.2cm}
%	\vspace{-0.2cm}
	\caption{Chlorine concentration in Pipe 11 (a), Pipes 31 and 110 (b) for 4 days under three uncertainty sources (unknown disturbance happens at 3000-$\mathrm{th}$ minute and it takes 50 and 200 minutes for Pipe 11 and Pipe 110 to recover).}
	\label{fig:concentration3}
\end{figure}

{Actuator placement problem is another topic that is out the scope of this paper, and the location of booster stations are assumed to be predetermined subjectively, and fixed at Junctions 11, 22, and 31, and there is a chlorine source ($1.0$ mg/L) at Reservoir 9 (marked as a green star in Fig.~\ref{fig:setup}b).} The initial chlorine concentrations at junctions and Tank 2 are $0.5$ mg/L and $1.0$ mg/L. {The simulation duration $T_d = 4$ days ($5760$ minutes), demand changes every $T_n = 2$ hour, hydraulic step is $T_h = 1$ hour, and $T_p = 5$ minutes. We apply  the same uncertainty and MPC parameters} as the ones in three-node network. The final control action is shown in Fig.~\ref{fig:ControlActionU_net1}, and the corresponding control effects are presented in Figs.~\ref{fig:concentration2} and~\ref{fig:concentration3}. The main control scheme and results are similar to three-node network. Hence, they are not repeated. Instead, the interesting and meaningful observations and differences are presented next.

First, multi-booster stations are installed, and we can see that the booster at Junction 11 shoulders the majority of burden compared with the other two boosters from Fig.~\ref{fig:ControlActionU_net1}. This result can be explained intuitively from Fig.~\ref{fig:setup} since \textit{(i)} the booster station installed at Junction 11 covers the majority of nodes except Junction 10. In fact, the concentration at Junction 10 can not be adjusted in this case; \textit{(ii)} the booster station installed at Junction 22 covers Junctions 22, 23, and 32; and \textit{(iii)} the booster station installed at Junction 31 only covers Junctions 31 and 32. Note that  cover range may vary slightly if  flow direction changes.

Second, Rather than showing concentrations of all nodes and links, we plot the range zone and average value in Figs.~\ref{fig:concentration2} and~\ref{fig:concentration3}, and only results for specific junctions and pipes ($\mathrm{J11}$, $\mathrm{J21}$, $\mathrm{J31}$, $\mathrm{P11}$, $\mathrm{P31}$, and $\mathrm{P110}$) are plotted. Note that \textit{(i)} concentrations from some nodes or links (Pump 9, Junction 10, and Pipe 10 in Fig.~\ref{fig:setup}) cannot be adjusted by MPC controller are excluded when plotting range zone and average value, and \textit{(ii)} the average value does not reach the reference due to soft constraints instead of hard constraints.

\begin{table}[t]
	\caption{Number of optimization variables (percentage reduction) of MPC (worst  case scenario defined as booster stations installed at each node).}
	\centering
	\label{tab:scale}
	\renewcommand{\arraystretch}{1.2}
	\normalsize
	%	\resizebox{\linewidth}{!}{%
	\begin{tabular}{c|c|c|c|c|c}
		\hline
		\multirow{2}{*}{\textit{Networks}} & \multirow{2}{*}{\textit{\makecell{\# of com-\\ponents$^*$}} } & \multirow{2}{*}{\textit{LP~\eqref{equ:MPCLP}}} & \multicolumn{2}{c|}{\textit{\makecell{Unconstrained QP~\eqref{equ:minquad}}}} & \multirow{2}{*}{\textit{\makecell{QP~\eqref{equ:mpcscalable-2}}}} \\ \cline{4-5}
		&  &  & \textit{Solver~\eqref{equ:minquad} } & \textit{Analytical~\eqref{equ:controllaw}} &  \\ \hline
		\textit{\textit{\makecell{Three-node\\network}}} & \makecell{\{1,1,1,\\1,1,0\}}  & 32,100 & 900 (97\%) & 0 (100\%)   & 900 (97\%) \\ \hline
		\textit{Net1} & \makecell{\{9,1,1,\\12,1,0\}}   & 366,900 &  3,300 (99\%)  &  0 (100\%) & 3,300 (99\%)  \\ \hline
		\textit{Net3} & \makecell{\{92,2,3,\\117,2,0\}}  & 3,568,800 & 29,100 (99\%)  & 0 (100\%)   & 29,100 (99\%)  \\ \hline \hline
		\multicolumn{6}{l}{\footnotesize{
				\makecell{$^*$Number of each component in WDN: \{$n_\mathrm{J}$, $n_\mathrm{R}$, $n_\mathrm{TK}$, $n_\mathrm{P}$, $n_\mathrm{M}$, $n_\mathrm{V}$\}.} }}
	\end{tabular}%
	%	}
\end{table}

Third, the oscillations either from control actions in Fig.~\ref{fig:ControlActionU_net1} or from the control effects in Figs.~\ref{fig:concentration2} and~\ref{fig:concentration3} are seen  because \textit{(i)} Net1 network is looped and complex, and there exists a delay time for a junction to reach the desired value. For example,  it takes Junction 13 (we do not show the result for this junction due to space limitation) about $178$ minutes to receive the chlorine injected from the booster station at Junction 11, while it only takes 1 minute for Junction 11. This is also the reason why Junction 13 has a delay, and Junction 11 has an overshoot, %Besides, this phenomenon is also common in other nodes;   
\textit{(ii)}  the flow rates changes frequently, %see Fig~\ref{fig:net1setup_b}, 
and before the concentrations are adjusted to the desired value, another state change of flow rates has arrived which leads to tiny oscillations.

Fourth,  the concentration is  near $1.0$ mg/L during some time intervals for Pipe 110 (see Fig.~\ref{fig:concentration3_b}). For Pipe 110 during [612,716] minutes in Fig.~\ref{fig:concentration3_b}, the water is flowing out from Tank 2 (see $\mathrm{P110}$ in Fig.~\ref{fig:ControlActionU_net1}), and the concentration of Pipe 110 is decided by its upstream node concentration in Tank 2 ($1.0$ mg/L). Similar situation happens 10 times when it is in tank draining period in Fig.~\ref{fig:concentration3_b} and this is due to \textit{(i)} the flow rates are not optimization variables, and \textit{(ii)} we cannot control the flow direction after assuming the hydraulic model is solved first.

Next, we discuss the scalability of our proposed method. As we mentioned in Sections~\ref{sec:LP-WQC} and~\ref{sec:scalable-WQC}, $\mathrm{WQC-MPC}$~\eqref{equ:MPCLP}  has huge dimension and computational burden,  while the~\eqref{equ:minquad}, or~\eqref{equ:mpcscalable-2} is scalable, and it can be seen in Tab.~\ref{tab:scale}. We discuss scalability for three networks. Here we consider the worst case scenario, that is, booster stations are assumed to be installed at all nodes thereby resulting in more optimization variables. We also consider $N_p = 5$ minutes when $\Delta t = 1$ second, or $N_p = 300$ time-steps, $s_{L}$ = 100 segments. According to~\eqref{equ:numberofvar} and~\eqref{equ:numberofvar2}, we present the number of variables and its corresponding percentage reduction in Tab.~\ref{tab:scale}. The minimum percentage reduction is $97\%$, and if we use~\eqref{equ:controllaw}, there is no optimization variable at all, which proves our method is scalable. {The result of scalability test is shown in next section.}

\subsection{Net3 network}\label{sec:net3} 
Net3 network~\cite{rossman2000epanet} shown in Fig.~\ref{fig:setup}c, composed of 92 junctions, 2 reservoirs, 3 tanks, 117 pipes, and 2 pumps, is used for scalability test, and the uncertainty is not considered in this test case. Similarly, before discussing the final control effects, the effectiveness of the proposed LDE model is verified and shown in Fig.~\ref{fig:Error_EPANET_LDE}. The maximum error is 7.4\%, and the relative error is in [0.5\%, 3\%] when the L-W scheme is stable. 

There are two chlorine sources (0.5 mg/L) at the LAKE and RIVER  (marked as a green star in Fig.~\ref{fig:setup}c). The actuator placement problem is out the scope of this paper, the number of mass booster stations and the corresponding locations are decided subjectively so far as stated. Hence, we simply put three mass boosters at Junctions 217, 237, and 247 to improve the chlorine concentrations of junctions and pipes in a certain area (the zoomed-in area in Fig.~\ref{fig:setup}c) to the reference value $\m y^\mathrm{ref} = 0.6$ mg/L. 

\begin{figure}[t]
	\centering
	\includegraphics[width=0.56\linewidth]{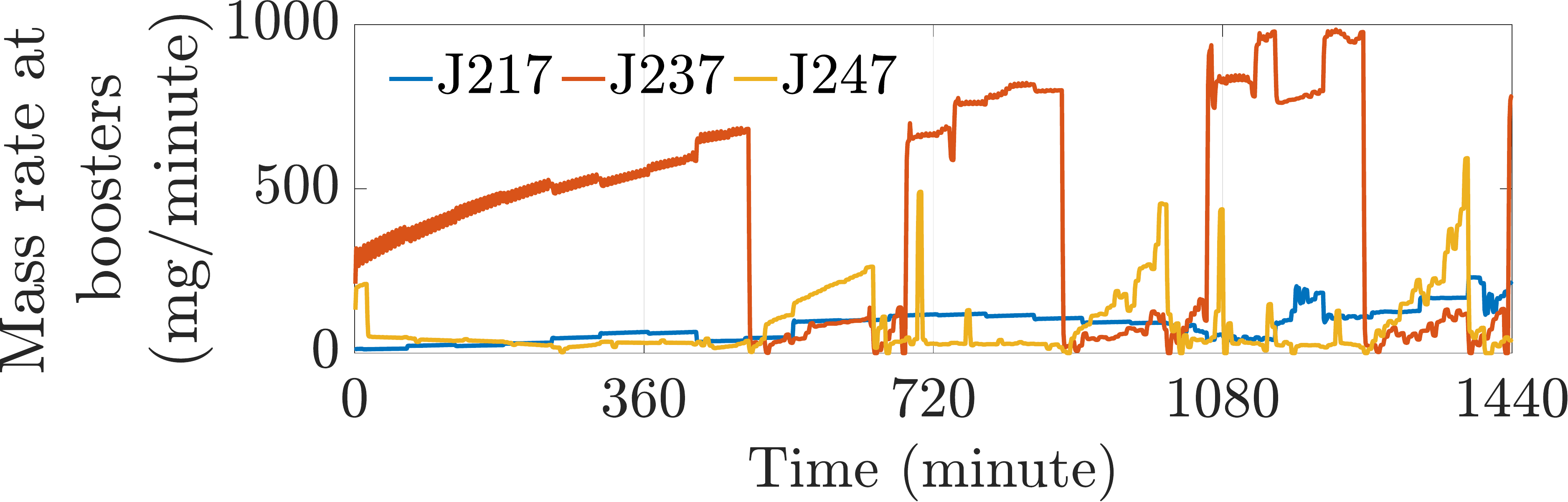}
	\caption{{Control action $\m u$ at three mass booster locations in Net3 (1 day).}}
	\label{fig:ControlActionU_net3}
\end{figure}

\begin{figure}[t]
	\centering
	\subfloat[\label{fig:net3_a}]{\includegraphics[keepaspectratio=true,scale=0.22]{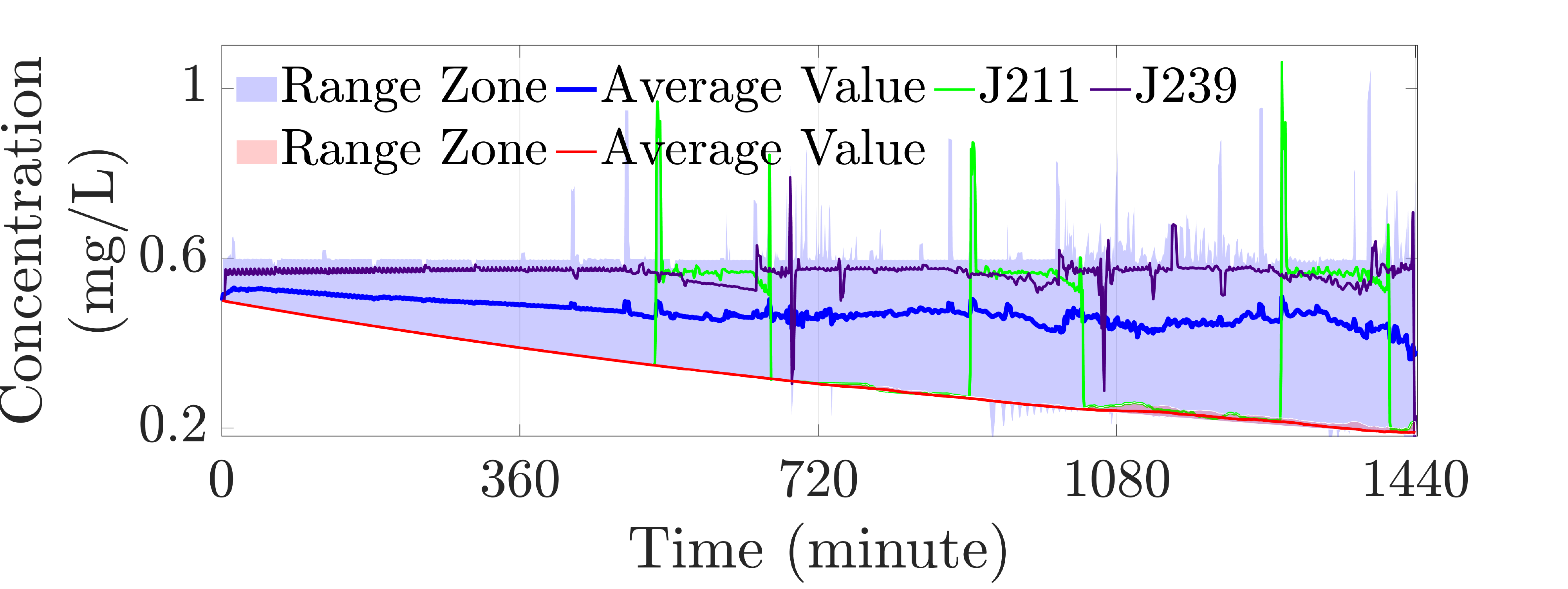}}{} \vspace{-0.1cm} 
	\subfloat[\label{fig:net3_b}]{\includegraphics[keepaspectratio=true,scale=0.22]{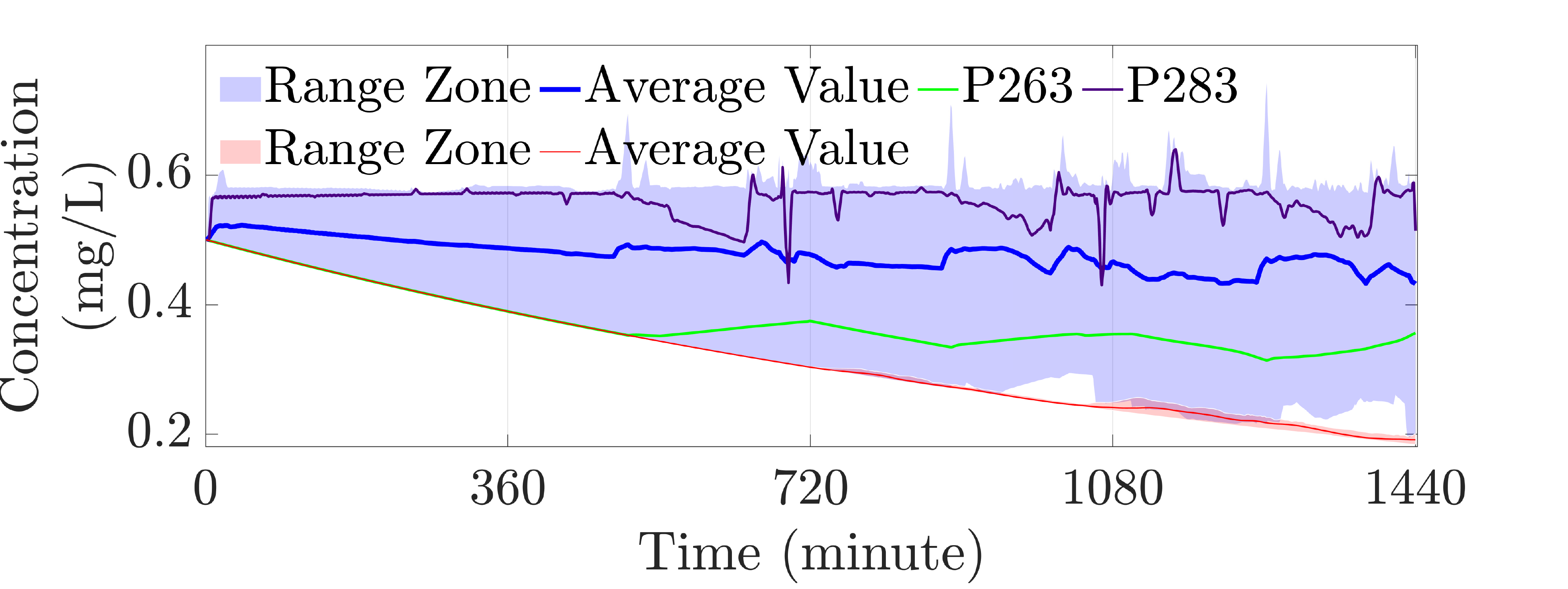}}{} \hspace{-0.2cm}\vspace{-0.2cm}
	\caption{{Chlorine concentrations in junctions (a) and pipes (b) in controlled area of Net3 (range zone and average value marked as red/blue are the results before/after applying MPC controller)}}
	\label{fig:net3}
\end{figure}

The time-steps for this simulation are $T_d = 24$ hours ($1440$ minutes), $T_n = T_h = 1$ hour, and $T_p = 5$ minutes.  The price $\lambda = 0.001 \$/\mathrm{mg}$, matrices  $\m Q$ and $\m R$ in MPC~\eqref{equ:controllaw} are $3\m I$ and $5\m I$ ($\m I $ is an identity matrix with corresponding dimensions). The MPC controller \eqref{equ:controllaw} is used to decide the optimal dosage of chlorine injected by booster stations (see control actions in Fig.~\ref{fig:ControlActionU_net3}) to maintain the proper concentration $\m y^{\mathrm{ref}}$ (see control effects in Fig.~\ref{fig:net3}) in the zoomed-in area. The main control scheme, results, and analysis are similar to previous test cases and not repeated. We simply present the difference in control effects with/without applying MPC controller next.

In Fig.~\ref{fig:net3_a}, the average concentration of all junctions in the zoomed-in area after applying MPC controller is maintained at about 0.52 mg/L (see blue line), while the average concentration of all junctions without applying MPC controller simply decays from 0.5 mg/L to 0.2 mg/L (see red line). This verifies the effectiveness of the designed MPC controller. Specifically, the concentrations at Junctions 211 and 239 are shown as examples, and we observe that they are maintained at $\m y^\mathrm{ref} = 0.6$ mg/L most time.  

In Fig.~\ref{fig:net3_b}, the average concentration of all pipes in the zoomed-in area with/without applying the MPC controller is shown as the blue/red line.  Comparing these two lines, we can observe the effectiveness of the MPC controller. In particular,  the concentration in Pipe 283 is almost maintained at the desired reference value 0.6 mg/L while the concentration in  Pipe 263 is only 0.4 mg/L which is far from the reference. This is due to Pipe 283, located between Junctions 237 and 247 (both are installed with mass booster stations), can be taken care of by both boosters while Pipe 263, located at the downstream of Junction 237 and far from it, can only be impacted by only one booster.

With the results from Net3, the scalability of our proposed method is verified.

\normalcolor

\subsection{Comparisons with rule-based control}

In this section, we perform thorough a case study to showcase the performance of our presented controller in comparison with traditional WQC through EPANET's built-in Rule-based Control (RBC). The simulation is performed for the three-node network in Fig.~\ref{fig:setup}a. 

In particular, RBC  can decide the chlorine dose with regarding to different deviation levels  thorough predefined rules. In our case, the deviation level is defined as $$\mathrm{deviation} = \frac{1}{2}\left( \sum_{s=1}^{100}( c_{23}^\mathrm{P} (s,t)-  y^\mathrm{ref}))/100 + ( c_{2}^\mathrm{J} (t)-  y^\mathrm{ref})\right)$$ which is always a number in $[- y^\mathrm{ref}, 0]$. Each rule is an IF-ELSE statement, that is, ``if $ a \leq \mathrm{deviation}  < b$, then $u=c$ mg/minute". This implies if the deviation is in $[a,b]$, then the booster should inject $c$ mg of chlorine at that time-step.

\begin{figure}[t]
	\centering
	\subfloat[\label{fig:Rule_a}]{\includegraphics[keepaspectratio=true,scale=0.22]{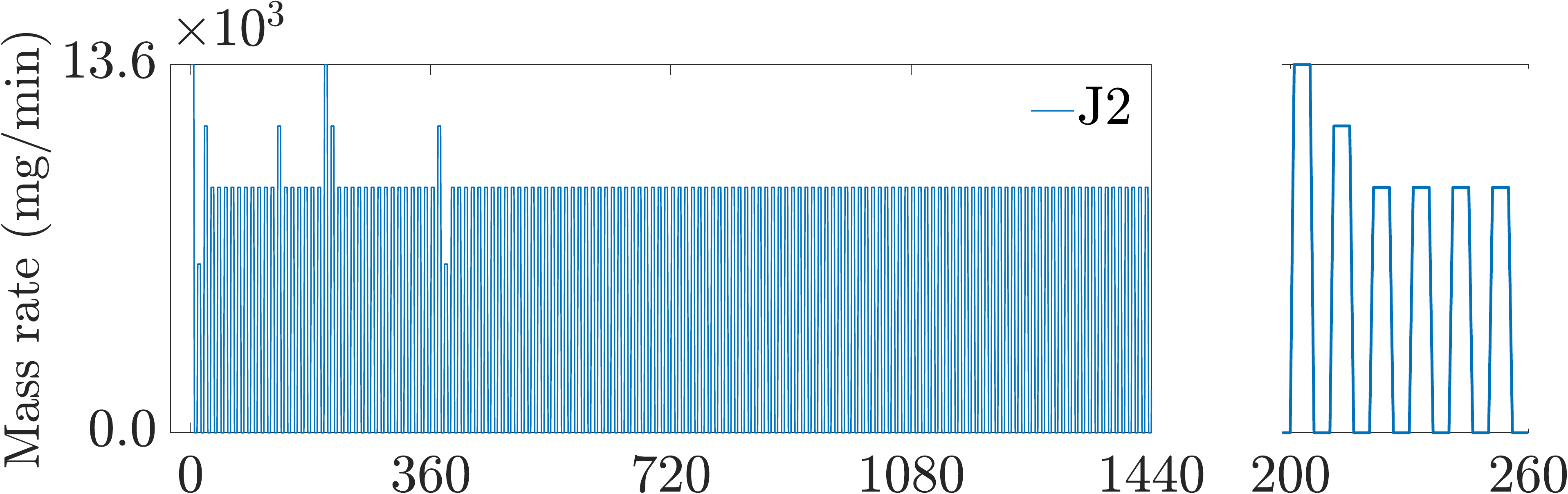}}{}\vspace{-0.1cm} \hspace{0.5cm}
	\subfloat[\label{fig:Rule_b}]{\includegraphics[keepaspectratio=true,scale=0.22]{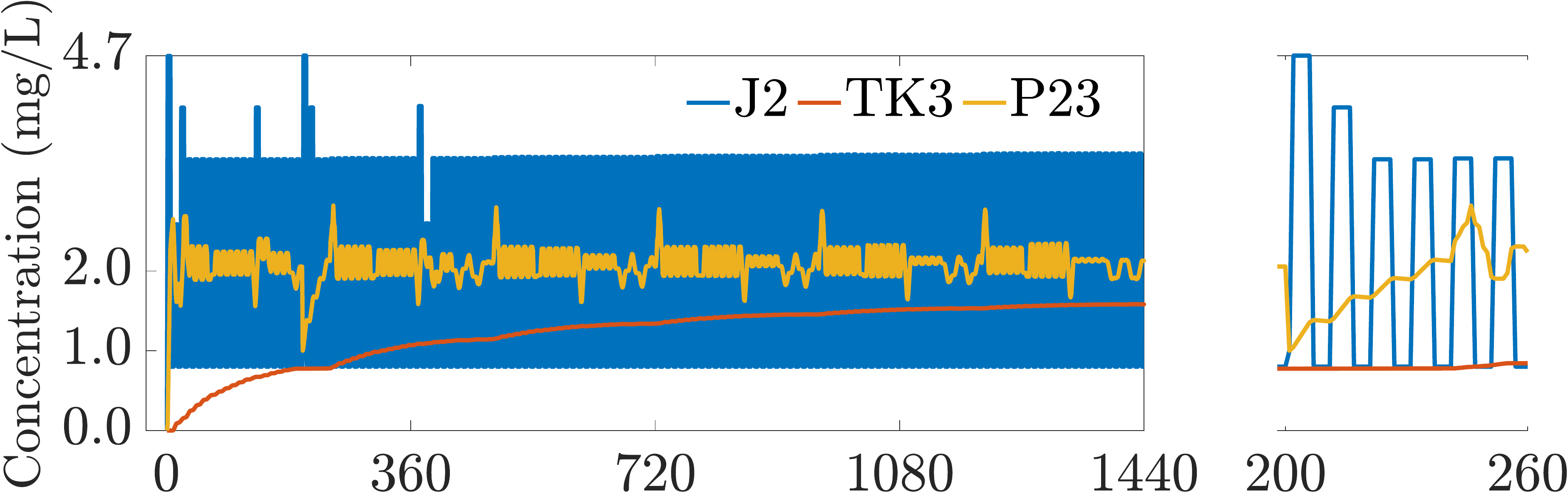}}{}\vspace{-0.2cm}
	\caption{(a) RBC control action $\m u$ during 1440 minutes for three-node network, (b) chlorine concentration at Junction 2, Tank 3, and Pipe 23 under three uncertainty sources, and unknown disturbance happens at 200-$\mathrm{th}$ minute.}
	\label{fig:Rule}
\end{figure}

To ensure fair comparison with the presented MPC approach, we apply all uncertainty sources to RBC as we did with MPC in the previous section. The RBC law is shown in Fig.~\ref{fig:Rule}. We observe that the overall result is inferior to the MPC law shown in Fig.~\ref{fig:concentration}, since \textit{(i)} RBC is not network-driven and cannot optimize a cost function as MPC can do, and \textit{(ii)} RBC exhibits larger oscillations due to its bang-bang nature. Besides that, conflicts exist among rules if we have multi-boosters, and we have to define priorities for rules. The objective function values for three-node network are compared and listed in Tab.~\ref{tab:compareresult}, and we observe that MPC outperforms RBC in all three metrics.

\begin{table}[t]
	\normalsize
	\setlength\tabcolsep{1pt}
	\centering
	\caption{{Comparison of performances for the three-node network. The three optimization objectives correspond with  $\frac{1}{2}(\m y^{\mathrm{ref}}-\m y_p)^{\top}\m Q(\m y^{\mathrm{ref}}-\m y_p)$ (for reference deviation), $\frac{1}{2} \Delta \m u_p^\top \m R \Delta \m u_p$ (for smoothness of the control action), and  $\m b^{\top} \Delta \m u_p$ (for the chlorine cost).  }}
	\setlength{\tabcolsep}{0.4em} % for the horizontal padding
	\renewcommand{\arraystretch}{1.5}
%	\resizebox{\linewidth}{!}{%
	\begin{tabular}{c|c|c|c|c}
		\hline
		\textit{Objectives}     &    \makecell{\textit{Reference Deviation} } & \makecell{\textit{Smoothness} }  & \makecell{\textit{Chlorine cost (\$)}} & \textit{Total $J(\m u)$} \\ \hline
		%\multirow{2}{*}{\textit{\makecell{Three-node\\network}}} 
		\textit{MPC} & $ 1.22\times  10^3 $ &$ 1.73\times  10^7 $&$ 5.99\times 10^3$ & $1.74\times 10^7 $  \\  \cline{1-5} 
		\textit{RBC} &  $ 3.73\times  10^3 $ & $ 2.42\times  10^{10} $& $6.64\times  10^3 $ &$ 2.42\times  10^{10} $ \\ 
		\hline \hline
	\end{tabular}%
%	}
	\label{tab:compareresult}
\end{table}

\section{\large Conclusions, Broader Impacts, and Paper Limitations}~\label{sec:limitations}
This paper presents the first high-fidelity and thorough state-space, control-oriented model of water quality dynamics in drinking water networks. A plug-and-play optimal model predictive control algorithm exploits this model to perform real-time water quality regulation. The case studies demonstrate that such control algorithm is fast in responding to network disturbances. 

The major broader impact of this work is two-fold. First, this study provides the control engineering research community a tool to compute the time-varying state-space matrices for water quality dynamics. This can be used to test many other more sophisticated control-theoretic techniques. Second, the proposed control algorithm can be used in future control centers when water quality measurements are performed online. This allows water utilities to remotely control chlorine dosages given real-time chlorine measurements in select network locations. 

This paper is not devoid of limitations. First, the paper only focuses on single-species reaction models. Multiple-species reaction models will then include nonlinear difference equations. Hence, an MPC for such a model will then involve a nonlinear, nonconvex constrained optimization problem that is hard to optimally solve. Second, although the linear MPC  shows promise in dealing with uncertainty, this controller cannot, in theory, guarantee a desired water quality performance although case studies show adequate performance. Third, the constrained state and input problem formulations is critical in practice, and an MPC formulation and implementation to guarantee the satisfaction of state and input bounds is required. Furthermore, the locations of booster stations and sensors and the parameters such as the reference value have a significant impact on the final control effects.  This study does not provide analysis on these aforementioned issues.  To that end, the authors' future work will focus on these important extensions: \textit{(i)} designing scalable computational methods to deal with a nonlinear MPC through convex relaxations; \textit{(ii)} \textit{robustifying} the proposed MPC to theoretically account for worst-case uncertainty; \textit{(iii)} exploring an MPC formulation guaranteeing the state and input bounds, the feasibility of the MPC problem with complex constraints, and the impacts of the locations of booster stations, sensors, and the MPC parameters such as reference value.

\section*{\large Acknowledgments}

This material is based upon work supported by the National Science Foundation under Grants CMMI-DCSD 1728629, 2015671, and 2015603. All the codes and tested networks are freely available in open data repository~\cite{wang_2020} for research reproducibility. 

\bibliographystyle{IEEEtran}
\bibliography{IEEEabrv,bibfile}

\newpage
\appendices
\section{\large An illustrative example of water quality model derivation of typical components}~\label{app:example}
A five-node illustrative example shown in Fig.~\ref{fig:example3} is used to help the reader understand the theoretical derivation of proposed water quality model in Section~\ref{sec:Ctrl-WQM}. Each node is installed with a booster station and each pipe will be divided into three segments for simplicity according to the Lax-Wendroff (L-W) scheme introduced next.  %When  of a component is hard to follow, we show a detailed example of that component from Fig.~\ref{fig:example3} to make our derivation further clear.

\begin{exmpl}~\label{exmp:pipe} According to L-W scheme, the concentration of each segment of $\mathrm{P}23$ in Fig.~\ref{fig:example3} is
	%\begin{linenomath*}
		\begin{subequations}~\label{equ:Pipe23}
			\begin{align}
			\hspace{-1.em}	\begin{split} \label{equ:Pipe23A}
			c^{\mathrm{P}}_{23}(1,t+ \Delta t)  &= \underline{\alpha}_{23}(t) c^{\mathrm{J}}_{2}(t)  + {\alpha}_{23}(t) 	c^{\mathrm{P}}_{23}(1,t)   + \overline{\alpha}_{23}(t) 	c^{\mathrm{P}}_{23}(2,t) + r_{23}(c^{\mathrm{P}}_{23}(1,t)),
			\end{split}\\
			\hspace{-1.em}	\begin{split} \label{equ:Pipe23B}
			c^{\mathrm{P}}_{23}(2,t+ \Delta t)  &= \underline{\alpha}_{23}(t) c^{\mathrm{P}}_{23}(1,t) +  {\alpha}_{23}(t) 	c^{\mathrm{P}}_{23}(2,t)  + \overline{\alpha}_{23}(t) 	c^{\mathrm{P}}_{23}(3,t)  +  r_{23}(c^{\mathrm{P}}_{23}(2,t)),
			\end{split}\\
			\hspace{-1.em}	\begin{split} \label{equ:Pipe23C}
			c^{\mathrm{P}}_{23}(3,t+ \Delta t)  &= \underline{\alpha}_{23}(t) c^{\mathrm{P}}_{23}(2,t)    + {\alpha}_{23}(t) 	c^{\mathrm{P}}_{23}(3,t)     +  \overline{\alpha}_{23}(t) 	c^{\mathrm{J}}_{3}(t) + r_{23}(c^{\mathrm{P}}_{23}(3,t)).
			\end{split}
			\end{align}
		\end{subequations}
	%\end{linenomath*}
	Similarly, we can obtain equations of each segment in $\m c^{\mathrm{P}}_{24}(t+\Delta t)$ and $\m c^{\mathrm{P}}_{52}(t+\Delta t)$ in Fig.~\ref{fig:example3}.
\end{exmpl}

\begin{exmpl} ~\label{exmp:junction}
	According to mass balance at $\mathrm{J3}$ and Assumption~\ref{asmp:mixing},  we have
	%\begin{linenomath*}
		\begin{subequations}~\label{equ:J3}
			\begin{align}
					& c^\mathrm{J}_{3}(t + \Delta t) = c^\mathrm{V}_{34}(t + \Delta t), \label{equ:JunctionBalanceExamplea}\\
				\begin{split} \label{equ:JunctionBalanceExampleb}
			& q^\mathrm{B}_3(t + \Delta t) c^\mathrm{B}_3(t + \Delta t)  + q_{23}(t + \Delta t)c^\mathrm{P}_{23}(3,t + \Delta t)   =  q^\mathrm{D}_3(t + \Delta t) c^\mathrm{J}_{3}(t + \Delta t) + q_{34}(t + \Delta t)  c^\mathrm{V}_{34}(t + \Delta t),
			\end{split} 
			\end{align}
		\end{subequations}
	%\end{linenomath*}
	where $c^\mathrm{P}_{23}(3,t+ \Delta t)$ is the third segment of $\mathrm{P}23$ in Fig.~\ref{fig:example3}, and is already derived as~\eqref{equ:Pipe23C}. After substituting~\eqref{equ:Pipe23C} and ~\eqref{equ:JunctionBalanceExamplea} into~\eqref{equ:JunctionBalanceExampleb}, we obtain
	%\begin{linenomath*}
		\begin{align*}
		c_{3}^\mathrm{J}(t + \Delta t) & =   a^\mathrm{J}_\mathrm{J}(t) c_3^\mathrm{J}(t)  + a^\mathrm{P}_\mathrm{J}(2,t) c_{23}^\mathrm{P}(2,t) + a^\mathrm{P}_\mathrm{J}(3,t) c_{23}^\mathrm{P}(3,t)+ b^\mathrm{B}_\mathrm{J}(t+\Delta t) c_3^\mathrm{B}(t + \Delta t) + R_3^\mathrm{J}(t + \Delta t) c^{\mathrm{P}}_{23}(3,t),
		\end{align*}
	%\end{linenomath*}
	where  the coefficients are
	%\begin{linenomath*}
		\begin{align*}
		\hspace{-0.5em} a^\mathrm{J}_\mathrm{J}(t) &= \overline{\alpha}_{23}(t) \beta_1(t + \Delta t),\  a^\mathrm{P}_\mathrm{J}(2,t) = \underline{\alpha}_{23}(t) \beta_1(t+\Delta t), \  a^\mathrm{P}_\mathrm{J}(3,t) = {\alpha}_{23}(t) \beta_1(t+\Delta t), \\ 
		 b^\mathrm{B}_\mathrm{J}(t) &= \beta_2(t + \Delta t),  \  R_3^\mathrm{J}(t + \Delta t) = \beta_1(t + \Delta t) k_{23}^\mathrm{P}.
		\end{align*}
	%\end{linenomath*} 
	Note that $\beta_1(t+\Delta t) = \frac{q_{23}(t+\Delta t)}{q_{34}(t+\Delta t) + q_3^\mathrm{D}(t+\Delta t)}$ and  $\beta_2(t+\Delta t) = \frac{q_3^{\mathrm{B}}(t+\Delta t)}{q_{34}(t+\Delta t) + q_3^\mathrm{D}(t+\Delta t)}$. Similarly, after listing mass balance equations at $\mathrm{J3}$ and $\mathrm{J4}$, $c_{2}^\mathrm{J}(t + \Delta t)$ and $c_{4}^\mathrm{J}(t + \Delta t)$ can be derived, but omitted here.
	
	%where $a^\mathrm{J}_\mathrm{J}(t) = \alpha_{ki}(+1,t) \beta_1(t+\Delta t)$ is an element of $A^\mathrm{J}_\mathrm{J}(t)$ in~\eqref{equ:abs-nodes-mass}, $a^\mathrm{P}_\mathrm{J}(2,t) = \alpha_{ki}(-1,t) \beta_1(t+\Delta t)$ and $a^\mathrm{P}_\mathrm{J}(3,t) = \alpha_{ki}(0,t) \beta_1(t+\Delta t)$ are elements of $A^\mathrm{L}_\mathrm{J}(t)$, and $b^\mathrm{B}_\mathrm{J}(t) = \beta_2(t)$ is an element of $B^\mathrm{J}(t)$. 
\end{exmpl}

\begin{exmpl} ~\label{exmp:selectmatrix4J}
	Consider that four links ($\mathrm{M}12$, $\mathrm{P}52$, $\mathrm{P}24$, $\mathrm{P}23$) are connected by a junction ($\mathrm{J}2$) in Fig.~\ref{fig:example3} and selection matrix $\m S^\mathrm{in}_\mathrm{J}$ selecting inflows can be obtained directly by changing $-1$ to $0$ in $-\m E^\mathrm{L}_\mathrm{J2}$ which is defined as $\begin{bmatrix}
	\m E_\mathrm{J2}^\mathrm{P} & \m E_\mathrm{J2}^\mathrm{M} & \m E_\mathrm{J2}^\mathrm{V} 
	\end{bmatrix}$;  $\m S^\mathrm{out}_\mathrm{J}$ selecting outflows can be obtained by changing $-1$ to $0$ in $\m E^\mathrm{L}_\mathrm{J2}$. That is, 
	%\begin{linenomath*}
		\begin{equation*} 
		\includegraphics[width=0.5\linewidth,valign = c]{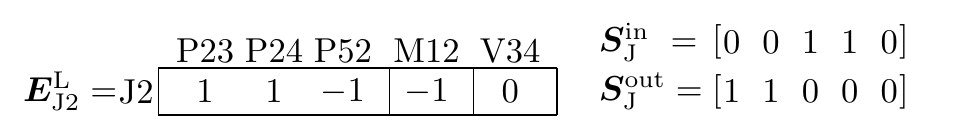}
		\end{equation*}
	%\end{linenomath*}
\end{exmpl}

\begin{exmpl} ~\label{exmp:tank}
	According to mass balance at $\mathrm{TK5}$ shown in Fig.~\ref{fig:example3} and Assumption~\ref{asmp:mixing}, we have
	%\begin{linenomath*}
		\begin{align*}
		c_{5}^{\mathrm{TK}}(t +  \Delta t) 
		&= a^\mathrm{TK}_\mathrm{TK}(t +  \Delta t) c_{5}^{\mathrm{TK}}(t) + b^\mathrm{B}_\mathrm{TK}(t +  \Delta t) c_{5}^{\mathrm{B}}(t +  \Delta t)  + R^{\mathrm{TK}}r^{\mathrm{TK}}(t +  \Delta t)(c_{5}^{\mathrm{TK}}(t)),\notag 
		\end{align*}
	%\end{linenomath*}
	where coefficients $a^\mathrm{TK}_\mathrm{TK}(t +  \Delta t) =   \frac{V_{5}^{\mathrm{TK}}(t) -q_{52}(t)\Delta t}{V_{5}^{\mathrm{TK}}(t  +  \Delta t)} $,  
	$b^\mathrm{B}_\mathrm{TK}(t +  \Delta t) =\frac{V_{5}^{\mathrm{B}}(t +  \Delta t)}{V_{5}^{\mathrm{TK}}(t +  \Delta t)}$, and $R^{\mathrm{TK}}(t +  \Delta t) = \frac{  \Delta t}{V_{5}^{\mathrm{TK}}(t +  \Delta t)}$.
\end{exmpl}

\begin{exmpl} \label{exmp:reservoir}
	The concentration difference equation of $\mathrm{R}1$ in Fig.~\ref{fig:example3}  is
	$c_{1}^{\mathrm{R}}(t +  \Delta t) = 	c_{1}^{\mathrm{R}}(t)$.
\end{exmpl}

\begin{exmpl} ~\label{exmp:selectmatrix}
	Consider that the pump $\mathrm{M12}$ and the valve $\mathrm{V34}$  are connected as shown in Fig.~\ref{fig:example3}, and Matrix  $\m S^\mathrm{N}_\mathrm{M}$ ($\m S^\mathrm{N}_\mathrm{V}$) is obtained simply by changing $-1$ to $0$ in $\m E^\mathrm{N}_\mathrm{M}$ ($\m E^\mathrm{N}_\mathrm{V}$) which is defined as $\begin{bmatrix}
	\m E^\mathrm{J}_\mathrm{M} & \m E^\mathrm{R}_\mathrm{M} & \m E^\mathrm{TK}_\mathrm{M} 
	\end{bmatrix}$ $\left(\begin{bmatrix}
	\m E^\mathrm{J}_\mathrm{V} & \m E^\mathrm{R}_\mathrm{V} & \m E^\mathrm{TK}_\mathrm{V} 
	\end{bmatrix}\right)$.  That is
	%\begin{linenomath*}
		\begin{equation*} 
		\includegraphics[width=0.5\linewidth,valign = c]{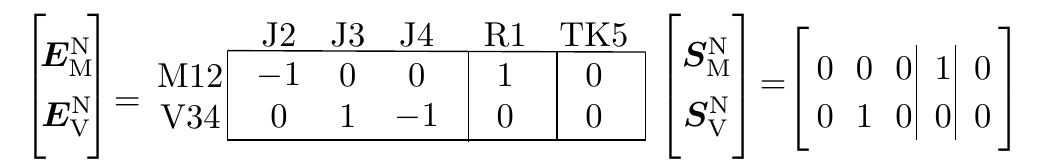}
		\end{equation*}
	%\end{linenomath*}

	We can see $\m S^\mathrm{N}_\mathrm{M}$ ($\m S^\mathrm{N}_\mathrm{V}$) select $\mathrm{R}1$ ($\mathrm{J3}$) that are the upstream node of $\mathrm{M12}$ ($\mathrm{V34}$).
\end{exmpl}

\begin{exmpl} ~\label{exmp:pump} 
	Since the upstream node of $\mathrm{M}12$ is $\mathrm{R1}$, and according to~\eqref{equ:pumpi},  $   c_{12}^\mathrm{M}(t+\Delta t) = c_{12}^\mathrm{M}(t).$
	
	Since the upstream node of $\mathrm{V}34$ is $\mathrm{J3}$, we have $c_{34}^\mathrm{V}(t+\Delta t) =  c_3^\mathrm{J}(t + \Delta t)$, and $c_3^\mathrm{J}(t + \Delta t)$ is available in Example~\ref{exmp:junction}. Hence,
	%\begin{linenomath*}
		\begin{align*}
		c_{34}^\mathrm{V}(t + \Delta t) & =   a^\mathrm{J}_\mathrm{J}(t) c_3^\mathrm{J}(t)  + a^\mathrm{P}_\mathrm{J}(2,t) c_{23}^\mathrm{P}(2,t) + a^\mathrm{P}_\mathrm{J}(3,t) c_{23}^\mathrm{P}(3,t)+ b^\mathrm{B}_\mathrm{J}(t+\Delta t) c_3^\mathrm{B}(t + \Delta t) + R^\mathrm{J} r_{23}(c^{\mathrm{P}}_{23}(3,t)),
		\end{align*}
	%\end{linenomath*}
	where all parameters are the same as the ones in Example~\ref{exmp:junction}.
\end{exmpl}
\begin{figure}[t]
%	\vspace{-1em}
	\centering
	\includegraphics[width=0.65\linewidth]{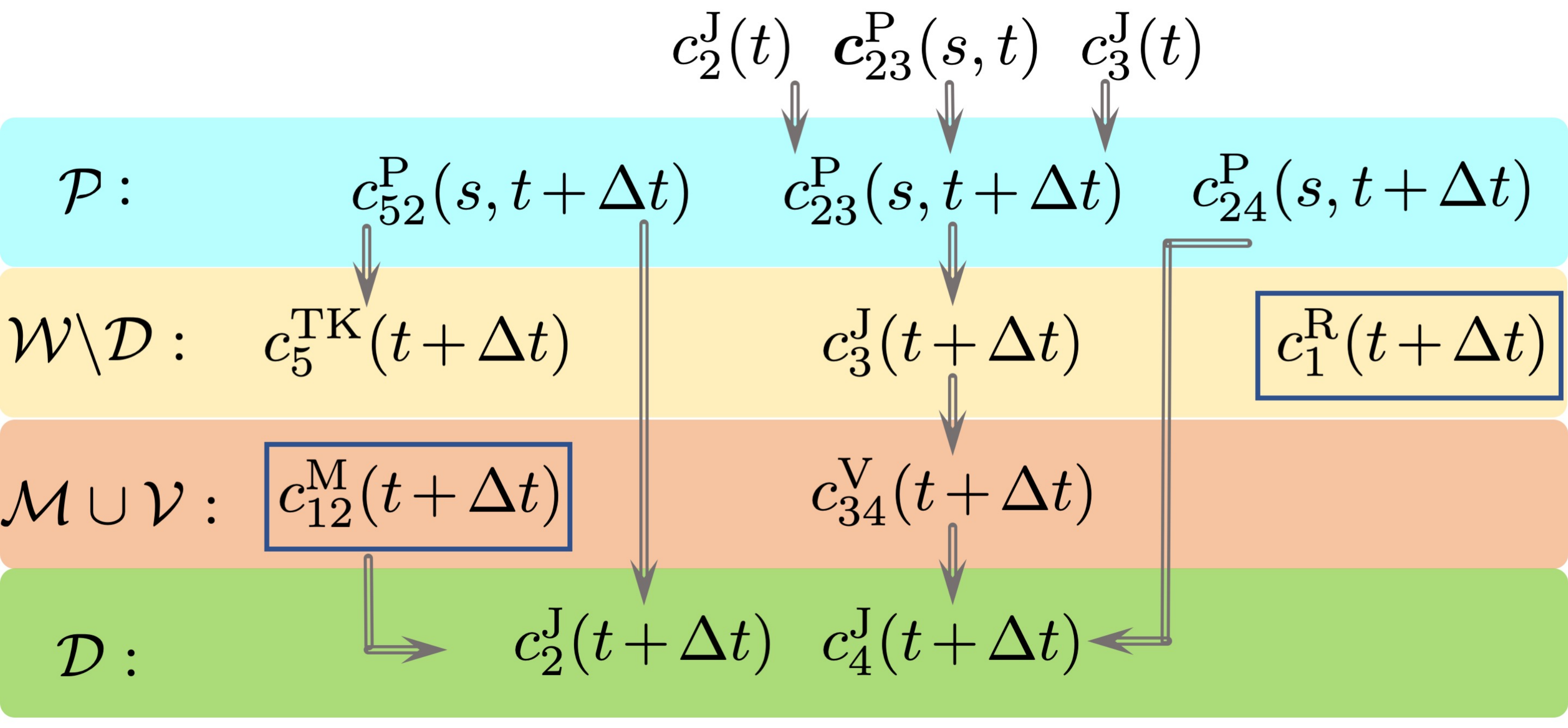}
	\caption{{Concentration dependence forest of the illustrative example in Fig.~\ref{fig:example3} (the segment in a pipe is denoted as $s$).}}
	\label{fig:Example2findMatrix}
	\vspace{-1em}
\end{figure}

\begin{exmpl} ~\label{exmp:dependence}
	According to the topology in Fig.~\ref{fig:example3},  we have 
	%\begin{linenomath*}
		\begin{align*}
		\mathcal{P} = \{\mathrm{P23}, \mathrm{P24}, \mathrm{P52} \}, \mathcal{M} = \{\mathrm{M12}\}, \mathcal{V} = \{\mathrm{V34}\}, 
		\mathcal{D} = \{\mathrm{J2}, \mathrm{J4}\}, \mathcal{W\char`\\ D} = \{\mathrm{R1}, \mathrm{J3},\mathrm{TK5} \}. 
		\end{align*}
	%\end{linenomath*}
	%We denote the first segments in pipes as $\mathrm{Pij}(1)$,  the set $\mathcal{F} = \{\mathrm{P23(1)}, \mathrm{P24(1)}, \mathrm{P52(1)} \}$, and the non-first segments of pipes are collected in $ \mathcal{P} \char`\\ \mathcal{F}$. Hence, the corresponding concentration  in $ \mathcal{P} \char`\\ \mathcal{F}$ is $c_{ij}^\mathrm{P}(s,t+\Delta t)$ for $s \neq 1$, $ij \in \mathcal{P}$, and the concentration  in $\mathcal{F}$ is $c_{ij}^\mathrm{P}(1,t+\Delta t)$ for $ij \in \mathcal{P}$.

	The concentration dependence forest for this case shown in Fig.~\ref{fig:Example2findMatrix} includes three independent trees.  This forest shows the details of how state variable of each component from $t$ to $t+\Delta t$ are connected.  Specifically, $\mathrm{Reservoir}$ $1$ ($c^\mathrm{R}_1$) is independent of the rest, and forms an independent tree.%, so does $\mathrm{Pump}$ $12$ ($c^\mathrm{M}_{12}$). 
	
	Let us take the tree in which the root is $\mathrm{Pipe}$ $23$ ($\mathrm{P}23$) as an example. According to~\eqref{equ:Pipe23} in Example~\ref{exmp:pipe}, $c^\mathrm{P}_{23}(s,t+\Delta t)$ can be expressed by $c^\mathrm{J}_2(t)$, $c^\mathrm{P}_{23}(s,t)$, and $c^\mathrm{J}_3(t)$. The downstream node of $\mathrm{Pipe}$ $23$ is $\mathrm{J3}$, and $c^\mathrm{J}_3(t+\Delta t)$ depends on $c^\mathrm{P}_{23}(s,t+\Delta t)$, see $Example~\ref{exmp:junction}$. Moreover, the downstream link of $\mathrm{J3}$ is $\mathrm{V34}$. Hence, $c^\mathrm{V}_{34}(t+\Delta t)$ depends on $c^\mathrm{J}_{3}(s,t+\Delta t)$, see $Example~\ref{exmp:pump}$. The downstream node of  $\mathrm{V34}$ is $\mathrm{J4}$, which also is the downstream node of $\mathrm{Pipe}$ $24$. Therefore,  $c^\mathrm{J}_4(t+\Delta t)$ depends on $c^\mathrm{V}_{34}(t+\Delta t)$ and  $c^\mathrm{P}_{24}(s,t+\Delta t)$ simultaneously.
	
	After all concentration at $t+\Delta t$ are available, the corresponding matrix form for each component type can be obtained. 
\end{exmpl}
\normalcolor
\section{\large Matrix form derivation  for junctions}~\label{app:derivation_Junction}
%\section{\large Matrix form derivation for junctions}~\label{app:derivation}
The Hadamard product or division is commutative, associative and distributive over addition~\cite{horn2012matrix}. For example, if $\m x$, $\m y$, and $\m z$ are vectors of the same size, and all elements in $\m y$ are non-zeros, then we have  $\m x  \oslash \m y \circ \m z= \m x  \circ \m z  \oslash \m y$. We present another property that is frequently used  in the later derivation of water quality modeling as
\begin{mypro}\label{prp:hadama}
	If $\m x$ and $\m y$ are vectors of the same size, and $\m E$ is a matrix with proper dimension, then 
	%\begin{linenomath*}
		$$(\m E \m x) \circ (\m E \m y)  =  \m E (\m x \circ \m y) = \m E \diag(\m x)  \m y. $$
	%\end{linenomath*}
\end{mypro}

With the above properties,  the detail of matrix form derivation of water quality modeling for various components (junction, tanks) is presented next.

Links $ki$ and $ij$ can be pipes, valves, pumps, or their combinations, that is, $\m c^\mathrm{L}(\m s_{L},t + \Delta t) \triangleq \{\m c^\mathrm{P}(\m s_{L},t + \Delta t), \m c^\mathrm{M}(t + \Delta t), \m c^\mathrm{V}(t + \Delta t),\}$, where $\m c^\mathrm{P}(\m s_{L},t + \Delta t)$ can be expressed by $\m S^\mathrm{P} \m c^\mathrm{P}$ where $\m S^\mathrm{P} \in \mbb{R}^{(n_{\mathrm{P}} \cdot s_{L}) \times {(n_{\mathrm{P}} \cdot s_{L}) }}$ is a matrix selecting the last segment of Pipes. Note that $\m c^\mathrm{L}(\m s_{L},t + \Delta t) = \m S^\mathrm{L} \m c^\mathrm{L}(t + \Delta t)$.

According to Assumption~\ref{asmp:mixing} or~\eqref{equ:node2pipe}, we know that $c^\mathrm{J}_{i}(t + \Delta t)$ equals $c^\mathrm{P}_{ij}(1,t + \Delta t)$, $c^\mathrm{M}_{ij}(t + \Delta t)$, or $c^\mathrm{V}_{ij}(t + \Delta t)$. Then we consider the simplification of~\eqref{equ:tempJunc} by moving $\m q^{\mathrm{out}}(t + \Delta t) \triangleq  \m q^{\mathrm{out}}_\mathrm{J}(t + \Delta t) + \m q^\mathrm{D}(t + \Delta t)$ to the right-hand side, and we have
%\begin{linenomath*}
	\begin{align} 
	\hspace{-0.8em}\m c^\mathrm{J}(t + \Delta t) &=    \m q^\mathrm{in}_\mathrm{J}(t + \Delta t) \oslash \m q^{\mathrm{out}}(t + \Delta t) \circ ({\diag(\m S^\mathrm{in}_\mathrm{J})} \m c^\mathrm{L}(\m s_{L},t +  \Delta t)) +  \m q^\mathrm{B}_\mathrm{J}(t + \Delta t) \oslash \m q^{\mathrm{out}}(t + \Delta t)   \circ   \m c^\mathrm{B}_\mathrm{J}(t + \Delta t).  ~\label{equ:tempJunction} 
	\end{align}
%\end{linenomath*}

Note that  According to  the  commutative property and Property~\ref{prp:hadama},  we derive the first item of~\eqref{equ:tempJunction} as
%\begin{linenomath*}
	\begin{align*}
	& \m q^\mathrm{in}_\mathrm{J}(t + \Delta t) \oslash \m q^{\mathrm{out}}(t + \Delta t)   \circ ({\diag(\m S^\mathrm{in}_\mathrm{J})} \m c^\mathrm{L}(\m s_{L},t + \Delta t))\\
	%& = \m q^\mathrm{in}_\mathrm{J}(t + \Delta t) \circ (\m S^\mathrm{L} \m c^\mathrm{P}(t + \Delta t)) \oslash \m q^{\mathrm{out}}(t + \Delta t)    \\
	& =  {\left(\diag(\m S^\mathrm{in}_\mathrm{J}) \diag(\m q^\mathrm{L}(t + \Delta t)) \m S^\mathrm{L} \m c^\mathrm{L}(t + \Delta t) \right) \oslash \m q^{\mathrm{out}}(t + \Delta t)}   \\
	& = {\underbrace{(\diag(\m q^{\mathrm{out}}(t + \Delta t)  ))^{-1} \diag(\m S^\mathrm{in}_\mathrm{J}) \diag(\m q^\mathrm{L}(t + \Delta t))  \m S^\mathrm{L}}_{\textstyle \m A^\mathrm{L}(t+\Delta t)}} \m c^\mathrm{L}(t + \Delta t)
	\end{align*}
%\end{linenomath*}
Similarly, for the second item of~\eqref{equ:tempJunction}, we have
%\begin{linenomath*}
	\begin{align*}
	& \m q^\mathrm{B}_\mathrm{J}(t+\Delta t) \oslash (\m q^{\mathrm{out}}(t+\Delta t))  \circ   \m c^\mathrm{B}_\mathrm{J}(t+\Delta t) = \underbrace{(\diag(\m q^{\mathrm{out}}(t+\Delta t))^{-1} \m E^\mathrm{B}_\mathrm{J} \diag(\m q^\mathrm{B}(t+\Delta t))}_{\textstyle \m B^\mathrm{J}(t+\Delta t)} \m c^\mathrm{B}(t+\Delta t).
	\end{align*}
%\end{linenomath*}
Note that $\m q^\mathrm{B}_\mathrm{J}(t+\Delta t) = \m E^\mathrm{B}_\mathrm{J} \m q^\mathrm{B}(t+\Delta t)$ and $\m c^\mathrm{B}_\mathrm{J}(t+\Delta t) = \m E^\mathrm{B}_\mathrm{J} \m c^\mathrm{B}(t+\Delta t)$, and  Property~\ref{prp:hadama} is applied to the above derivation.  
Hence,~\eqref{equ:tempJunction} can be written as
%\begin{linenomath*}
	\begin{align*}
	\hspace{-0.8em}\m c^\mathrm{J}(t + \Delta t) &=  \m A^\mathrm{L}(t +  \Delta t)  \m c^\mathrm{L}(t +  \Delta t) +   \m B^\mathrm{J}(t+\Delta t)  \m c^\mathrm{B}_\mathrm{J}(t + \Delta t).
	\end{align*}
%\end{linenomath*}
Since $\m c^\mathrm{L}(t + \Delta t)$ lumping $\m c^\mathrm{P,M,V}(t + \Delta t)$ can be expressed by $\m c^\mathrm{L}(t)$, $\m c^\mathrm{J}(t)$, etc. After substituting, we can get the matrix form.  Next, we only give an example of pipes, that is, 
%\begin{linenomath*}
	\begin{align} ~\label{equ:tempJunction2}
	\hspace{-0.8em}\m c^\mathrm{J}(t + \Delta t) =  \m A^\mathrm{P}(t +  \Delta t)  \m c^\mathrm{P}(t +  \Delta t) + \m B^\mathrm{J}(t+\Delta t)  \m c^\mathrm{B}_\mathrm{J}(t + \Delta t). 
	\end{align}
%\end{linenomath*}

According to~\eqref{equ:abs-pipe-mass}, we have
%\begin{linenomath*}
	\begin{align}~\label{equ:pipe-mass}
	\hspace{-0.7em}\m c^\mathrm{P}(t+\Delta t) = \m A^\mathrm{J}_\mathrm{P}(t) \m c^\mathrm{J}(t) + \m A^\mathrm{P}_\mathrm{P}(t) \m c^\mathrm{P}(t) + \m r^\mathrm{P}(\m c^\mathrm{P}(t)).
	\end{align}
%\end{linenomath*}

After substituting~\eqref{equ:pipe-mass} to~\eqref{equ:tempJunction2}, we have matrix form of Junction $i$ at $t+\Delta t$ when links are all pipes as  
%\begin{linenomath*}
	\begin{align}
	\m c^\mathrm{J}(t+\Delta t)  &= \m A_\mathrm{J}^\mathrm{J}(t+\Delta t)   \m c^\mathrm{J}(t) + \m A_\mathrm{J}^\mathrm{P}(t+\Delta t ) \m c^\mathrm{P}(t )   +  \m B^\mathrm{J}(t + \Delta t)   \m c^\mathrm{B}(t + \Delta t) + \m R^\mathrm{J}(t + \Delta t) \m r(\m x(t)),  ~\label{equ:junction-pipe}
	\end{align}
%\end{linenomath*}
where $ \m A_\mathrm{J}^\mathrm{J}(t+\Delta t)  =  \m A^\mathrm{P}(t+\Delta t) \m A^\mathrm{J}_\mathrm{P}(t)$, $ \m A_\mathrm{J}^\mathrm{P}(t+\Delta t )  = \m A^\mathrm{P}(t+\Delta t) \m A^\mathrm{P}_\mathrm{P}(t)$, $\m R^\mathrm{J}(t + \Delta t)$ is $\begin{bmatrix}\m A^\mathrm{P}(t + \Delta t) & \m O\end{bmatrix}$ and $\m r(\m x(t))$ is defined as $\begin{bmatrix} \m r(\m c^\mathrm{P}(t)) \\ \m r(\m c^\mathrm{TK}(t)) \end{bmatrix}$ according to~\eqref{equ:reactionVector}.

The final result considering pumps and valves is similar to~\eqref{equ:junction-pipe}, and after replacing $\m A_\mathrm{J}^\mathrm{P}(t+\Delta t ) \m c^\mathrm{P}(t)$ with $\m A_\mathrm{J}^\mathrm{L}(t + \Delta t) \m c^\mathrm{L}(t)$ in~\eqref{equ:junction-pipe}, the final result for junction is~\eqref{equ:abs-nodes-mass} in Tab.~\ref{table:wqp}, where $\m A_\mathrm{J}^\mathrm{L}$ is \textit{contribution matrix} from links for the concentration of junctions.

Note that each element in $\m q^{\mathrm{out}}_\mathrm{J} + \m q^\mathrm{D}$ is non-zero. Otherwise, it means there is a junction consuming no water or transporting zero water to next junction, and it is not necessary to calculate the concentrations at junctions in such case.   Hence, $\diag(\m q^{\mathrm{out}}_\mathrm{J} + \m q^\mathrm{D}))^{-1}$ exists, and Equation~\eqref{equ:junction-pipe} holds true.

\section{\large Matrix form derivation  for tanks}~\label{app:derivation_Tank}
We consider the simplification of~\eqref{equ:abs-tank-temp} by moving $\m V^{\mathrm{TK}}(t+\Delta t)$ to the right-hand side, and we have
%\begin{linenomath*}
	\begin{align} 
 \m c^{\mathrm{TK}}(t+\Delta t) 	&  = (\m V^{\mathrm{TK}}(t) - \Delta t \,  \m q^{\mathrm{out}}_\mathrm{TK}(t) ) \oslash \m V^{\mathrm{TK}}(t+\Delta t) \circ  \m c^{\mathrm{TK}}(t)  + \Delta t\ \m q^\mathrm{in}_\mathrm{TK}(t) \oslash \m V^{\mathrm{TK}}(t+\Delta t) \circ \m c^{\mathrm{P}}(t) \notag\\
	&+  \m V^\mathrm{B}_\mathrm{TK}(t+\Delta t)   \oslash \m V^{\mathrm{TK}}(t+\Delta t)  \circ \m c^{\mathrm{B}}_\mathrm{TK}(t+\Delta t)   + \Delta t\ \m r^\mathrm{TK}(\m c^\mathrm{TK}), \notag
	\end{align}
%\end{linenomath*}
where the first item is
%\begin{linenomath*}
	\begin{align*}
& (\m V^{\mathrm{TK}}(t) - \Delta t \,  \m q^{\mathrm{out}}_\mathrm{TK}(t) ) \oslash \m V^{\mathrm{TK}}(t+\Delta t) \circ  \m c^{\mathrm{TK}}(t) = \underbrace{\left(\diag(\m V^{\mathrm{TK}}(t+\Delta t))\right)^{-1} \diag(\m V^{\mathrm{TK}}(t) - \Delta t \,  \m q^{\mathrm{out}}_\mathrm{TK}(t)) }_{\textstyle \m A_{\mathrm{TK}}^\mathrm{TK}(t)  }\m c^{\mathrm{TK}}(t),
	\end{align*}
%\end{linenomath*}
and the second and the third items are
%\begin{linenomath*}
	\begin{align*}
	&\Delta t\ \m q^\mathrm{in}_\mathrm{TK}(t) \oslash \m V^{\mathrm{TK}}(t+\Delta t) \circ \m c^{\mathrm{P}}(t) = \underbrace{\Delta t  \left(\diag(\m V^{\mathrm{TK}}(t+\Delta t))\right)^{-1}  \diag(\m q^\mathrm{in}_\mathrm{TK}(t) ) }_{\textstyle \m A_{\mathrm{TK}}^\mathrm{P}(t)  }\m c^{\mathrm{P}}(t),
	\end{align*}
%\end{linenomath*}
and 
%\begin{linenomath*}
	\begin{align*}
	& \m V^\mathrm{B}_\mathrm{TK}(t+\Delta t)   \oslash \m V^{\mathrm{TK}}(t +\Delta t)  \circ \m c^{\mathrm{B}}_\mathrm{TK}(t +\Delta t) = \underbrace{\left(\diag(\m V^{\mathrm{TK}}(t +\Delta t))\right)^{-1}  \m E^\mathrm{B}_\mathrm{TK} \diag(\m V^\mathrm{B}_\mathrm{TK}(t +\Delta t) ) }_{\textstyle \m B^\mathrm{TK}(t) }\m c^{\mathrm{B}}(t +\Delta t).
	\end{align*}
%\end{linenomath*}

The fourth item $\Delta t\  \m r^\mathrm{TK}(\m c^\mathrm{TK}(t))$ can be rewritten in matrix form $\m R^\mathrm{TK}(t) \m r(\m x(t))$, and $\m R^\mathrm{TK}(t) = \Delta t \begin{bmatrix}
\m O & \m I
\end{bmatrix}$. Hence, 
%\begin{linenomath*}
	\begin{align*}
	\hspace{-0.1em}   \m c^{\mathrm{TK}}(t&+\Delta t) =  \m A_{\mathrm{TK}}^\mathrm{TK}(t) \m c^{\mathrm{TK}}(t) + \m A_{\mathrm{TK}}^\mathrm{P}(t) \m c^{\mathrm{P}}(t) +  \m B^\mathrm{TK}(t+\Delta t) \m c^{\mathrm{B}}(t+\Delta t) + \m R^\mathrm{TK}(t) \m r(\m x(t)),
	\end{align*}
%\end{linenomath*}
that is~\eqref{equ:abs-tank-mass} in Tab.~\ref{table:wqp}. Note that \textit{(i)} the $\m A_{\mathrm{TK}}^\mathrm{TK}(t)$, $\m A_{\mathrm{TK}}^\mathrm{P}(t)$, and $\m B^\mathrm{TK}(t)$ are state-space, \textit{contribution} matrices for the concentration of tanks from tanks, pipes, and boosters, \textit{(ii)} the tanks, not like junctions, always connect with pipes, that's why we can use $\m c^{\mathrm{P}}$ in~\eqref{equ:abs-tank-mass} instead of $\m c^{\mathrm{L}}$ in~\eqref{equ:abs-nodes-mass}, and \textit{(iii)} the tanks are assumed not empty, that is, inverse of $\diag(\m V^{\mathrm{TK}}(t+\Delta t))$ always exists.

\end{document}